\documentclass{gtpart}
\agtart
\usepackage{amsmath,amsthm,latexsym,amssymb}
\usepackage{eucal}
\usepackage{mathrsfs}
\usepackage[all]{xy}
\allowdisplaybreaks

\title [Loop space homology of a homotopy fiber]{An algebraic model for the loop space homology of a homotopy fiber}
\author {Kathryn Hess}
\givenname{Kathryn}
\surname{Hess}
\address{Institut de g\'eom\'etrie, alg\`ebre et topologie (IGAT) \\
    \'Ecole Polytechnique F\'ed\'erale de Lausanne \\
    CH-1015 Lausanne \\
    Switzerland}
    \email{kathryn.hess@epfl.ch}
    
\author {Ran Levi}
\givenname{Ran}
\surname{Levi}
\address {Department of mathematical sciences\\ University of Aberdeen\\ 
Meston Building 339\\ Aberdeen AB24 3UE\\ United Kingdom}
\email {ran@maths.abdn.ac.uk}
\date {\today}
 \volumenumber{}
\issuenumber{}
\publicationyear{}
\papernumber{}
\startpage{}
\endpage{}
\doi{}
\MR{}
\Zbl{}
\received{}
\revised{}
\accepted{}
\published{}
\publishedonline{}
\proposed{}
\seconded{}
\corresponding{}
\editor{}
\version{}

\numberwithin{equation}{section}

\newtheorem{thm}{Theorem}[section]
\newtheorem{lem}[thm]{Lemma}
\newtheorem{prop}[thm]{Proposition}
\newtheorem{cor}[thm]{Corollary}

\theoremstyle{definition}
\newtheorem{defn}[thm]{Definition}
\theoremstyle{remark}
\newtheorem{rmk}[thm]{Remark}
\newtheorem{ex}[thm]{Example}

\newtheorem{summary}[thm]{Summary}

\newcommand \Om{\Omega}

\newcommand \del{\partial}
\newcommand \vp{\varphi}

\newcommand\tom{\widetilde \Omega}
\newcommand \si{s\sp{-1}}
\newcommand\cotor{\operatorname {Cotor}}
\newcommand\bC{\overline C}
\newcommand\bX{\overline X}
\newcommand\cbc{C\oplus \bC}

\newcommand\bc{\bar c}
\newcommand\bx{\bar x}
\newcommand\tD{\widetilde \Delta}
\newcommand \td{\widetilde d}
\newcommand\tpsi{\widetilde\psi}
\newcommand\tPsi{\widetilde{\Psi}}

\newcommand\ttheta{\widetilde \theta}
\newcommand\tTheta{\widetilde \theta}

\newcommand\hth{\widetilde \theta}
\newcommand\hn{\nu }
 
\newcommand\shc{\mathbf F}
\newcommand\lp{\mathfrak {PL}}
\newcommand\pl{\mathfrak {PL}}
\newcommand\dl{ \mathfrak {L}_2}

\newcommand\hf{\mathfrak {LF}}

\renewcommand\H{\operatorname{H}_*}
\newcommand\op{\mathcal}
\newcommand\cat{\mathbf}
\newcommand\acoalg{\coalg{A}}
\newcommand\aalg{\op A\text{-}\cat {Alg}}
\newcommand\acirc{\underset {\op A}{\diamond}}
\newcommand\pcirc{\underset {\op P}{\diamond}}
\newcommand\ind{\operatorname{Ind}}
\newcommand\ract{\triangleleft}
\newcommand\aract{\underset {\op A}{ \ract}}
\newcommand\tC{\widetilde C}
\newcommand\twedge{\curlywedge}
\newcommand{\coalg}[1]{\op{{#1}}\text{-}\cat{Coalg}}

\newcommand{\Coalg}[2]{(\op{#1},\op{#2})\text{-}\cat{Coalg}}
\newcommand{\ob}{\operatorname{Ob}}

\begin{document}
 \begin{abstract} Let $F$ denote the homotopy fiber of a map $f:K\to L$ of $2$-reduced simplicial sets. Using as input data the strongly homotopy coalgebra structure of the chain complexes of $K$ and $L$, we construct a small, explicit chain algebra, the homology of which is isomorphic as a graded algebra to the homology of $GF$, the 
 simplicial (Kan) loop group  on $F$.  To construct this model, we develop machinery for modeling the homotopy fiber of a morphism of chain Hopf algebras.
 
 Essential to our construction is a generalization of the operadic description of the category $\cat {DCSH}$ of chain coalgebras and of strongly homotopy coalgebra maps given in \cite {HPS} to strongly homotopy morphisms of comodules over Hopf algebras.   This operadic description is expressed in terms of a  general theory of monoidal structures in categories with morphism sets parametrized by co-rings, which we elaborate here.
 \end{abstract}
 
 \maketitle


\section* {Introduction}\label{sec:introduction}

In this article we propose a ``neoclassical" approach to computing the homology algebra of double loop spaces, based on developing a deep, operadic understanding of ``strongly homotopy" structures for coalgebras and comodules, a notion that goes back  more than 30 years, to work of Gugenheim, Halperin, Munkholm and Stasheff \cite {GM}, \cite {HS}.  We also make extended use of one-sided cobar constructions, which we apply in innovative ways. 

Let $\cotor ^C(-,-)$ denote the derived functor of the cotensor product $-\underset C\square -$, for any coalgebra $C$. Eilenberg and Moore proved long ago \cite{EM} that for any (Serre) fibration $E\to B$ with fiber $F$ such that $B$ is connected and simply connected and for any commutative ring $R$, there is an  $R$-linear isomorphism
\begin{equation}\label{eqn:moore-iso1}
\H(F;R)\cong \cotor ^{C_{*}(B;R)}\big(C_{*}(E;R),R\big).
\end{equation}
In particular, if $B$ is actually $2$-connected, then 
\begin{equation}\label{eqn:moore-iso2}
\H(\Om F;R)\cong \cotor ^{C_{*}(\Om B;R)}\big(C_{*}(\Om E;R),R\big),
\end{equation}
so that for any $2$-connected space $X$,  there is a linear isomorphism
\begin{equation}\label{eqn:moore-iso3}
\H(\Om ^2 X;R)\cong \cotor ^{C_{*}(\Om X;R)}(R,R).
\end{equation}
The duals of Theorems 4.1 and  5.1 in \cite{fht} imply that  $\cotor ^{C_{*}(\Om B;R)}\big(C_{*}(\Om E;R),R\big)$ admits an algebra structure with respect to which the isomorphism (\ref{eqn:moore-iso2}) can be assumed to be an algebra morphism.  

In this article we define a functor $\hf$ that associates to any map $f$ of $2$-reduced simplicial sets  a chain algebra $\hf (f)$ such that $\H \big(\hf (f)\big)$ is isomorphic as an algebra to $\H (GF)$, where $F$ is the homotopy fiber of $f$ and  $G$ is the Kan loop group functor on simplicial sets (Theorem \ref{thm:loopfib-model}).  The algebra isomorphism  $\H \big(\hf (f)\big)\cong\H (GF)$ is realized on the chain level by a zig-zag of quasi-isomorphisms of chain algebras
$$\hf(f)\xleftarrow\simeq \bullet\xrightarrow\simeq \cdots\xleftarrow\simeq\bullet\xrightarrow\simeq C(GF).$$ 
As a special case, we obtain a functor $\dl$ from the category of $2$-reduced simplicial sets to the category of connected chain algebras over a principal ideal domain $R$ such that $\H \big(\dl (K)\big)$ is isomorphic as a graded algebra to $\H (G^2K;R)$. 

The model that we propose for the loop homology of a homotopy fiber offers certain advantages.  First, there are no extension problems to be solved:  the homology algebra of the model is exactly isomorphic to the homology algebra of the loops on the homotopy fiber.  Second,  our model is functorial, so that it can be applied to determining the homomorphism induced on double loop space homology by a simplicial map.

 Finally, our model is ``small''  and therefore amenable to explicit computations.   More precisely, if $K$ is a simplicial set with exactly $n$ nondegenerate simplices of positive degree, where $n<\infty$, then our model $\dl(K)$ of $G^2K$ is a subalgebra of a free algebra on $2n$ generators.  Like the differential in the cobar construction on $C(K)$, the differential in $\dl(K)$ depends only on the differential and the comultiplication on $C(K)$.    In particular, in section \ref{subsec:formal}, we provide an explicit, relatively simple description of our model when $K$ is either formal or a double suspension.

By way of comparison, note that the iterated cobar construction on the chains on $K$, which is another model of $G^2K$, is free as an algebra on an infinite number of generators.  Its differential depends not only on the differential and the comultiplication on $C(K)$, but also on the natural comultiplication on the cobar construction on $C(K)$, which can be very involved.    Another possible model, the cobar construction on $C(GK)$, is also free, but on a generating set that is infinite in each degree, and, in addition, has a very complicated differential.  Finally, both the multiplication and the differential on the chain Hopf algebra $C(G^2K)$ itself are extremely complex.

In the general case of loops on the homotopy fiber $F$ of a simplicial map $K\to L$, the dual of Theorem 5.1 in \cite {fht} states that there is a quasi-isomorphism of chain algebras 
$$C(GF)\xrightarrow\simeq C(GK)\otimes _{t_{\Om}}\Om C(GL),$$
where $C(GK)\otimes _{t_{\Om}}\Om C(GL)$ denotes the one-sided cobar construction of Definition \ref {defn:one-sided}, endowed with the multiplication of Corollary \ref{cor:cobar-alg}.   The chain algebra  $C(GK)\otimes _{t_{\Om}}\Om C(GL)$ is not of finite type, even if $K$ and $L$ have only a finite number of nondegenerate simplices, and both its differential and its multiplication are quite complicated.  On the other hand, if $K$ and $L$ have exactly $m$ and $n$ nondegerate simplices of positive degree, respectively, then the chain algebra model constructed here, $\hf(f)$, is a subalgebra of a free algebra on $m+2n$ generators, so that its multiplicative structure is relatively simple.  Its differential is also much easier to give explicitly than that of $C(GK)\otimes _{t_{\Om}}\Om C(GL)$. 

To construct our models, we need the full \emph{Alexander-Whitney coalgebra} structure of the normalized chains on a simpicial set.  The category $\cat F$ of Alexander-Whitney coalgebras (cf., Definition \ref{defn:F}) was introduced and studied in \cite {HPST}.  The objects of $\cat F$ are connected chain coalgebras $(C,\Delta)$  such that  the comultiplication $\Delta$ is itself a coalgebra map \emph{up to strong homotopy}, i.e., up to a coherent, infinite family of homotopies, which we denote $\Psi$.  Furthermore,  there is a functor $\tom:\cat F\to \cat H$, where $\cat H$ is the category of chain Hopf algebras, such that the chain algebra underlying $\tom (C, \Psi)$ is $\Om C$, the cobar construction on $C$.

An \emph{Alexander-Whitney model} of a chain Hopf algebra $H$ consists of an Alexander-Whitney coalgebra $(C,\Psi)$ together with a quasi-isomorphism $\theta:\tom (C,\Psi)\xrightarrow\simeq H$ of chain algebras that is also a map of coalgebras up to strong homotopy, where the homotopies are appropriately compatible with the multiplicative structure (cf., Definition \ref{defn:awmodel}).  As illustrated by the results in this article, Alexander-Whitney models can be useful tools for homology calculations in $\cat H$.

The topologist's motivation for considering the category $\cat F$ is the existence of a natural Alexander-Whitney model of the chain Hopf algebra $C(GK)$, where $K$ is a reduced simplicial set.  As shown in \cite {HPST}, there is a functor $\widetilde C:\cat {sSet}_{0}\to \cat F$ from the category of reduced simplicial sets to the category of Alexander-Whitney coalgebras such that for all simplicial sets $K$, there is a natural quasi-isomorphism of chain algebras
$$\theta_{K}:\tom\widetilde C(K)\xrightarrow\simeq C(GK),$$
which is also a map of coalgebras up to strong homotopy. 

Given the existence of natural Alexander-Whitney models, the most important steps on the path to constructing the model $\hf(f)$ and to proving that its homology algebra is isomorphic to $\H(GF)$ are the following.
 
\begin{enumerate}
\item For any chain Hopf algebra $H$ and any right $H$-comodule algebra $B$, we observe that  $\cotor ^H(B,R)$  admits a natural graded algebra structure (Corollary \ref{cor:cotor-alg} and, more generally, Proposition \ref{prop:mult-cotor}).  In particular, for any morphism  $\varphi: H'\to H$ of chain Hopf algebras, $\cotor ^H(H',R)$, which can be seen as the homology of the ``homotopy fiber'' of $\varphi$, admits a natural graded algebra structure.
\item We show that the category  $\cat F$ admits a natural  ``based-path'' construction, i.e., a functor $\widetilde {\mathfrak P}:\cat F\to \cat F$ such that $\widetilde{\mathfrak P}(C,\Psi)$ is acyclic for all $(C,\Psi)$, together with a natural ``projection'' morphism in $\cat F$ from $\widetilde{\mathfrak P}(C,\Psi)$ to $(C,\Psi)$ (Definition \ref{defn:based-path}).
\item For any morphism $\omega: (C',\Psi')\to (C,\Psi)$ in $\cat F$, we prove that  the chain Hopf algebra
$$\tom\bigg((C',\Psi')\coprod \widetilde{\mathfrak P}(C,\Psi)\bigg)$$ 
is cofree over $\tom (C,\Psi)$ on a cobasis that is itself a sub chain algebra, denoted $\hf(\omega)$, of $\tom\big((C',\Psi')\coprod \widetilde{\mathfrak P}(C,\Psi)\big)$ (Corollary \ref{cor:pl-cofree}).
\item Given an Alexander-Whitney model $\omega: (C',\Psi')\to (C,\Psi)$ of a morphism of chain Hopf algebras $\varphi:H'\to H$, we prove that   $\H \big(\hf(\omega)\big)\cong\cotor^H(H',R)$ as graded algebras (Theorem \ref{thm:homol-loopfib}).
\end{enumerate}

Let $f:K\to L$ be a simplicial morphism of $2$-reduced simplicial sets with homotopy fiber $F$. Applying (4) to  $\widetilde C(f):\widetilde C(K)\to\widetilde C(L)$, we obtain an isomorphism of algebras
$$\H (GF)\cong \H \bigg(\hf\big(\widetilde C(f)\big)\bigg),$$
thanks to the algebra isomorphism
$$\H (GF)\cong\cotor ^{C(GL)}\big(C(GK), R\big)$$
that follows from the dual of Theorem 5.1 in \cite{fht}. 

To make this article as self-contained as possible and to establish our notation, we begin in section \ref{sec:prelims} by recalling the rather extensive foundations on which our current research is built.  Section \ref{subsec:dcsh}, in which we describe the operadic approach  to strongly-homotopy coalgebra structures of \cite{HPS}, is particularly important for the later sections of this paper and essential to providing a clean description of the yoga of Alexander-Whitney coalgebras.  Readers unfamiliar with the role of co-rings in monoidal categories as parametrizing objects for enlarged morphism sets or with operads will find all of the necessary definitions in sections \ref{subsec:cats} and \ref{subsec:operads}.

Section \ref{sec:ext-mult-nat} concerns the naturality of multiplicative structure on 
$\cotor$, which plays an important role in the proofs of Theorems \ref{thm:lift-props} and \ref{thm:homol-loopfib}, the key elements of step (4) in the plan outlined above.  Given chain Hopf algebras $H$ and $H'$, as well as a right $H$-comodule algebra $B$ and a right $H'$-comodule algebra $B'$, there is an obvious notion of ``morphisms'' from $(H;B)$ to $(H';B')$:  the set of pairs $(f;g)$, where $f:H\to H'$ is a chain Hopf algebra map and $g:B\to B'$ is a chain algebra map respecting the coactions of $H$ and $H'$.  It is easy to see that any such pair induces an algebra map $\cotor^H(B,R)\to \cotor^{H'}(B',R)$.  There is however a more general type of ``morphism'' from $(H;B)$ to $(H';B')$, which we call a \emph{comodule algebra map up to strong homotopy (CASH map)}, that also induces an algebra map $\cotor^H(B,R)\to \cotor^{H'}(B',R)$.  

In section \ref{subsec:comodules} we define CASH maps and establish existence results (Propositions \ref{prop:rel-free-ext} and \ref{prop:cash-coproduct}) that we use afterwards to prove Theorems \ref{thm:lift-props} and \ref{thm:homol-loopfib}.  Before verifying the existence results, we provide an equivalent, operadic definition of CASH maps, modeled on the operadic approach to strongly homotopy coalgebra structures, that facilitates considerably the bookkeeping involved in working with the infinite family of homotopies associated to a CASH map. The general study of monoidal structures and parametrizations by co-rings in section  \ref{sec:mon-struct} is essential to the development of this operadic approach.

Sections \ref{sec:path-loop} and \ref{sec:homol-htpyfib} are  devoted to the study of Alexander-Whitney coalgebras and their use in calculations of the homology of homotopy fibers in the category $\cat H$ of chain Hopf algebras. Topology comes into play again in section \ref{sec:topology}, where we apply the purely algebraic results of the preceeding sections to constructing our loop-homotopy fiber model. In particular, we study the special cases of  double suspensions and  of formal spaces,  obtaining a simplified model for the homology of their double loop spaces, which is a free algebra on a set of generators we describe completely.

The first author would like to thank the University of Aberdeen for its kind hospitality during the initial phase of research on this project, while the second author would like to thank the EPFL for hosting him during the completion of the project.  Both authors would like to thank 
referees for pointing out the relevance of \cite {fht} to their work and for providing helpful organizational advice.

\subsection*{Notation and conventions}
\begin{itemize}
\item Given objects $A$ and $B$ of a category $\cat C$, we let $\cat C(A,B)$ denote the set of morphisms with source $A$ and target $B$. 
\item Throughout this paper we are working over a principal ideal domain $R$.  We denote the category of graded $R$-modules by $\cat{grMod}_R$ and the category of chain complexes over $R$ by $\cat{Ch}_R$.   The underlying graded modules of all chain (co)algebras are assumed to be $R$-free.  
\item The degree of an element $v$ of a graded module $V$ is denoted either $|v|$ or simply $v$, when used as an exponent, and no confusion can arise. 
\item Throughout this article we apply the Koszul sign convention for commuting elements  of a graded module or for commuting a morphism of graded modules past an element of the source module.  For example,  if $V$ and $W$ are graded algebras and $v\otimes w, v'\otimes w'\in V\otimes W$, then $$(v\otimes w)\cdot (v'\otimes w')=(-1)^{|w|\cdot |v'|}vv'\otimes ww'.$$ Futhermore, if $f:V\to V'$ and $g:W\to W'$ are morphisms of graded modules, then for all $v\otimes w\in V\otimes W$, 
$$(f\otimes g)(v\otimes w)=(-1)^{|g|\cdot |v|} f(v)\otimes g(w).$$
\item A graded module $V$ is \emph {bounded below} if there is some $N\in \mathbb Z$ such that $V_k=0$ for all $k<N$. It is \emph { $n$-connected} if, in 
particular, $V_{k}=0$ for all $k\leq n$.  We write $V_{+}$ for $V_{>0}$.
\item The \emph {suspension} endofunctor $s$ on the category of graded modules is defined on objects $V=\bigoplus _{i\in \mathbb Z} V_ i$ by
$(sV)_ i \cong V_ {i-1}$.  Given a homogeneous element $v$ in
$V$, we write $sv$ for the corresponding element of $sV$. The suspension $s$ admits an obvious inverse, which we denote $\si$.
\item Given chain complexes $(V,d)$ and $(W,d)$, the notation
$f:(V,d)\xrightarrow{\simeq}(W,d)$ indicates that $f$ induces an isomorphism in homology. 
In this case we refer to $f$ as a \emph { quasi-isomorphism}.
\item Let $V$ be a positively-graded $R$-module.  The free associative algebra on $V$ is denoted 
$TV$, i.e., 
$$TV\cong R\oplus V\oplus (V\otimes V)\oplus (V\otimes V\otimes V)\oplus\cdots .$$
A typical basis element of $TV$ is denoted $v_ 1\cdots v_ n$.
\item Given a comodule 
$(M,\nu )$ over a coalgebra $(C,\Delta)$, we let $\Delta ^{(i-1)}$ denote the 
iterated comultiplication $C\to  C^{\otimes i}$ and $\nu ^{(i)}$ the 
iterated coaction $M\to M\otimes C^{\otimes i}$.   The reduced comultiplication is denoted $\overline \Delta$.
\item If $C$ is a simply connected chain 
coalgebra with reduced comultiplication $\overline \Delta$ and differential $d$, then $\Om C$ denotes  the \emph {cobar construction} on $C$, i.e., the chain algebra $(T\si (C_{+}),d_{\Om})$, where 
$d_{\Om}=-\si ds+(\si \otimes \si)\overline \Delta s$ on generators.

Furthermore, for every pair of simply-connected chain coalgebras $C$ and $C'$ 
\begin{equation}\label{eqn:milgram}
q:\Om \big( C\otimes C'\big )\xrightarrow\simeq \Om C\otimes \Om C'
\end{equation}
denotes the quasi-isomorphism of chain algebras defined by Milgram in \cite{Mi}.

\end{itemize}

\tableofcontents

\section {Preliminaries}\label{sec:prelims}

For the convenience of the reader, we recall here certain algebraic foundations of our work. We begin by reminding the reader how co-rings in monoidal categories can act as parametrizing objects for categories of modules with enlarged morphism sets, as described in \cite{HPS} and \cite {HPST}.  We then review the theory of operads, seen as monoids  with respect to a certain nonsymmetric monoidal structure on the category of symmetric sequences of objects in a given symmetric monoidal category.  In particular we analyze the category of right modules over a given operad $\op P$, comparing it to  the category of $\op P$-coalgebras.  Finally, we summarize briefly results in \cite{HPS} and \cite {HPST} that provide an operadic description of  the category $\cat {DCSH}$ of chain coalgebras and of strongly homotopy coalgebra maps, in terms of a certain co-ring over the associative operad.

\subsection {Co-rings in monoidal categories}\label{subsec:cats}

Let $(\cat C, \otimes, I)$ be a monoidal category, and let $(A,\mu, \eta)$ be a monoid in $\cat C$.  
If the category $\cat C$ admits coequalizers and $M\otimes -$ and  $-\otimes N$ preserve colimits for all objects $M$ and $N$, then  the category of $A$-bimodules, ${}_{A}\cat {Mod}_{A}$, is a monoidal category also, with monoidal product $-\underset A\otimes -$.  If $M$ and $M'$ are $A$-bimodules, then
$M\underset A\otimes M'$ is the coequalizer of the diagram 
$$M\otimes A\otimes M'\underset {\rho\otimes Id_{M'}}{\overset {Id_{M}\otimes \lambda '}{\rightrightarrows}}M\otimes M'.$$
The unit object with respect to $-\underset A\otimes -$ is $A$ itself, where the right and left $A$-actions on $A$ are given by the multiplication map $\mu$.

\begin{defn}\label{defn:co-ring} An \emph{$A$-co-ring} is a comonoid in the monoidal category $({}_{A}\cat {Mod}_{A},\underset A\otimes,A) $.  An $A$-co-ring thus consists of an $A$-bimodule $M$, together with two morphisms of $A$-bimodules
$$\psi: M\to M\underset A\otimes M\qquad\text{and}\qquad \varepsilon: M\to A$$
such that $\psi$ is coassociative and counital with respect to $\varepsilon$.
\end{defn}

Examples of co-rings abound in algebra and topology.  In particular, any Frobenius algebra is a co-ring over itself, while the Hopf algebroids of stable homotopy theory are co-rings with extra structure.  Moreover, any ring homomorphism $\vp: A\to B$ induces a canonical $B$-co-ring structure on $M=B\otimes _{A}B$, where the comultiplication is
$$M\to M\otimes _{B}M: b\otimes b'\mapsto (b\otimes 1)\otimes (1\otimes b').$$

Co-rings play an important role in this article, as they induce natural enlargements of categories of modules, leaving the objects fixed and expanding the morphism sets.    Allowing larger morphism sets translates into weakening the notion of morphism of modules.
In this sense a co-ring plays the role of a family of parameters, with respect to which such a weaker notion is coherently defined.

\begin{defn} \label{defn:fat} Let $\cat {Mod}_{A}$ denote the category of right $A$-modules. Given an $A$-co-ring $(M,\psi, \varepsilon)$, let $\cat{Mod}_{A,M}$ denote the category with $\ob\cat {Mod}_{A,M}=\ob\cat {Mod}_{A}$ and 
$$\cat {Mod}_{A,M}(N,N'):=\cat {Mod}_{A}(N\underset A\otimes M, N').$$
Given $f\in\cat {Mod}_{A,M}(N,N')$ and $f'\in  \cat {Mod}_{A,M}(N',N'')$, their composite $f'f\in \cat {Mod}_{A,M}(N,N'')$ is equal to the composite in $\cat{Mod}_{A}$ of the following sequence of morphisms of right $A$-modules.
$$N\underset A\otimes M\xrightarrow{Id_{N}\underset A\otimes \psi} N\underset A\otimes M\underset A\otimes M\xrightarrow {f\underset A\otimes Id_{M}} N'\underset A\otimes M\xrightarrow {f'} N''$$
\end{defn}

Composition in $\cat {Mod}_{A,M}$ is associative and unital, since $\psi$ is coassociative and counital.  Furthermore, there is a natural, faithful functor
\begin{equation}\label{eqn:emb-fat}
\mathfrak I_{M}:\cat {Mod}_{A}\to \cat {Mod}_{A,M},
\end{equation}
which is the identity on objects and which sends a morphism $f:N\to N'$ of right $A$-modules to 
$$\mathfrak I_{M}(f)=f\underset A\otimes \varepsilon: N\underset A\otimes M\to N'\underset A\otimes A\cong N'.$$
The category $\cat {Mod}_{A,M}$ is therefore truly an enlargement of $\cat{Mod}_{A}$.

We conclude this section by clarifying our vision of a co-ring as a family of a parameters. 

\begin{defn}\label{defn:parameter}  Let $(M, \psi, \varepsilon)$ be an $A$-co-ring, endowed with a strict morphism of left $A$-modules $\eta :A\to M$.  Let $N,N'\in \ob \cat {Mod}_{A}$.  A morphism $f\in \cat C(N,N')$ is a \emph{morphism of right $A$-modules up to $M$-parametrization} if there is a strict morphism of right $A$-modules $g:N\underset A\otimes M\to N'$ such that the following diagram in $\cat C$ commutes.
$$\xymatrix{N\cong N\underset A\otimes A\ar[drr]^f\ar[d]_{Id_{N}\underset A\otimes \eta}\\
N\underset A\otimes M\ar[rr]^g &&N'
}$$
\end{defn}  

There is an analogous enlargement of the category of left $A$-modules.  For the experts, we note that these enlargements are, of course,  coKleisli constructions, induced by the comonads $-\underset A\otimes M$ and $M\underset A\otimes -$.

\subsection {Operads and their modules and coalgebras}\label{subsec:operads}

Let $(\cat C, \otimes , I)$ be a symmetric monoidal category such that $\cat C$ admits coequalizers and countable coproducts and has an initial object $0$.  Let $\mathbf C^\Sigma $ denote the \emph{category of symmetric sequences} in $\mathbf C$.  An object $\mathcal X$ of $\mathbf C^\Sigma$ is  a family $\{\mathcal X(n)\in \mathbf C\mid n\geq 0\}$ of objects in $\mathbf C$ such that $\mathcal X(n)$ admits a right action of the symmetric group $\Sigma_n$, for all $n$.  The object $\op X(n)$ is called the \emph{$n^{\text{th}}$ level} of the symmetric sequence $\op X$.
 
For all $\op X, \op Y\in \cat C^\Sigma$, a \emph{morphism of symmetric sequences} $\varphi:\op X\to \op Y$ consists of a family 
$$\{\varphi_{n}\in \cat C\big(\op X(n), \op Y(n)\big)\mid \varphi_{n} \quad\text{is $\Sigma_{n}$-equivariant}, n\geq 0\}.$$
More formally, $\cat C^\Sigma$ is the category of contravariant functors from the symmetric groupoid $\boldsymbol \Sigma$ to $\cat C$, where $\ob \boldsymbol \Sigma =\mathbb N$, the set of natural numbers, and $\boldsymbol\Sigma (m,n)$ is empty if $m\ne n$, while $\boldsymbol \Sigma (n,n)=\Sigma_{n}$.

The category $\cat C$ can be ``linearly'' embedded in the category $\cat C^\Sigma$, via a functor
\begin{equation}\label{eqn:emb-lin}
 \op L:\cat C\to \cat C^\Sigma,
 \end{equation}
which is defined on $A\in \ob \cat C$ by $\op L(A)(1)=A$ and $\op L(A)(n)=0$ for all $n\not=1$ and similarly for morphisms.

There is another important embedding of $\cat C$ into $\cat C^\Sigma$
\begin{equation}\label{eqn:emb-t}
\mathcal T: \mathbf C\to\mathbf C^\Sigma
\end{equation} 
defined by $\mathcal T(A)(n)=A^{\otimes n}$ for all $n$. The right action of $\Sigma _{n}$ on $\op T(A)(n)=A^{\otimes n}$ is given by permutation of the factors, using iterates of the natural symmetry isomorphism $\tau:A\otimes A\xrightarrow\cong A\otimes A$ in $\cat C$.   
For example, if $\cat C$ is the category of graded modules, then
 $$(a_{1}\otimes \cdots \otimes a_{n})\cdot \sigma =a_{\sigma (1)}\otimes \cdots \otimes a_{\sigma (n)}$$  for all $a_{1},..., a_{n}\in A$.

As a first indication of the role of differential structure in symmetric sequences, we introduce the following useful operation on symmetric sequences of chain complexes in the image of $\op T$.  The analogy with the notion of a derivation on an algebra is evident.  

\begin{defn}\label{defn:derivation} Let $f,g,t:A\to B$ be morphisms of graded $R$-modules, homogeneous of degrees $0$, $0$ and $m$, respectively. The \emph {$(f,g)$-derivation of  symmetric sequences} induced by $t$ is the morphism of symmetric sequences 
$$\op D_{(f,g)}(t):\op T(A)\to \op T(B)$$ 
that is of degree $m$ in each level and that is defined as follows in level $n$.
$$\op D_{(f,g)}(t)_{n}=\sum _{j=0}^{n-1} f^{\otimes j}\otimes t\otimes g^{\otimes n-j-1}$$
When $A=B$ and $f=Id_A=g$, we simplify notation and write $\op D(t)$ for the $(Id_A, Id_A)$-derivation induced by $t$.\end{defn}

\begin{ex}If $C$ is chain complex with differential $d$, then the levelwise differential on $\op T(C)$ is $\op D(d)$.
\end{ex}

In this article we use the following two monoidal structures on the category of symmetric sequences.

\begin{defn}\label{defn:level-tensor}
The \emph{level tensor product} of two symmetric sequences
$\op{X}$ and $\op{Y}$ is the symmetric sequence given by
\[
    ( \op{X} \otimes \op{Y} )(n) = \op{X}(n) \otimes \op{Y}(n)
    \quad (n \geq 0),
\]
endowed with the diagonal action of $\Sigma_{n}$.
\end{defn}

The following, well-known result is very easy to prove.

\begin{prop}\label{prop:level-tensor}
Let $\op{C} = \{ \op{C}(n) \}_{n \geq 0}$ be the symmetric
sequence with $\op{C}(n) = I$ and trivial $\Sigma_{n}$-action,
for all $n \geq 0$.  Then $(\cat{C}^{\Sigma},\otimes,\op{C})$ is a
closed symmetric monoidal category, called the \emph{level
monoidal structure} on $\cat{C}^{\Sigma}$.
\end{prop} 

A (co)monoid in $\cat C^\Sigma$ with respect to the level monoidal structure is called a \emph{level (co)monoid}.

Note that the functor $\mathcal T$ is strong monoidal with respect to the level monoidal structure on symmetric sequences, i.e., for all $C,C\in \ob\cat C$, there is a natural isomorphism $\op T(C\otimes C')\cong \op T(C)\otimes \op T(C')$, given in each level by iterated application of the natural symmetry isomorphism in $\cat C$.

 The category $\mathbf C^\Sigma$ also admits a nonsymmetric, right-closed monoidal structure, defined as follows.
 
 \begin{defn}\label{defn:comp-prod} The \emph{composition tensor product} of two symmetric sequences $\op X$ and $\op Y$ is the symmetric sequence $\op X\diamond \op Y$ given by
$$(\op X\diamond\op Y)(n)=\coprod_{\substack{ k\geq 1\\ \vec n\in I_{k,n}}} \op X(k)\underset {\Sigma _k}{\otimes} \big(\op Y(n_1)\otimes \cdots\otimes \op Y(n_k)\big)\underset {\Sigma _{\vec n}}{\otimes} I[\Sigma _n],$$
where $I_{k,n}=\{\vec\imath=(n_1,...,n_k)\in \mathbb N^k\mid \sum _j n_j=n\}$ and $\Sigma_{\vec n}=\Sigma _{n_1}\times \cdots\times \Sigma _{n_k}$, seen as a subgroup of $\Sigma _n$.
The left action of $\Sigma _k$ on  $\coprod _{\vec n\in I_{k,n}}\op Y(n_{1})\otimes \cdots \otimes \op Y(n_{k})$ is given by permutation of the factors, using the natural symmetry isomorphism $A\otimes B\cong B\otimes A$ in $\cat C$. 
 \end{defn}

 \begin{prop}  Let $\op J$ denote the symmetric sequence with $\mathcal J(1)=1$ and $\mathcal J(n)=0$ otherwise, with trivial $\Sigma _{n}$-action.  Then $(\mathbf C^\Sigma,\diamond , \mathcal J)$ is a right-closed monoidal category, called the \emph{composition monoidal structure} on $\cat C^\Sigma$.
 \end{prop} 
 
A proof of this result can be found  in~\cite[section
II.1.8]{MSS}.
 
Unwrapping the definition of the composition product of symmetric sequences, we obtain the next, well-known lemma, which tells us which data determine a morphism with source a composition of symmetric sequences.  
 
\begin{lem}\cite{Markl}\label{lem:classical} Let $\op X$, $\op Y$ and $\op Z$ be symmetric sequences in $\cat C$.  Let $I_{m,n}= \{\vec n=(n_{1},...,n_{m})\mid \sum _{j}n_{j}=n\}$.  Let 
$$\mathfrak F=\big\{ \op X(m)\otimes \big(\op Y(n_{1})\otimes \cdots \otimes \op Y(n_{m})\big)\xrightarrow{ \theta _{\vec n}} \op Z(n)\mid \; n\geq 0, m\geq 1, \vec n\in I_{m,n}\big \}$$
be a family of morphisms in $\cat C$.  If the following diagrams commute for all $m$, $n$, $\vec n$, $\sigma \in \Sigma _{m}$ and $\tau _{j}\in \Sigma_{n_{j}}$ for $1\leq j\leq m$, then $\mathfrak F$ induces a morphism of symmetric sequences $\theta: \op X\diamond \op Y\to\op Z$.
$$\xymatrix{
 \op X(m)\otimes \big(\op Y(n_{1})\otimes \cdots \otimes \op Y(n_{m})\big)\ar [d]_{\theta _{\vec n}}\ar [r]^{\sigma \otimes \sigma ^{-1}}& \op X(m)\otimes \big(\op Y(n_{\sigma (1)})\otimes \cdots \otimes \op Y(n_{\sigma (m)})\big)\ar[d] _{\theta _{\sigma^{-1} \vec n}}\\
 \op Z(n)\ar [r]^{\sigma (n_{\sigma (1)},...,n_{\sigma (m)})}&\op Z(n)}$$
 $$\xymatrix{
 \op X(m)\otimes \big(\op Y(n_{1})\otimes \cdots \otimes \op Y(n_{m})\big)\ar [d]_{\theta _{\vec n}}\ar [rr]^{1\otimes \tau _{1}\otimes \cdots\otimes \tau _{m}}&& \op X(m)\otimes \big(\op Y(n_{1})\otimes \cdots \otimes \op Y(n_{m})\big)\ar[d] _{\theta _{\vec n}}\\
 \op Z(n)\ar [rr]^{\tau _{1}\oplus \cdots \oplus \tau _{m}}&&\op Z(n)}$$
 \end{lem}
 
 In the statement above, $\sigma ^{-1} \vec n:=(n_{\sigma (1)},...,n_{\sigma (m)})$, which defines a left action of $\Sigma _{m}$ on $I_{m,n}$.

\begin{rmk}\label{rmk:intertwine} For any objects $\mathcal X, \mathcal X', \mathcal Y, \mathcal Y'$ in $\mathbf C^\Sigma$, there is an obvious, natural \emph{intertwining} map
$$\xymatrix{\mathfrak i: (\mathcal X\otimes \mathcal X')\diamond (\mathcal Y\otimes \mathcal Y')\ar [r]&(\mathcal X\diamond \mathcal Y)\otimes (\mathcal X'\diamond\mathcal Y')}.$$
\end{rmk}

\begin{defn}An \emph {operad} in $\mathbf C$ is a monoid with respect to the composition product, i.e., a triple $(\op P, \gamma,\eta )$, where $\gamma: \op P\diamond \op P\to \op P$ and $\eta:\op J\to \op P$ are morphisms in $\mathbf C^\Sigma$, and $\gamma$ is appropriately associative and unital with respect to $\eta$.   A morphism of operads is a monoid morphism in the category of symmetric sequences.
\end{defn} 

The most important example of an operad in this paper is the \emph {associative operad} $\mathcal A$, given by $\mathcal A(n)=I[\Sigma _n]$ for all $n$, endowed with the obvious multiplication map, induced by permutation of blocks. 

Operads derive their importance from their role in parametrizing $n$-ary (co)operations and governing the identites among them.  In this article we focus on cooperations and thus on coalgebras over an operad $\op P$.  A \emph {$\op{P}$-coalgebra} is an object $A$ of $\cat{C}$
along with a sequence of structure morphisms
\[
    \theta_{n} : A \otimes \op{P}(n) \rightarrow A^{\otimes n},\quad n\geq 0
\]
that are appropriately associative, equivariant, and unital.  We refer the reader to e.g., \cite{MSS}, for the full definition. 

A morphism of $\op P$-coalgebras is a morphism in $\cat C$ that commutes with the coalgebra structure maps. The category of
$\op{P}$-coalgebras and their morphisms is
denoted $\coalg{P}$.

\begin{rmk}  Algebraists are used to thinking of coalgebras as modules with additional structure.  It is important to note that if  $\op P$ is an operad, then a $\op P$-module (in the sense defined in section \ref{subsec:cats}) is a object of $\cat C^\Sigma$ with additional structure, while a $\op P$-coalgebra is an object of $\cat C$ with additional structure.

On the other hand, as explained in section 2.2 of \cite{HPS}, the functor $\mathcal T$ restricts to a faithful functor 
 $$\xymatrix@1{\mathcal T: \coalg{P}\ar [r]&\mathbf{Mod}_{\mathcal P}}$$
 from the category of $\mathcal P$-coalgebras to the category of right $\mathcal P$-modules (with respect to the composition product $\diamond$), i.e,  $\op P$-coalgebra structure on an object $A$  in $\cat C$ induces a right $\op P$-action map $\op T(A)\diamond \op P\to \op T(A)$ in $\cat C^\Sigma$.
 \end{rmk}
 
 Let $(\op M, \psi, \varepsilon)$ be a $\op P$-co-ring, and consider $\cat {Mod}_{\op P, \op M}$, the enlarged version of $\cat{Mod}_{\op P}$ described in the section \ref{subsec:cats}.  Define an enlarged version $\Coalg PM$ of $\coalg P$ by
 $\ob \Coalg PM=\ob \coalg P$
\begin{equation}\label{eqn:fatcoalg}
\Coalg PM(A, A'):=\cat{Mod}_{\op P,\op M}\big (\op T(A),\op T(A')\big)=\cat {Mod}_{\op P}\big(\op T(A)\pcirc \op M, \op T(A')\big),
\end{equation}
 for all $A,A'\in \ob \Coalg PM$, with composition defined as in $\cat {Mod}_{\op P,\op M}$.  
 
 Let $A$ and $A'$ be $\op P$-coalgebras. In keeping with Definition \ref{defn:parameter}, we say that a morphism $f:A\to A'$ is a \emph{morphism of $\op P$-coalgebras up to $\op M$-parametrization} if $\op T(f):\op T(A)\to \op T(A')$ is a morphism of right $\op P$-modules up to $\op  M$-parametrization.
 
 From this formulation, it follows that co-rings over operads are, in a strong sense, \emph {relative operads}.  They parametrize higher, ``up to homotopy'' structure on morphisms of $\op P$-coalgebras and govern relations among the higher homotopies and the $n$-ary cooperations on the source and target.

\subsection {Strongly-homotopy coalgebra structures}\label{subsec:dcsh}

The category $\mathbf {DCSH}$ of coassociative chain coalgebras and of coalgebra morphisms up to strong homotopy was first defined by Gugenheim and Munkholm in the early 1970's \cite {GM}, when they were studying extended naturality of the functor $\operatorname{Cotor}$.   Its objects are simply connected, augmented, coassociative chain coalgebras, and a morphism from $C$ to $C'$ is a map of chain algebras $\Omega C\to \Omega C'$.  The category $\mathbf {DCSH}$ plays an important role in topology (cf., Theorem \ref{thm:loopmodel}).

In a slight abuse of terminology, we say that a chain map between chain coalgebras $f:C\to C'$ is a \emph {DCSH map} if there is a morphism  in $\mathbf {DCSH}(C,C')$ of which $f$ is the linear part. In other words, there is a map of chain algebras $g:\Om C\to \Om C'$ such that 
$$g|_{s^{-1}C_+}=\si f s +\quad\text{higher-order terms}.$$

Let $\op A$ denote the associative operad in the category of chain complexes. In \cite {HPS} the authors constructed an $\mathcal A$-co-ring $\mathcal F$, called the \emph {Alexander-Whitney co-ring}, which can applied in the framework of section \ref{subsec:cats} to providing an operadic description of $\mathbf{DCSH}$.  The co-ring $\op F$ also admits a level comultiplication $\Delta_{\op F}:\op F\to \op F\otimes \op F$ that is compatible with its composition comultiplication  $\psi_{\op F}: \op F\to \op F\acirc \op F$ and that plays an important role in development of monoidal structure in $\cat {DCSH}$ (cf., section \ref{subsec:alex-whit}).

The symmetric sequence of graded modules underlying $\mathcal F$ is $\op A\diamond \op S\diamond \op A$, where, for all $n\geq 1$, $\op S(n)= R[\Sigma _n]\cdot z_{n-1}$, the free $R[\Sigma _n]$-module on a generator of degree $n-1$, and $\op S(0)=0$.   We refer the reader to pages 853 and 854 in  \cite{HPST} for the explicit formulas for the differential $\del_{\op F}:\op F\to \op F$, the composition comultiplication $\psi _{\op F}$  and the level comultiplication $\Delta_{\op F}$.  We remark that $\op F$ admits a natural filtration with respect to which both $\psi_{\op F}$ and $\Delta _{\op F}$ are filtration-preserving, while $\del_{\op F}$ is filtration-decreasing.

Consider $\Coalg AF$ (cf., equation \ref{eqn:fatcoalg}). Any morphism $\theta\in \Coalg AF(C,C')$ gives rise to a family $\mathfrak F(\theta)$ of linear maps from $C$ into $\op T(C')$, defined as follows.
\begin{equation}\label{eqn:family}
\mathfrak F(\theta):=\{\theta _k=\theta (-\otimes z_{k-1}):C\to  (C')^{\otimes k}=\op T(C')(k)\mid k\geq 1\}
\end{equation}
The existence of such a family $\mathfrak F(\theta)$ is equivalent to the existence of a morphism of symmetric sequences of graded modules $\op L(C)\diamond \op S\to \op T(C')$, where $\op L:\cat {grMod}_{R}\to \cat {grMod}_{R}^\Sigma$ is the ``linear'' embedding  (\ref {eqn:emb-lin}). We show below (Proposition \ref{prop:extend-map}) that, under certain conditions, the existence of such a family implies that of a corresponding map in $\Coalg AF$.

The  important result below follows immediately from the Cobar Duality Theorem in \cite {HPS}. 

\begin{thm}\label{thm:induction} \cite{HPS} 
There is a full and faithful functor, called the \emph {induction functor},
$$\ind : \Coalg AF\to \aalg$$
defined on objects by $\ind (C)=\Om C$ for all $C\in \ob \Coalg AF$ and on morphisms by
$$\ind(\theta)|_{\si C}=\sum _{k\geq 1} (\si )^{\otimes k} \theta (-\otimes z_{k-1} ) s :\si C_{+} \to  \Om C'$$
for all  $\theta\in \Coalg AF(C,C')$.
\end{thm}

As an easy consequence of Theorem \ref{thm:induction}, we obtain the following crucial operadic characterization of $\mathbf {DCSH}$.

\begin{thm}\cite{HPS} There is an isomorphism of categories 
$$\Coalg AF \xrightarrow{\cong}\mathbf {DCSH}$$
defined to be the identity on objects and to be $\ind$ on morphisms.
\end{thm}

\begin{rmk}  Thanks to this operadic description of $\cat{DCSH}$, we see that strongly homotopy coalgebra maps are exactly morphisms of $\op A$-coalgebras up to $\op F$-parame\-triza\-tion.
\end{rmk}


\section {Monoidal structures and modules over operads}\label{sec:mon-struct}

We carry out in this section a detailed study of monoidal structures on categories of modules and of coalgebras over a fixed operad $\op P$, in both their usual and enlarged, ``up-to-parametrization'' forms, with respect to some $\op P$-co-ring $\op Q$. We devote particular attention to the monoids in these categories, which we call $\op P$-rings (in $\cat {Mod}_{\op P}$), pseudo $\op P$-rings (in $\cat {Mod}_{\op P, \op Q}$), $\op P$-Hopf algebras (in $\coalg P$) and pseudo $\op P$-Hopf algebras (in $\Coalg PQ$). We begin by treating the general case, then specialize to $\Coalg AF$.

\subsection{Monoidal structures and  co-ring parametrizations}\label{subsec:monoidal-paramet}

Let $(\cat C, \otimes , I)$ be any symmetric monoidal category admitting coequalizers and coproducts.  Thanks to the existence and naturality of the intertwining map $\mathfrak i$ (Remark \ref{rmk:intertwine}), the level tensor product of symmetric sequences induces a symmetric monoidal structure $\wedge$ on the category of operads.  If $(\op P, \gamma)$ and $(\op P', \gamma')$ are operads, then $(\op P,\gamma)\wedge (\op P',\gamma'):=(\op P\otimes \op P', \gamma'')$, where $\gamma''$ is the composite
$$(\op P\otimes \op P')\diamond (\op P\otimes \op P')\xrightarrow{\mathfrak i}(\op P\diamond \op P)\otimes (\op P'\diamond \op P')\xrightarrow{\gamma\otimes \gamma'}\op P\otimes \op P'.$$
The unit object with respect to the monoidal product $\wedge$ is $\op C$ (cf., Proposition \ref{prop:level-tensor}).

Henceforth, let $(\op P, \gamma, \Delta,\epsilon)$ be a \emph{Hopf operad}, i.e., a level comonoid in the category of operads: $\Delta :\op P\to \op P\wedge \op P$ is a coassociative morphism of operads that is counital with respect to $e: \op P\to \op C$, which is also a morphism of operads.  The category $\mathbf {Mod}_{\op P}$ of right $\op P$-modules then admits a symmetric monoidal product, also denoted $\wedge$, which is defined as follows.  If $(\op M, \rho)$ and $(\op M',\rho')$ are two right $\op P$-modules, then $(\op M,\rho)\wedge(\op M',\rho'):=(\op M\otimes \op M', \rho'')$, where $\rho''$ is the composite
$$(\op M\otimes \op M')\diamond \op P\xrightarrow{1\diamond \Delta} (\op M\otimes \op M')\diamond (\op P\otimes \op P)\xrightarrow{\mathfrak i}(\op M\diamond \op P)\otimes (\op M'\diamond \op P)\xrightarrow{\rho\otimes \rho'} \op M\otimes \op M'.$$
The unit object is $\op C$, endowed with the right $\op P$-action given by the composite
$$\op C\diamond \op P\xrightarrow{1\diamond \varepsilon}\op C\diamond \op C\xrightarrow{\gamma_{\op C}} \op C,$$
where $\gamma _{\op C}$ is the usual multiplication on $\op C$.

There is an induced, symmetric monoidal structure $(\coalg{P},\wedge, I)$ such that there is a natural isomorphism of functors $\op T(-\wedge -)\cong \op T(-)\wedge \op T(-)$ from $\cat C$ into $\cat {Mod}_{\op P}$.

The category of monoids in $(\mathbf{Mod}_{\op P},\wedge, \op P)$ and morphisms thereof is denoted $\mathbf {Ring}_{\op P}$.  We call the objects of this category \emph { $\op P$-rings}.  Restricting to monoids in  $(\coalg{P},\wedge, I)$, we obtain the category $\op P\text{-}\cat {Hopf}$ of \emph { $\op P$-Hopf algebras}. 

The categories ${}_{\op P}\mathbf {Mod}$ of left $\op P$-modules and ${}_{\op P}\mathbf {Mod}_{\op P}$ of $\op P$-bimodules over $(\op P, \gamma, \Delta)$ also admit symmetric, level-monoidal structures, defined analogously to that on $\mathbf {Mod}_{\op P}$.  The category of $\op P$-bimodules is endowed with a second, nonsymmetric monoidal structure derived from the composition structure. Given two $\op P$-bimodules $\op M$ and $\op N$, their \emph {composition product over $\op P$}, denoted $\op M\underset {\op P} {\diamond} \op N$, is defined to be the obvious coequalizer.  Naturality arguments show that the intertwining map induces a natural morphism of $\op P$-bimodules
$$\xymatrix{\mathfrak i: (\mathcal X\wedge \mathcal X')\pcirc (\mathcal Y\wedge \mathcal Y')\ar [r]&(\mathcal X\pcirc \mathcal Y)\wedge (\mathcal X'\pcirc\mathcal Y')}$$
intertwining $\wedge$ and $\pcirc$.

\begin{defn}\label{defn:comon-coring} A \emph{level-comonoidal $\op P$-co-ring} is a $\op P$-co-ring $(\op Q, \psi_{\op Q}, \varepsilon_{\op Q})$  endowed with a coassociative, level comultiplication
$$\Delta _{\op Q}: \op Q\to \op Q\wedge \op Q,$$ 
which is counital with respect to 
$$e _{\op Q}:\op Q\to \op C.$$
Furthermore, the diagrams
$$\xymatrix{
\op Q\ar [dd]^{\Delta _{\op Q}}\ar [rr] ^{\psi_{\op Q}}&&\op Q\pcirc \op Q\ar [d]^{\Delta _{\op Q}\pcirc \Delta _{\op Q}}\\
&&(\op Q\wedge \op Q)\pcirc (\op Q\wedge \op Q)\ar [d] ^{\mathfrak i}\\
\op Q\wedge \op Q\ar [rr]^{\psi_{\op Q}\wedge \psi _{\op Q}}&&(\op Q\pcirc \op Q)\wedge (\op Q\pcirc \op Q)}$$
and
$$\xymatrix{
\op Q\ar[r]^{\Delta _{\op Q}}\ar [d]^{\varepsilon_{\op Q}}&\op Q\wedge \op Q\ar [d]^{\varepsilon_{\op Q}\wedge \varepsilon_{\op Q}}\\
\op P\ar[r]^{\Delta_{\op P}}&\op P\wedge \op P }$$
must commute.
\end{defn}

Recall  the ``inclusion" functor (\ref{eqn:emb-fat})
$$\mathfrak I_{\op Q}:\mathbf {Mod}_{\op P}\to \mathbf {Mod}_{\op P, \op Q}.$$
Restricting $\mathfrak I_{\op Q}$ to $\coalg{P}$ defines an ``inclusion" functor
\begin{equation}\label{eqn:incl-coalg}\mathfrak I_{\op Q}:\coalg{P}\to \Coalg PQ,\end{equation}
where $\mathfrak I_{\op Q}$ is defined on a morphism $f:C\to C'$ by $\mathfrak I_{\op Q}(f):=\op T(f)\pcirc \varepsilon _{\op Q}$.

Bringing the level comultiplication $\Delta_{\op Q}$ into play, we can define a symmetric monoidal product $\curlywedge$ on $\mathbf {Mod}_{\op P, \op Q}$ as follows, so that the restriction to $\cat{Mod}_{\op P}$ is the same as $\wedge$.   On objects, $\op M\curlywedge \op N$ is the same  as $\op M\wedge \op N$ in $\mathbf {Mod}_{\op P}$, while the monoidal product $\theta \curlywedge \theta'$ of $\theta \in  \mathbf {Mod}_{\op P, \op Q}(\op M, \op N)$ and $\theta' \in  \mathbf {Mod}_{\op P, \op Q}(\op M',\op N')$ is defined to be the following composite of morphisms of strict $\op P$-bimodules.
$$(\op M\wedge \op M')\pcirc\op Q\xrightarrow{1\pcirc\Delta _{\op Q}}(\op M\wedge \op M')\pcirc(\op Q\wedge \op Q)\xrightarrow{\mathfrak i} (\op M \pcirc\op Q)\wedge  (\op M '\pcirc\op Q)\xrightarrow{\theta \wedge \theta '} \op M\wedge \op M'$$
The compatibility of the two comultiplications implies that $-\curlywedge -$ is indeed a bifunctor. The coassociativity of $\Delta_{\op Q}$ ensures the associativity of $\curlywedge$, while the counit of $\Delta _{\op Q}$ gives rise to the unit of $\curlywedge$.  
By restriction, and using that $\op T$ is strong monoidal, we obtain a monoidal structure $\curlywedge$ on $\Coalg{P}{Q}$, which is the usual monoidal product of $\op P$-coalgebras on objects.

\begin{defn}\label{defn:pseudo}The category of monoids in $\mathbf {Mod}_{\op P, \op Q}$ with respect to the monoidal product $\curlywedge$ and of their morphisms is denoted $\mathbf{PsRing}_{\op P,\op Q}$.  We call its objects \emph {pseudo $\op P$-rings}, suppressing explicit mention of the governing comultiplication $\psi _{\op Q}$. Restricting to $\curlywedge$-monoids in  $\Coalg PQ$, we obtain the category $(\op P,\op Q)\text{-}\cat {PsHopf}$ of \emph {pseudo $\op P$-Hopf algebras}.  
\end{defn}

If  $(\op B,\mu)$ is a $\op P$-ring, where $\mu\in \mathbf {Mod}_{\op P}(\op B\wedge \op B, \op B)$ is the product map, it is clear that $\big(\mathfrak I_{\op Q}(\op B), \mathfrak I_{\op Q}(\mu)\big)$ is a pseudo $\op P$-ring.  In other words, $\mathfrak I_{\op Q}$ induces an ``inclusion" functor
$$\mathfrak I_{\op Q}:\mathbf{Ring}_{\op P}\to \mathbf{PsRing}_{\op P, \op Q}.$$
Similarly, there is an induced, ``inclusion" functor
$$\mathfrak I_{\op Q}:\op P\text{-}\cat {Hopf}\to (\op P,\op Q)\text{-}\cat {PsHopf}.$$

When the $\op P$-bimodule $\op Q$ is a free bimodule, there exist ``free" constructions in the category of pseudo $\op P$-Hopf algebras, as specified in the next proposition.  Before stating the proposition, we state and prove a crucial lemma, which is useful elsewhere in this article as well, then introduce one necessary definition.

Restricting to $\cat C= \cat {grMod}_{R}$ or $\cat C=\cat {Ch}_{R}$, let $\op L:\mathbf C\to \mathbf C^\Sigma$ be the ``linear'' embedding of (\ref{eqn:emb-lin}).  Let $\mathfrak u:\op L\to \op T$ denote the obvious ``inclusion on level $1$" natural transformation.

\begin{lem}\label{lem:inducedT} Let $A$ and $B$ be graded $R$-modules, and let $\op X$ be a symmetric sequence of graded $R$-modules.  Any morphism $\theta:\op L(A)\diamond \op X\to \op T(B)$ in $\cat{grMod}_R^\Sigma$ gives rise naturally to a morphism $\widehat \theta: \op T(A)\diamond \op X\to \op T(B)$ of symmetric sequences such that $\widehat\theta(\mathfrak u\diamond Id)=\theta$.
\end{lem}

 \begin{proof}  Recall from Definition \ref{defn:comp-prod} that, in the definition of the composition product of symmetric sequences $\op X$ and $\op Y$,  the left action of $\Sigma _{k}$ on $\coprod_{\substack{ k\geq 1\\ \vec n\in I_{k,n}}} \op Y(n_1)\otimes \cdots\otimes \op Y(n_k)$ is given by
 $$\sigma \cdot(y_{1}\otimes \cdots \otimes y_{k})=y_{\sigma ^{-1} (1)}\otimes \cdots \otimes y_{\sigma ^{-1}(k)}$$
 for all $\sigma \in \Sigma _{k}$ and $y_{i}\in\op Y(n_{i})$, $1\leq i\leq k$.
 
 For all $m$, $n$ and $\vec n\in I_{m,n}$, define 
 $$\widehat\theta _{\vec n}: \op T(A)(m)\otimes \big (\op Y(n_{1})\otimes\cdots \otimes \op Y(n_{m})\big)\longrightarrow\op T(B)(n)$$
 to be the composite
 $$\xymatrix{A^{\otimes m}\otimes \big (\op Y(n_{1})\otimes\cdots \otimes \op Y(n_{m})\big)\ar[r]^{\cong}&(A\otimes \op Y (n_{1}))\otimes \cdots\otimes (A\otimes \op Y (n_{m}))\ar [d]_{\theta ^{\otimes m}}\\
&B^{\otimes n_{1}}\otimes \cdots \otimes B^{\otimes n_{m}}\ar  [d]_{=}\\
&B^{\otimes n}.}$$

Since $\theta$ is a morphism of symmetric sequences, the second diagram in Lemma \ref{lem:classical} commutes for $\op X=\op T(A)$ and $\op Z=\op T(B)$.  The first diagram commutes  in this case as well because for all $a_{1}\otimes \cdots \otimes a_{m}\in A^{\otimes m}$ and all $y_{1}\otimes \cdots\otimes y_{m}\in \op Y(n_{1})\otimes \cdots \otimes \op Y(n_{m})$,
\begin{align*}
&\widehat\theta _{\sigma^{-1}\vec n}(\sigma \otimes \sigma ^{-1})\big((a_{1}\otimes\cdots\otimes a_{m})\otimes (y_{1}\otimes\cdots \otimes y_{m})\big)\\
&=\widehat\theta _{\sigma^{-1}\vec n}\big((a_{\sigma(1)}\otimes \cdots \otimes a_{\sigma (m)})\otimes (y_{\sigma(1)}\otimes\cdots \otimes y_{\sigma (m)})\big)\\
&=\pm \theta (a_{\sigma (1)}\otimes y_{\sigma (1)})\otimes \cdots \otimes \theta (a_{\sigma (m)}\otimes y_{\sigma (m)})\\
&=\pm \sigma (n_{\sigma (1)},...,n_{\sigma (m)})\left(\theta (a_{1}\otimes y_{1})\otimes \cdots \otimes \theta (a_{m}\otimes y_{m})\right)\\
&= \sigma (n_{\sigma (1)},...,n_{\sigma (m)})\widehat \theta _{\vec n}\big((a_{1}\otimes \cdots \otimes a_{m})\otimes (y_{1}\otimes \cdots \otimes y_{m})\big).
\end{align*}
\end{proof}

\begin{defn}  Let $H$ be a $\op P$-Hopf algebra in $\cat{grMod}_R$ or $\cat{Ch}_R$. A \emph {free algebraic $\op P$-Hopf extension} of $H$ by a generator $v$ consists of a morphism of $\op P$-Hopf algebras $j:H\to H'$ such that 
the underlying morphism of graded algebras is the inclusion map $H\hookrightarrow H\coprod Tv$,  where $H\coprod Tv$ is the coproduct in the category of graded algebras of $H$ and of the free algebra on $v$.
\end{defn}

We first explain in what sense free algebraic $\op P$-Hopf extensions truly are free, in the nondifferential setting.  

\begin{prop}\label{prop:free-ext} Let $\cat C=\cat{grMod}_{R}$.  Let $(\op Q,\psi_{\op Q},\varepsilon_{\op Q}, \Delta _{\op Q},e_{\op Q})$ be a level-comonoidal $\op P$-co-ring (cf., Definition \ref{defn:comon-coring}) that is free as a $\op P$-bimodule, generated by $\op X$.  Let $H\coprod Tv$ be a free algebraic $\op P$-Hopf extension of  a $\op P$-Hopf algebra $H$, and let $H'$ be another $\op P$-Hopf algebra.
  
For all $\theta \in (\op P,\op Q)\text{-}\cat {PsHopf}(H,H')$ and $\lambda \in \cat{Mod}_R^\Sigma\big(\op L(R\cdot v)\diamond \op X, \op T(H')\big)$, there is a unique morphism 
$$\widehat{(\theta+\lambda)}\in (\op P,\op Q)\text{-}\cat {PsHopf}(H\coprod Tv,H')$$
extending $\theta$ and $\lambda$.
\end{prop}

The proof of this proposition, which is somewhat technical, is in the appendix.

\begin{cor} Let $(\op Q,\psi_{\op Q},\varepsilon_{\op Q}, \Delta _{\op Q},e_{\op Q})$ be as in the statement of Proposition \ref{prop:free-ext}.  
Let $H$ and $H'$ be $\op P$-Hopf algebras.  If the underlying algebra of $H$ is free on a free graded $R$-module $V$ that is bounded below, then  
for all $\lambda \in \cat{grMod}_R^\Sigma\big(\op L(V)\diamond \op X, \op T(H')\big)$, there is a unique morphism 
$$\widehat\lambda\in (\op P,\op Q)\text{-}\cat {PsHopf}(H,H')$$ 
extending  $\lambda$.
\end{cor}

More informally, we can say that  if $\op P$ is free as a bimodule, then a pseudo-$\op P$-Hopf algebra map with domain free as an algebra is specified by its values on generators of $\op P$ and of the domain.

\begin{proof}  The proof proceeds by induction on degree of elements in a basis of $V$, starting in the lowest degree $k$  for which $V_k\not=0$, applying Proposition \ref{prop:free-ext} at each step. Here, $\theta$ is taken to be the unique morphism with domain $0$.
\end{proof}

\subsection{Application to the Alexander-Whitney co-ring}\label{subsec:alex-whit}

Specializing to the case where $\op P=\op A$, the associative operad, and $\op Q=\op F$, the Alexander-Whitney co-ring (cf., Section \ref{subsec:dcsh}), we explain how to verify that a pseudo-$\op A$-Hopf morphism with free domain respects differential structure. 
The proof of Theorem \ref{thm:induction} relies implicitly on the following proposition, which comes in handy later in this article as well.  

Recall the notion of a family $\mathfrak F(\theta)$ induced by $\theta \in \Coalg AF(C,C')$ from (\ref{eqn:family}) and of a derivation $\op D(t):\op T(A)\to \op T(B)$ induced by a morphism $t:A\to B$ of graded $R$-modules (Definition \ref{defn:derivation}).

\begin{prop}\label{prop:extend-map}  Let $\cat C=\cat {Ch}_{R}$.  Let $H\coprod Tv$ be a free algebraic $\op A$-Hopf extension of an $\op A$-Hopf algebra $H$, and let $H'$ be another $\op A$-Hopf algebra. Let $\Delta$ and $\Delta'$ denote the comultiplications and $d$ and $d'$ the differentials on $H\coprod Tv$ and on $H'$. Let 
$$\theta \in (\op A,\op F)\text{-}\cat {PsHopf}(H,H')$$  
with induced family $\mathfrak F(\theta)=\{\theta _k\mid k\geq 1\} $.

For any set $\{\lambda _k\in \op T(H')(k)\mid k\geq 1\}$ such that for all $k$, 
\begin{equation}
\op D(d')_{k}\lambda _k -\op D(\overline \Delta ')_{k-1}\lambda _{k-1}=\theta _k(dv)-\sum _{i+j=k} (\theta _i\otimes \theta _j)\overline \Delta (v), 
\end{equation}
$\theta $ can be extended to 
$$\widehat \theta\in (\op A,\op F)\text{-}\cat {PsHopf}(H\coprod Tv, H')$$ such that $\widehat\theta (v\otimes z_{k-1})=\lambda _k$ for all $k$.
\end{prop}

Thanks to this proposition, if $TV$ is a chain Hopf algebra with free underlying algebra and $H'$ is any chain Hopf algebra,  it is possible to construct monoidal morphisms in $\Coalg AF$ from $TV$ to $H'$ by induction on the generators  $V$.

\begin{proof} The family $\{\lambda_k\mid k\geq 1\}$ is equivalent to a morphism of symmetric sequences of graded $R$-modules $\lambda :\op L(R\cdot v)\diamond \op S\to \op T(H')$.
We can therefore apply Proposition \ref {prop:free-ext} to obtain $\widehat \theta$ as a morphism of nondifferential objects.  On the other hand, as we can see from the definition of $\del _{\op F}$, the hypothesis on the family $\{\lambda_{k}\}$ is exactly the condition that must be satisfied for $\widehat\theta$ to be a differential map. 
 \end{proof}

We recall now the relationship between the functor $\ind$ (Theorem \ref {thm:induction}) and the monoidal structures on the source and target categories, as developed in \cite{HPS} and \cite {HPST}.

\begin{lem} \cite {HPST} The induction functor $\ind : \Coalg AF\to \aalg$ is comonoidal, i.e., there is a natural transformation of functors into associative chain algebras
$$q: \ind(-\curlywedge-)\to \ind(-)\otimes \ind (-),$$
which is given by the Milgram equivalence (\ref{eqn:milgram}) on objects.
\end{lem}

Throughout the remainder of this article, we consider objects in 
the following category derived from $\Coalg AF$.  Recall $z_{k}$ is the generator of $\op F$ in level $k+1$, which is of degree $k$.

\begin{defn}  The objects of the \emph {weak Alexander-Whitney category} $\mathbf {wF}$ are pairs $(C,\Psi)$, where $C$ is a object in $\acoalg$ and $\Psi\in \Coalg AF(C, C\otimes C)$ such that
$$\Psi (-\otimes z_0):C\to C\otimes C$$
is exactly the comultiplication on $C$, while
$$\mathbf {wF}\big( (C,\Psi), (C',\Psi')\big)= \{ \theta \in \Coalg AF(C,C')\mid \Psi' \theta =(\theta\curlywedge \theta )\Psi\}.$$ 
The objects of $\cat {wF}$   are called \emph {weak Alexander-Whitney coalgebras}.
\end{defn}

As noted in the next lemma, the cobar construction provides an important link between the weak Alexander-Whitney category and the following category of algebras endowed with comultiplications, which are not necessarily coassociative.

\begin{defn} The objects of the \emph {weak Hopf algebra category} $\mathbf {wH}$ are pairs $(A,\psi)$, where $A$ is a chain algebra over $R$ and $\psi:A\to A\otimes A$ is a map of chain algebras, while
$$\mathbf{wH}\big((A,\psi), (A',\psi ')\big)=\{ f\in \aalg(A,A')\mid \psi' f=(f\otimes f)\psi\}.$$
\end{defn}

\begin{lem} \cite {HPST} \label{lem:extcobar}The cobar construction extends to a functor 
$$\tom : \mathbf {wF}\to \mathbf {wH},$$ 
given by $\tom (C,\Psi)=\big(\Om C, q \ind(\Psi)\big)$, where $\ind(\Psi):\Om C\to \Om(C\otimes C)$, as in Theorem \ref{thm:induction}, $q:\Om (C\otimes C)\to \Om C\otimes \Om C$ is Milgram's equivalence (\ref{eqn:milgram}) and $\tom \theta =\ind (\theta):\Om C\to \Om C'$ for all  $\theta \in \mathbf {wF}\big( (C,\Psi), (C',\Psi')\big)$.
\end{lem}

Motivated by topology, we are particularly interested in those objects $(C, \Psi)$ of $\cat {wF}$ for which $\tom (C, \Psi)$ is actually a strict Hopf algebra, i.e., such that $q\ind(\Psi)$ is coassociative.

\begin{defn}\label{defn:F} The \emph {Alexander-Whitney category} $\cat F$ is the full subcategory of $\cat {wF}$ such that $(C,\Psi)$ is an object of $\cat F$ if and only if $q\ind(\Psi)$ is coassociative.  We call the objects of $\cat F$ \emph {Alexander-Whitney coalgebras}.
\end{defn}

From Lemma \ref{lem:extcobar}, it is clear that $\tom$ restricts to a functor 
\begin{equation}\label{eqn:tom}
\tom:\cat F\to \cat H
\end{equation}
where $\cat H=\op A\text{-}\mathbf {Hopf}$ is the usual category of chain Hopf algebras. 

We can now explain the topological importance of the category $\mathbf F$. 

\begin{thm}\label{thm:loopmodel}\cite{HPST} There is a functor $\tC:\cat {sSet}_1\to \cat {F}$ from the category of $1$-reduced simplicial sets to the Alexander-Whitney category such that the coassociative chain coalgebra underlying $\tC(K)$ is $C(K)$, the normalized chains on $K$.  Furthermore, there is a natural quasi-isomorphism of chain algebras
$$\tom\tC(K)\xrightarrow{\simeq}C(GK)$$
that is also itself a DCSH map.
\end{thm}

\section {Extended multiplicative naturality of Cotor}\label{sec:ext-mult-nat}

Let $C$ be a chain coalgebra with comultiplication $\Delta$.  If $(M,\nu)$ and $(M',\nu')$ are right and left $C$-comodules, respectively, then their \emph{cotensor product over $C$} is
$$M\square_{C}M'=\ker (M\otimes M'\xrightarrow{\nu \otimes 1-1\otimes \nu'}M\otimes C\otimes M').$$
In particular, if we endow the ground ring $R$ with its trivial left $C$-comodule structure, then 
$$M\square_{C} R\cong\{x\in M\mid \nu(x)=x\otimes 1\},$$
so that $M\square_{C}R$ can be seen as a graded submodule of $M$, which we can think of as the ``cofixed points'' of the coaction $\nu$. 

In this section we study the derived functor of cotensor product, $\cotor_{C}(M,M')$.   We begin by recalling the formula for Cotor in terms of one-sided cobar constructions, from which it is immediately clear that Cotor is natural in all three variables, with respect to morphisms of comodules over a fixed coalgebra and with respect to morphisms of coalgebras.  In \cite{GM} Gugenheim and Munkholm proved an ``extended naturality'' result for Cotor, i.e., that Cotor is actually functorial with respect to a much larger class of morphisms. In section \ref{subsec:comodules} we reformulate Gugenheim and Munkholm's result in operadic language.

We show in section \ref{subsec:one-sided} that if $H$ is a chain Hopf algebra and $B$ is a (left) $H$-comodule algebra, then $\cotor_{H}(R,B)$ admits a graded multiplicative structure, which is natural in both variables, with respect to morphisms of comodule algebras over fixed Hopf algebras and with respect to morphisms of Hopf algebras.  We then prove in section \ref{subsec:comodules} that there is a larger class of morphisms, the class of \emph{comodule-algebra maps up to strong homotopy (CASH maps)}, with respect to which the multiplicative structure of $\cotor_{H}(R,B)$ is natural, i.e., we establish ``extended multiplicative naturality'' of Cotor, which plays an important role in sections \ref{sec:homol-htpyfib} and \ref{sec:topology}.

\subsection {Cotor: definition and naturality}\label{subsec:one-sided}

The derived functor of cotensor product, Cotor, can be calculated in terms of the following complex.

\begin{defn}\label{defn:one-sided} Let $C$ be a 
simply-connected chain coalgebra, and let $M$ be a right $C$-comodule.  The \emph{one-sided cobar construction} $M\otimes _{t_{\Om}}\Om C$ is the chain complex with underlying graded $R$-module $M\otimes Ts^{-1} C_{+}$ and with differential $D_{\Om}$ given by
$$D_{\Om }(x\otimes w)=dx\otimes w +(-1)^x x\otimes d_{\Om }w+(-1)^{x_{i}}x_{i}\otimes (s^{-1} c^{i}\cdot w),$$
where $x\in M$, $w\in \Om C$, $d$ is the differential on $M$, $d_{\Om }$ is the cobar construction differential (cf., Notation and conventions), $\nu (x)=x_{i}\otimes c^{i}$ and $s^{-1} 1 := 0$.  

There is an analogous definition of $\Om C\otimes _{t_{\Om}} N$ for any left $C$-comodule $N$. 
\end{defn}

\begin{rmk}  If $M=C$ or $N=C$, we obtain the usual \emph{acyclic cobar constructions}:
$$C\otimes _{t_{\Om}}\Om C\qquad \text{and}\qquad \Om C\otimes _{t_{\Om}} C.$$ 
\end{rmk}

\begin{rmk}  The formula in the definition above makes it clear that there are functors
$$-\otimes \Om C:\cat {Comod}_{C}\to \cat {Ch}_{R}$$
and
$$\Om C\otimes -:{}_{C}\cat {Comod}\to \cat{Ch}_{R}.$$
\end{rmk}

One-sided cobar constructions can be applied  to calculations of $\cotor$, the derived functor of the cotensor product.  Let $C$ be a connected coalgebra, and let $M$ and $N$ be right and left comodules over $C$, respectively.  Then, as shown in e.g., \cite {EM}, 
\begin{equation}\label{eqn:cotor}
\cotor ^C(M,N)=\H \big((M\otimes _{t_{\Om}}\Om C)\otimes _{\Om 
C}(\Om C\otimes _{t_{\Om }}N)\big).
\end{equation}
It follows from the previous remark that Cotor is a bifunctor
$$\cotor: \cat {Comod}_{C}\times {}_{C}\cat {Comod}\to \cat {grMod}_{R}.$$
We think of this as the \emph{linear naturality} of Cotor.

\begin{rmk} We can also use the cobar construction to define the homotopy fiber of a morphism of coaugmented chain coalgebras $f:C'\to C$.  Consider the projection $\pi:C\otimes _{t_{\Om}}\Om C\to C$, which is a surjective morphism of chain complexes with contractible source, and therefore an acceptable candidate for a fibrant replacement of the coaugmentation $\eta: R\to C$.  Consequently, we can define the \emph{homotopy fiber} of $f$ to be the pullback of
$$C\otimes _{t_{\Om}}\Om C\xrightarrow\pi C\xleftarrow f C',$$
i.e., $C'\otimes _{t_{\Om}}\Om C$.  The homology of the homotopy fiber of $f$ is thus exactly $\cotor ^C(C',R)$.
\end{rmk}

Let $H$ be a chain Hopf algebra. Recall that a chain algebra $B$ that is also an $H$-comodule is an \emph{$H$-comodule algebra} if the $H$-coaction map is a morphism of chain algebras. In \cite {Mil}, Miller proved the existence of a natural chain algebra structure on the one-sided cobar construction $\Om H\otimes _{t_{\Om }}B$, for any commutative Hopf algebra $H$ and any left $H$-comodule algebra $B$.   Here we dualize Theorem 4.1 of \cite {fht}, obtaining a generalization of Miller's result to any chain Hopf algebra $H$. As a consequence, $\cotor ^H(R,B)$ admits a natural multiplicative structure for any Hopf algebra $H$ and any $H$-comodule algebra $B$. 

We begin by considering a special case: the acyclic cobar construction.  Though it would be possible to prove the next proposition and its corollaries by appealing to Theorem 4.1 in \cite {fht} and then dualizing, we prefer to give a direct, constructive proof, since the explicit formulas we provide are much simpler than those in the dual case and prove quite useful.

\begin{prop}\label{prop:cobar-alg} If $H$ is any chain Hopf algebra, then the 
free left $\Om H$-module structure on $\Om H\otimes _{t_{\Om }}H$ can be extended to a chain algebra 
structure  such that
$$(1\otimes c)(\si a\otimes 1)=(-1)^{(a+1)c} \si a\otimes c +(-1)^c\si 
(c\cdot a)\otimes 1 +(-1)^{c+ac^i} \si (c_{i}\cdot a)\otimes c^i$$
for all $a,c\in H$, where  $\overline \Delta (c)=c_{i}\otimes c^i$ and 
$$(1\otimes c)(1\otimes e)=1\otimes c\cdot e$$
for all $c,e\in H$.
\end{prop}

\begin{proof} Given the multiplication as partially defined in the 
statement of the proposition, we extend it to all of $\Om H\otimes 
_{t_{\Om }}H$ by associativity, which is possible since $\Om H$ is 
free as an algebra on $\si H_{+}$.  Hence all that we must do is 
verify that $D_{\Om}$ is a derivation with respect to this product.  
We do the second case and leave the first, the proof of which is quite 
similar, to the reader.

If image of $c$ and $e$ under the reduced comultiplication are $c_{i}\otimes c^i$ and 
$e_{j}\otimes e^j$, respectively, then the image of 
$c\cdot e$ under the \emph{unreduced} comultiplication is
\begin{align*}
(c\cdot e)_{k}\otimes (c\cdot e)^k=&c\cdot e\otimes 1+ 1\otimes c\cdot e\\
&+ c\otimes e +(-1)^{ec}e\otimes c\\
&+c_{i}\otimes c^i\cdot e+(-1)^{ec_{i}}c_{i}\cdot e\otimes c^i\\
&+c\cdot e_{j}\otimes e^j+(-1)^{e_{j}c}e_{j}\otimes c\cdot e^j\\
&+ (-1)^{e_{j}c^i}c_{i}\cdot e_{j}\otimes c^i\cdot e^j.
\end{align*} 
Consequently, 
\begin{align*}
D_{\Om}(1\otimes c\cdot e)=& 1\otimes d(c\cdot e)-\si (c\cdot 
e)\otimes 1\\
&-\si c\otimes e -(-1)^{ec}\si e\otimes c\\
&-\si c_{i}\otimes c^i\cdot e-(-1)^{ec_{i}}\si (c_{i}\cdot e)\otimes c^i\\
&-\si (c\cdot e_{j})\otimes e^j-(-1)^{e_{j}c}\si e_{j}\otimes c\cdot e^j\\
&- (-1)^{e_{j}c^i}\si (c_{i}\cdot e_{j})\otimes c^i\cdot e^j.
\end{align*}
On the other hand
$$D_{\Om}(1\otimes c)\cdot (1\otimes e)=1\otimes dc\cdot e-\si c\otimes e-\si 
c_{i}\otimes c^i\cdot e ,$$
while 
\begin{align*}
(-1)^c (1\otimes c)\cdot D_{\Om }(1\otimes e)=&(-1)^c 1\otimes c\cdot 
de +(-1)^c(1\otimes 
c)\cdot (-\si e\otimes 1-\si e_{j}\otimes e^j)\\
=&-(-1)^{ce}\si e\otimes c -\si (c\cdot  e)\otimes 1\\
&-(-1)^{ec^i} \si 
(c_{i}\cdot e)\otimes c^i
-(-1)^{e_{j}c}\si e_{j}\otimes c\cdot e^j\\
&-\si (c\cdot e_{j})\otimes 
e^j-(-1)^{e_{j}c^i}\si (c_{i}\cdot e_{j})\otimes c^i\otimes 
e^j.\end{align*}
It is now obvious that 
$$D_{\Om}(1\otimes c\cdot e)=D_{\Om}(1\otimes c)\cdot (1\otimes 
e)+(-1)^c (1\otimes c)\cdot D_{\Om }(1\otimes e). $$
\end{proof}

This proposition admits the following generalization, the proof of 
which is essentially identical to that of the proposition.

\begin{cor}\label{cor:cobar-alg} If $H$ is any chain Hopf algebra and $B$ is a 
left $H$-comodule algebra, with coaction map $\nu$, then the 
free left $\Om H$-module structure on $\Om H\otimes _{t_{\Om }}B$ can be extended to a chain algebra 
structure  such that
$$(1\otimes b)(\si a\otimes 1)=(-1)^{(a+1)b} \si a\otimes b +(-1)^c\si 
(c\cdot a)\otimes 1 +(-1)^{b+ab^i} \si (c_{i}\cdot a)\otimes b^i$$
for all $a\in H$, $b\in B$, where $\nu (b)= c\otimes 1 + 1\otimes b +c_{i}\otimes 
b^i $ and 
$$(1\otimes b)(1\otimes b')=1\otimes b\cdot b'$$
for all $b,b'\in B$.
\end{cor}

An analogous results clearly holds for right comodule algebras as well.

The multiplicative structure defined above is easily seen to be natural, in the following sense.  Let $f:H\to H'$ be a morphism of chain Hopf algebras.  Let $(B,\nu)$ be a left $H$-comodule algebra, and let $(B',\nu')$ be a left  $H'$-comodule algebra.  Let $g:B\to B'$ be morphism of chain algebras such that $(f\otimes g)\nu=\nu 'g$.  It is easy to check that the chain map 
$$\Om f\otimes g:\Om H\otimes _{t_\Om }B\to  \Om H'\otimes _{t_\Om }B'$$
respects the multiplicative structure defined in Corollary \ref {cor:cobar-alg}.   

Let  $\cat{CA}^\ell$ denote the following category. Objects are pairs $(H;B)$, where $H$ is a chain Hopf algebra and $B$ is a left $H$-comodule algebra, both over the fixed PID $R$.  A morphism from $(H;B)$ to $(H';B')$ is a pair $(f;g)$, where $f:H\to H'$ is a morphism of chain Hopf algebras and $g:B\to B'$ is a chain algebra map such that $(f\otimes g)\nu=\nu 'g$, where $\nu$ and $\nu'$ are the coactions on $B$ and $B'$, respectively.  The analogous category for right comodule algebras is denoted $\cat {CA}^r$.

\begin{cor}\label{cor:natl-ca} Let $\cat{chA}$ denote the category of chain algebras over $R$. There are functors
$$\Om^{\ell} (-;-):\cat {CA}^\ell\to \cat{chA}$$
defined by $\Om (H;B)= \Om H\otimes _{t_{\Om }}B$ (with the multiplicative structure of Corollary \ref {cor:cobar-alg}) and $\Om (f;g)=\Om f\otimes g$ and 
$$\Om^{r} (-;-):\cat {CA}^r\to \cat{chA}$$
defined by $\Om (H;B)= B\otimes _{t_{\Om }}\Om H$ (with the multiplicative structure analogous to that of Corollary \ref {cor:cobar-alg}) and $\Om (f;g)=g\otimes \Om f$.
\end{cor}

We show in section \ref{subsec:comodules} that  the algebra structure on $\Om H\otimes _{t_\Om }B$ is actually natural with respect to a bigger class of morphisms than those of the category $\cat{CA}^\ell $.  This extended naturality of the algebra structure of $\Om H\otimes _{t_\Om }B$ plays a crucial role in section \ref{sec:topology}.

\begin{cor}\label{cor:cotor-alg}  Let $\cat{grA}$ denote the category of graded algebras over $R$. The functor $\cotor$ restricts and corestricts to functors
$$\cotor^{(-)} (R;-):\cat {CA}^\ell\to \cat{grA}$$
and
$$\cotor^{(-)} (-;R):\cat {CA}^r\to \cat{grA}$$
\end{cor}

This is the \emph{multiplicative naturality} of Cotor.

\begin{proof} From the definition of $\cotor$ (\ref{eqn:cotor}), we see that if $B$ is a left $H$-comodule algebra, then 
$$\cotor^H(R,B)=\H\big( \Om H\otimes _{t_{\Om }}B ).$$
The previous corollary implies that there is a natural graded algebra structure on $\H\big( \Om H\otimes _{t_{\Om }}B ).$  The right-comodule case works similarly.
\end{proof}

\begin{rmk} As a consequence of Proposition \ref{prop:cobar-alg}, we obtain that $\Om H\otimes _{t_{\Om }}H$ is itself an $H$-comodule algebra.  To establish this fact, we must show that the following diagram commutes.
$$\xymatrix{(\Om H\otimes _{t_{\Om }}H)^{\otimes 2}\ar [d]_{(1\otimes \Delta )^{\otimes 2}}\ar [r] ^\mu &\Om H\otimes _{t_{\Om }}H\ar [d]_{(1\otimes \Delta )}\\
\big (\Om H\otimes _{t_{\Om }}(H\otimes H)\big)^{\otimes 2}\ar [r]^\mu &\Om H\otimes _{t_{\Om }}(H\otimes H)}$$
Here, $\Delta$ is the comultiplication on $H$ and $\mu $ denotes the multiplication defined in the statements of Proposition \ref{prop:cobar-alg} and Corollary \ref{cor:cobar-alg}.  The left $H$-comodule structure on $H\otimes H$ is given by $\Delta \otimes 1$.  It suffices to check that this diagram commutes on elements of the form $1\otimes c\otimes \si a\otimes 1$, which is not a difficult computation. The coassociativity of $\Delta $ plays a crucial role in this verification.

Analagously, $H\otimes 
_{t_{\Om}}\Om H$ is also an $H$-comodule algebra.
\end{rmk}

\subsection {Maps of comodules up to strong homotopy}\label{subsec:comodules}

In this article we need relative versions of the results  from \cite {HPST} cited in section \ref{subsec:dcsh}, to establish conditions under which there is a multiplicative map between one-sided cobar constructions of the sort considered in Corollary \ref{cor:cobar-alg}.   As a consequence, we obtain Proposition \ref{prop:mult-cotor}, which is both a multiplicative generalization of the extended naturality of $\cotor$, due to Gugenheim and Munkholm \cite {GM} (Theorem \ref{thm:gm-cotor} below), and an extended version of  the multiplicative naturality of Cotor (Corollary \ref{cor:cotor-alg} above).

Let $C$ and $C'$ be connected chain coalgebras. Recall from section \ref{subsec:dcsh} that though every chain coalgebra map $f:C\to C'$ induces a chain algebra map $\Om f:\Om C\to \Om C'$, not all chain algebra maps $\varphi:\Om C\to \Om C'$ are so induced.   In particular, the category of chain coalgebras can be seen as a wide, but not full, subcategory of the category $\cat{DCSH}$, so that the morphisms in $\cat {DCSH}$ can be considered as ``weak'' chain coalgebra morphisms.

In this section we consider an analogous weakening of the morphisms in the category $\cat {CA}^r$ defined  in section \ref{subsec:one-sided}, for which we provide equivalent chain-level and operadic definitions, both of which are quite useful. The operadic definition serves to facilitate the proofs of the existence results (Propositions \ref{prop:rel-free-ext} and \ref{prop:cash-coproduct}) that play a key role in section \ref{sec:homol-htpyfib}.

The reader who is not interested in the fine details of our constructions and existence results can safely limit his perusal of this section only to the definition of the category $\cat {CASH}$ (Definition \ref{defn:cash}) and to Proposition \ref{prop:cash-coproduct}.

Let $\mathbf C$ denote either the category of graded $R$-modules or the category of chain complexes over $R$. Let $\mathbf C^{\Sigma_+} $denote the category of \emph {shifted symmetric sequences} in $\mathbf C$.  An object $\mathcal X$ of $\mathbf C^{\Sigma_{+}}$ is  a family $\{\mathcal X(n)\in \mathbf C\mid n\geq 0\}$ of objects in $\mathbf C$ such that $\mathcal X(n)$ admits a right action of the symmetric group $\Sigma_{n-1}$, for all $n>0$ and such that $\op X(0)=0$.  A morphism in $\cat C^{\Sigma_{+}}$ from $\op X$ to $\op Y$ consists of a family 
$$\{\varphi_{n}\in \cat C\big(\op X(n),\op Y(n)\big)\mid \varphi_{n}\quad\text{is $\Sigma_{n-1}$-equivariant}, n\geq 1\}.$$
There is a faithful functor 
$$\mathcal T_r: \mathbf C\times \mathbf C\to\mathbf C^{\Sigma_+}$$ 
where, for all $n$, $\mathcal T_r(A,B)(n)=A\otimes B^{\otimes n-1}$, where $\Sigma _{n-1}$ acts by permuting the factors of $B$. 

The following useful operation on symmetric sequences in the image of $\op T_r$ is a shifted version of the notion of derivation of symmetric sequences in the image of $\op T$ (Definition \ref{defn:derivation}).

\begin{defn} \label{defn:shift-deriv}Let $f,s:A\to B$ and $g,h, t:C\to D$ be homogenous linear maps of graded $R$-modules, such that $f$, $g$ and $h$ are homogeneous of degree $0$, and $s$ and $t$ are homogeneous of degree $m$. The \emph {$(f,g,h)$-derivation of  shifted symmetric sequences} induced by $s$ and $t$ is the morphism of symmetric sequences 
$$\op D_{(f,g,h)}(s,t):\op T_r(A,C)\to \op T_r(B,D)$$ 
that is of degree $m$ in each level and that is defined as follows in level $n$.
$$\op D_{(f,g,h)}(s,t)_{n}=s\otimes h^{\otimes n-1}+ \sum _{j=0}^{n-2} f\otimes g^{\otimes j}\otimes t\otimes h^{\otimes n-j-2}$$
When $A=B$, $C=D$ and $f=Id_A$ and $g=Id_C=h$, we simplify notation and write $\op D(s,t)$ for the $(Id_A, Id_C, Id_C)$-derivation induced by $s$ and $t$.\end{defn}

It is obvious that there is again a \emph {level monoidal structure}  $(\mathbf C^{\Sigma_+}, \otimes, \op C)$, where $(\op X\otimes \op Y)(n)=\op X(n)\otimes\op Y(n)$, endowed with the diagonal action of $\Sigma_{n-1}$, and $\op C(n)=R$, endowed with the trivial $\Sigma _{n-1}$-action. By proofs analogous to those in \cite[section
II.1.8]{MSS}, we can show that  the category $\mathbf C^{\Sigma_+}$ also admits a right action by the monoidal category $(\cat C^\Sigma,\diamond, \op J)$, i.e., there is a bifunctor
$$\cat C^{\Sigma _+}\times \cat C^\Sigma\to \cat C^{\Sigma _+}: (\op X,\op Y)\mapsto \op X\ract \op Y$$
defined by
$$(\op X\ract \op Y)(n)=\coprod_{\substack{ k\geq 1\\ \vec\imath\in I_{k,n}}} \op X(k)\underset {\Sigma _{k-1}}\otimes \big(Y(i_1)\otimes \cdots\otimes Y(i_k)\big)\underset {\Sigma _{\vec\imath-\vec e_1}}  \otimes R[\Sigma _{n-1}].$$
Here, $\Sigma _{k-1}$ acts on $Y(i_1)\otimes \cdots\otimes Y(i_k)$ by permuting $Y(i_2)\otimes \cdots\otimes Y(i_k)$, while $\vec\imath-\vec e_1=(i_1-1, i_2,...,i_k)$.
Furthermore, there is a natural isomorphism $\op X\ract (\op Y\diamond \op Z)\cong (\op X\ract\op Y)\ract \op Z$ for all shifted symmetric sequences $\op X$ and all symmetric sequences $\op Y$ and $\op Z$.  

Let $\op P$ be an operad, with multiplication map $\gamma: \op P\diamond \op P\to \op P$, and let $\op X$ be a shifted symmetric sequence.  We say that $\op X$ is a \emph {shifted right $\op P$-module} if there is a morphism of shifted symmetric sequences $\rho: \op X\ract \op P\to  \op X$ such that 
$$\rho (\rho\ract 1)=\rho(1\ract \gamma):(\op X\ract\op P)\ract P\cong \op X\ract (\op P\diamond \op P)\to \op X.$$ 
A morphism of shifted right $\op P$-modules is a morphism of the underlying shifted symmetric sequences that commutes with the right action maps. We write $\cat{Mod} ^+_{\op P}$ for the category of shifted right $\op P$-modules and their morphisms. Given a shifted symmetric sequence $\op X$ that is a shifted right $\op P$-module and a symmetric sequence $\op Y$ that is a left $\op P$-module (in the usual sense), we define $\op X\underset {\op P} \ract \op Y$ to be the obvious coequalizer.  

\begin{defn}\label{defn:comod-sh} Let $\theta \in \Coalg AF (C,C')$, inducing $\ind (\theta)\in \aalg(\Om C,\Om C')$ and therefore the structure of a right $\Om C$-module on $\Om C'$.    
Suppose that $M$ is a right $C$-comodule and $M'$ 
is a right $C'$-comodule. A  map 
of right $\Om C$-modules
$$h:M\otimes _{t_{\Om}}\Om C\to M'\otimes _{t_{\Om}}\Om C',$$ 
is a \emph {comodule map up to strong homotopy with respect to $\theta$} from $M$ to $M'$.

 Abusing terminology somewhat, we say  that a chain map $g:M\to M'$ is a  \emph{comodule map up to strong homotopy} if there is such an $h$ satisfying 
 $$h(x\otimes 1)-g(x)\otimes 1\in M'_{<\deg x}\otimes \Om C'$$
 for all $x\in M$.   
\end{defn}

A ``module'' version of Theorem \ref{thm:induction} holds for right comodules.  The proof proceeds by straightforward generalization of the absolute case.  Before stating the theorem, we remark that if $C$ is a coassociative chain coalgebra and $M$ is a right $C$-comodule, then $\op T_r(M,C)$ is naturally a shifted right $\op A$-module.  Observe that there is an isomorphism of graded $R$-modules
$$\op T_r(M,C)(k)\underset {\Sigma _{k-1}}  \otimes \op A(i_1)\otimes \cdots\otimes \op A(i_k)\cong \big(M\otimes \op A(i_1)\big)\otimes \big(\op T(C)(k-1)\underset {\Sigma _{k-1}}  \otimes\op A(i_2)\otimes \cdots\otimes \op A(i_k)\big).$$
To define a shifted right $\op A$-module structure on $\op T_r(M,C)$, we use the right $\op A$-module structure on $\op T (C)$ coming from its coalgebra structure, as well as the fact that for all $m\geq 1$, the comodule map $\nu: M\to M\otimes C$ induces a $\Sigma _{m-1}$-equivariant map 
$$\nu ^{(m)}:M\otimes \op A(m)\to  M\otimes C^{\otimes m-1}:(x\otimes \delta ^{(m)})\mapsto (Id_M\otimes \Delta ^{(m-2)})\nu (x).$$ 
In other words, $\op T_r$ induces a functor from the category of pairs $(M,C)$, where  $C$ is a coassociative coalgebra and $M$ is a right $C$-comodule, to the category of shifted right $\op A$-modules.

We can now state the ``module'' version of Theorem \ref{thm:induction}.

\begin{prop}\label{prop:rel-induction}  Let $\theta \in \Coalg AF (C,C')$, inducing $\ind (\theta)\in \aalg(\Om C,\Om C')$.  If $M$ is a right $C$-comodule and $M'$ 
is a right $C'$-comodule, then there is a natural bijection 
$$\ind ^+: \cat{Mod}^+_{\op A}\big(\op T_r(M,C)\aract \op F, \op T_r(M',C')\big)\to \cat{Mod}_{\Om C}\big (M\otimes _{t_{\Om}}\Om C,M'\otimes _{t_{\Om}}\Om C'\big)$$
specified by
$$\ind^+ (\omega)(x)=\sum_{k\geq 1}\big (Id_{M'}\otimes (\si)^{\otimes k-1}\big)\omega (x\otimes z_{k-1})$$
for all $x\in M$.
\end{prop}

To see why the formula above makes sense, note that 
$$\omega (x\otimes z_{k-1})\in \op T_r(M',C')(k)=M'\otimes (C')^{\otimes k-1}.$$  
Furthermore, since $M\otimes _{t_{\Om}}\Om C$ is a free right $\Om C$-module, the specification in Proposition \ref{prop:rel-induction} suffices to imply that $\ind^+ (\omega)(x\otimes v)=\ind^+(\omega)(x)\cdot \ind(\theta)(v)$ for all $x\in M$ and $v\in \Om C$. 

A morphism $\omega\in  \cat{Mod}^+_{\op A}\big(\op T_r(M,C)\aract \op F, \op T_r(M',C')\big)$ gives rise to a family 
\begin{equation}
\mathfrak F^+(\omega)=\{ \omega _k=\omega (-\otimes z_k): M\to M'\otimes (C')^{\otimes k}\mid  \deg \omega _k= k, k\geq 0\}. 
\end{equation}
Specifying a family $\mathfrak F^+(\omega)$ is equivalent to specifying a morphism of shifted symmetric sequences $\op L(M)\aract \op F\to  \op T_r(M',C')$.

Maps of left comodules up to strong homotopy are defined analogously.  A version of the proposition above,  expressed  in terms of a functor $\op T_\ell$, holds for left comodules as well.

Comodule maps up to strong homotopy are interesting because of their role in extending the linear naturality of $\cotor$ (cf.,  (\ref{eqn:cotor})), first established by Gugenheim and Munkholm in  \cite {GM} (dual of Theorem 3.5).   This extended naturality can be expressed as follows in the language we have developed above.

\begin{thm}\label{thm:gm-cotor} \cite {GM} Let $\theta \in \Coalg AF (C,C')$, where $C$ and $C'$ are simply connected chain coalgebras.  Let 
$g:M\to M'$ and $h:N\to N'$ be maps of right and left comodules, respectively, up to strong homotopy 
with respect to $\theta$, where $M$ and $N$ are $C$-comodules, and $M'$ and $N'$ are $C'$-comodules.  Then there is a natural induced morphism of graded $R$-modules
$$\cotor ^{\theta}(g,h):\cotor ^{C}(M,N)\to \cotor ^{C'}(M',N').$$
Furthermore if all the underlying graded modules are $R$-flat 
and $\theta (-\otimes z_0)$, $g$, and $h$ are all quasi-isomorphisms, then $\cotor ^{\theta}(g,h)$ 
is an isomorphism.
\end{thm}

We sketch a proof of Theorem \ref{thm:gm-cotor}, based on  Proposition \ref{prop:rel-induction}.   Let   
$$\xi :\op T_r(M,C)\aract \op F\to \op T(M',C')\quad\text{and}\quad\zeta :\op T_r(N,C)\aract \op F\to \op T(N',C')$$
be the morphisms of shifted $\op A$-modules of chain complexes, corresponding to $g$ and $h$. 
Thus, under the hypotheses of the theorem, we can set
$$\cotor ^{\theta}(g,h)=\H \big(\ind ^+(\xi)\otimes _{\Om C}\ind^+(\zeta)\big):\cotor ^{C}(M,N)\to \cotor 
^{C'}(M',N').$$
A standard spectral sequence argument then shows that $\cotor ^{\theta }(g,h)$ is 
an isomorphism if all modules are $R$-flat and if $\H \big(\theta (-\otimes z_0)\big)$, $\H g$, and $\H h$ are all isomorphisms.

We devote the remainder of this section to establishing a framework in which to state and prove a multiplicative version of Theorem \ref{thm:gm-cotor}.  Recall that if $H$ is a Hopf algebra, then an algebra $B$ 
is an \emph {$H$-comodule algebra} if it is an $H$-comodule and the comodule 
structure maps are algebra maps.   Furthermore, as seen in Corollary \ref{cor:cobar-alg}, a right $H$-comodule algebra $B$ naturally gives rise to a chain algebra, $B\otimes _{t_\Om } \Om H$.

We can now enlarge the category $\cat {CA}^{r}$  by weakening the definition of morphisms, in analogy with the passage from the category of chain coalgebras to the category $\cat {DCSH}$.

\begin{defn}\label{defn:cash} Let $\cat {CASH}$ be the category specified as follows.
\begin{enumerate}
\item Objects are pairs $(H;B)$, where $H$ is a chain Hopf algebra and $B$ is a right $H$-comodule algebra.
\item A morphism from an object $(H;B)$ to an object $(H';B')$ is a pair $(\theta;\gamma)$, where 
$$\theta \in (\op A,\op F)\text{-}\cat {PsHopf}(H,H')$$
and 
$$\gamma:B\otimes _{t_{\Om}}\Om H\to B'\otimes _{t_{\Om}}\Om H'$$
is a morphism of both chain algebras and $\Om H$-modules, where the right $\Om H$-module structure on $B'\otimes _{t_{\Om}}\Om H'$ is given by the algebra morphism $\ind(\theta):\Om H\to \Om H'$.
\end{enumerate}
Composition and identities are defined in the obvious manner.  The morphisms in $\cat {CASH}$ are called \emph{comodule-algebra maps up to strong homotopy}.

Given a morphism $(\theta;\gamma):(H;B)\to (H';B')$ in $\cat {CASH}$, let $\gamma_{0}$ denote the composite 
$$B\hookrightarrow B\otimes _{t_{\Om}}\Om H\xrightarrow \gamma B'\otimes _{t_{\Om}}\Om H'\xrightarrow{\pi}B',$$
where $\pi$ denotes the obvious projection.  We say that a chain map $g:B\to B'$ is a \emph{CASH map} if there is a  morphism $(\theta;\gamma):(H;B)\to (H';B')$ in $\cat {CASH}$ such that $\gamma_{0}=g$.
\end{defn}

\begin{rmk} Corollary \ref{cor:natl-ca} implies that the category $\cat {CA}$ embeds into $\cat{CASH}$ as a wide, but not necessarily full, subcategory.
\end{rmk}

The following relative version of Proposition \ref{prop:extend-map} is a crucial tool for construction of CASH maps.  The proof proceeds by direct, but somewhat cumbersome, generalization of Proposition \ref{prop:free-ext}, the details of which we spare the reader.  We use here the shifted derivations of Definition \ref{defn:shift-deriv}.

\begin{prop}\label{prop:rel-free-ext}Fix chain Hopf algebras $H$ and $H'$ and 
$$\theta \in (\op A,\op F)\text{-}\cat {PsHopf}(H,H')$$ 
with $\mathfrak F(\theta)=\{ \theta _k\mid k\geq 1\}$.  Let $B$ be a right $H$-comodule algebra and $B'$ a 
right $H'$-comodule algebra  such that $B$ is free as an algebra on an $R$-free 
graded module $V$.  Then any family of morphisms of graded $R$-modules
$$\Xi=\{\xi_{k}:V\to B'\otimes (H')^{\otimes k}\mid \deg \xi_{k}=k, k\geq 0\}$$
naturally induces a unique morphism of shifted right $\op A$-modules of graded $R$-modules
$$\widehat\xi:\op T_r(B,H)\aract \op F\to \op T_r(B',H')$$
such that $\widehat\xi (v\otimes z_k)=\xi _k(v)$ for all $v\in V$ and such that  $$\ind^+(\,\widehat\xi\,):B\otimes _{t_{\Om}}\Om H\to B'\otimes _{t_{\Om}}\Om H'$$ is a map of graded algebras and of $\Om H$-modules.

If, furthermore, for all $k\geq 0$ and for all $v\in V$,
\begin{equation}
\op D(d_{B'},d_{H'})_{k+1}\xi _k(v) -\op D(\overline\nu', \overline\Delta ')_{k}\xi _{k-1}(v)=\widehat\xi_k(dv)-\sum _{i+j=k}(\widehat\xi _i\otimes \theta _j)\overline \nu (v), 
\end{equation}
where $\mathfrak F^+(\widehat \xi )=\{\widehat\xi _k\mid k\geq 0\}$, then $\widehat\xi$ is a differential map.
\end{prop}

The next proposition, which explains how to construct a CASH map as a sort of coproduct of CASH maps when the underlying algebras of the sources are free, is essential to the proof in section 
\ref{sec:homol-htpyfib} that our algebraic ``homotopy fiber" has the right homology.  Before stating the proposition, we need one observation about coproducts and tensor products of comodule algebras.

\begin{rmk} Suppose that $A$ and $A'$ are right $H$-comodule algebras with coaction maps $\nu: A\to A\otimes H$ and $\nu':A'\to A'\otimes H$. Let  $A\coprod A'$ denote the coproduct of $A$ and $A'$ in the category of chain algebras. Since $\nu $ and $\nu'$ are algebra maps, they together induce an algebra map
$$\nu'':A\coprod A'\to (A\otimes H)\coprod (A'\otimes H)\to  (A\coprod A')\otimes H,$$
which satisfies the axioms of a coaction because $\nu $ and $\nu'$ do.
In other words, the algebra coproduct of $H$-comodule algebras is naturally an $H$-comodule algebra.

It is easy to check that the tensor product $A\otimes A'$, with its usual algebra structure, also admits a natural $H$-coaction
$$A\otimes A'\xrightarrow{\nu \otimes \nu'} (A\otimes H)\otimes (A'\otimes H) \xrightarrow{\cong}(A\otimes A')\otimes (H\otimes H)\xrightarrow{Id\otimes \mu}(A\otimes A')\otimes H$$
that is an algebra map, where $\mu$ denotes the multiplication in $H$. 
\end{rmk}

\begin{prop}\label{prop:cash-coproduct} Let $(\theta;\gamma):(H;A)\to (K;B)$ and $(\theta;\gamma'):(H;A')\to (K;B')$ be morphisms in $\cat{CASH}$.  Endow $A\coprod A'$ and $B\otimes B'$ with their natural $H$-comodule and $K$-comodule algebra structures.  If the algebras underlying $A$ and $A'$ are free on free graded $R$-modules $V$ and $V'$, respectively, then there exists a morphism 
$$(\theta; \gamma''):(H;A\coprod A')\to (K;B\otimes B')$$
in $\cat {CASH}$ such that $\gamma''_{0}(v)=\gamma_{0}(v)\otimes 1$ and $\gamma''(v')=1\otimes \gamma'_{0}(v')$ for all $v\in V$ and $v'\in V'$.
\end{prop} 
 
In the situation of the proposition above, we write
$$\gamma''_{0}=\gamma_{0}\star\gamma'_{0}.$$

\begin{proof} Let 
$$\xi :\op T_r(A,H)\aract \op F\to \op T(B,K)\quad\text{and}\quad\xi' :\op T_r(A',H)\aract \op F\to \op T(B',K)$$
be the morphisms of shifted $\op A$-modules of chain complexes, corresponding to $g$ and $g'$, with corresponding families $\mathfrak F^+(\xi )=\{\xi _m\mid m\geq 0\}$ and $\mathfrak F^+(\xi' )=\{\xi' _m\mid m\geq 0\}$.

Define a family of linear maps
$$\Xi''=\{\xi ''_m: V\oplus V'\to (B\otimes B')\otimes K^{\otimes m}\mid m\geq 0\}$$
by $\xi_m ''(v)=\iota _1\circ \xi_m(v)$ for all $v\in V$ and $\xi _m''(v')=\iota _2\circ \xi'_m(v')$ for all $v'\in V'$, where
$$\iota _1:B\otimes K^{\otimes m}\to (B\otimes B')\otimes K^{\otimes m}:x\otimes y_1\otimes \cdots \otimes y_m\mapsto x\otimes 1\otimes y_1\otimes \cdots \otimes y_m$$
and 
$$\iota _2:B'\otimes K^{\otimes m}\to (B\otimes B')\otimes K^{\otimes m}:x'\otimes y_1\otimes \cdots \otimes y_m\mapsto 1\otimes x'\otimes y_1\otimes \cdots \otimes y_m.$$
Note that the algebra map induced by the linear map $\xi '' _0$ is indeed $\gamma_{0}\star \gamma_{0} '$.
 
It is easy to check that the family $\Xi''$ satisfies the conditions of Proposition \ref{prop:rel-free-ext}, since $\xi$ and $\xi'$ are morphisms of shifted $\op A$-modules of chain complexes.
 \end{proof}

Inspired by the sketch of the proof of Theorem \ref{thm:gm-cotor}, we can easily verify the following result, establishing extended \emph{multiplicative} naturality of $\cotor$, generalizing both Corollary \ref{cor:cotor-alg} and Theorem \ref{thm:gm-cotor}.

\begin{prop} \label{prop:mult-cotor} Let $\cat{grA}$ denote the category of graded algebras over $R$. The functor $\cotor^{(-)}(-;R)$ of Corollary \ref{cor:cotor-alg} extends to a functor
$$\cotor^{(-)} (-;R):\cat {CASH}\to \cat{grA}.$$
\end{prop}

\begin{proof} Let $H$ and $H'$ be  simply-connected chain Hopf algebras, and let 
$$\theta \in (\op A,\op F)\text{-}\cat {PsHopf}(H,H').$$ 
Let $g:M\to M'$ be a CASH map with respect to $\theta$, 
where $M$ is a right $H$-comodule, and $M'$ is a right $H'$-comodule. Recall the  graded algebra structure on $\cotor ^H(M,R)$ from Corollary \ref{cor:cotor-alg} and its proof.  

 Let   
$$\xi :\op T_r(M,H)\aract \op F\to \op T(M',H')$$
be the morphism of shifted right $\op A$-modules corresponding to $g$. Since $g$ is a CASH map, $\ind^+ (\xi )$ is a chain algebra map and hence $\cotor ^\theta(g,Id)$ is a map of graded algebras. 
\end{proof}





\section {Path objects and homotopy fibers in $\mathbf F$}\label{sec:path-loop}

In this section we define  a functor
$$\xymatrix@1{\pl: \shc \ar [r]&\mathbf H},$$
called the \emph {path-loop functor}. For every $(C,\Psi)$ in $\shc$, there is a natural surjection of chain Hopf algebras $\pl (C,\Psi)\to  \tom( C,\Psi)$.  The definition of $\pl$ is the first step towards building a particularly nice chain algebra from which we can compute $\cotor ^H(R, R)$, when $H$ is a  chain Hopf algebra endowed with an \emph {Alexander-Whitney model}, as defined in section \ref{sec:homol-htpyfib}. As we explain in section \ref{sec:topology}, the terminology chosen is justified by the fact that  the homotopy fiber of the natural surjection $\pl \tC (K)\to \tom \tC(K)$ is indeed a model for $G^2K$, where $\tC$ denotes the functor of Theorem \ref{thm:loopmodel}.

We begin by more general considerations.  Given any graded module $X$, let $\bX$ denote $\si X$, and let $\bx$ denote an element $\si x$.  Let $\sigma : X\to X\oplus \bX$ be defined by $\sigma (x)=\bx$, and let $\iota : X\to X\oplus \bX$ denote the inclusion.  

Let $(X,d)$ be any chain complex.  The \emph {based-path object} on $(X,d)$, denoted $\mathfrak P(X,d)$, is defined to be the acyclic chain complex $(X\oplus \bX, \td)$, where $\td x=\iota dx-\bx$ and $\td \bx=-\overline{dx}$, i.e., 
\begin{equation}\label{eqn:cone-diffl}
\td\iota =\iota d-\sigma\quad\text {and}\quad \td\sigma=-\sigma d. 
\end{equation}
There is an obvious factorization
$$\xymatrix{
0\ar[rr]\ar[dr]_{\simeq}&&(X,d)\\
&\mathfrak P(X,d)\ar[ur]^\pi}$$
where $\pi$ denotes the obvious projection,
justifying the name we have given to the chain complex $\mathfrak P(X,d)$.

The based-path construction is clearly natural, i.e., there is a functor 
$$\mathfrak P: \cat{Ch}_R\to \cat{Ch}_R.$$  
Furthermore, the functor $\mathfrak P$ is comonoidal, where the natural transformation $\mathfrak j:\mathfrak P(-\otimes -)\to  \mathfrak P(-)\otimes \mathfrak P(-)$ is defined for chain complexes $X$ and $Y$ to be the injection
$$\mathfrak j_{X,Y}:(X\otimes Y)\oplus\overline {(X\otimes Y)}\cong (X\otimes Y)\oplus (\bX\otimes Y)\oplus (X\otimes \overline Y)\hookrightarrow (X\oplus \bX)\otimes (Y\oplus \overline Y).$$
In particular, if $(X,d,\Delta)$ is a coassociative coalgebra, then $\big(\mathfrak P(X,d), \tD\big)$ is also a coassociative coalgebra, where $\tD=\mathfrak j_{X,X}\mathfrak P(\Delta)$.  Note that the comultiplication on $\mathfrak P(X)$  is specified by $\tD\iota=(\iota\otimes \iota)\Delta$ and $\tD \sigma=(\sigma \otimes \iota +\iota \otimes \sigma)\Delta$ and that the projection map $\pi: \mathfrak P(X)\to X$ is a morphism of coalgebras.

The morphisms of graded $R$-modules $\iota$ and  $\sigma$ induce a morphism of symmetric sequences $\op D_{\iota, \iota}(\sigma):\op T(X)\to \op T\big(\mathfrak P(X)\big)$ that is of degree $-1$ in each level (cf., Definition \ref{defn:derivation}), while the differentials $d$ and $\td$ induce $\op D(d):\op T(X)\to \op T(X)$ and $\op D(\td):\op T\big(\mathfrak P(X)\big)\to \op T\big (\mathfrak P(X)\big)$.  It is a matter of straightforward calculation to show that (\ref{eqn:cone-diffl}) implies that
\begin{equation}\label{eqn:cone-deriv}
\op D(\td)\op D_{\iota, \iota}(\sigma)=-\op D_{\iota, \iota}(\sigma)\op D(d)\text{ and } \op D(\td)\op T(\iota)=\op T(\iota)\op D(d)-\op D_{\iota, \iota}(\sigma).
\end{equation}

\begin{prop}\label{prop:ext-morph} Let $C$ and $C'$ be coassociative chain coalgebras. Any morphism 
$$\theta :\op T (C)\acirc \op F\to  \op T (C')$$
of right $\op A$-modules of chain complexes lifts naturally to a morphism 
$$\widetilde\theta :\op T \big(\mathfrak P(C)\big)\acirc \op F\to  \op T\big(\mathfrak P (C')\big),$$
i.e., there is a commuting diagram of morphisms in $\Coalg AF$
$$\xymatrix{\mathfrak P(C)\ar [r]^{\widetilde \theta }\ar [d]_{\pi}&\mathfrak P(C')\ar [d]_{\pi}\\
C\ar[r]^\theta &C'.}$$

Here, the coalgebra morphisms denoted $\pi$ are considered as morphisms in $\Coalg AF$ via the inclusion functor.
\end{prop}

\begin{proof} We remark first that since $\theta$ is a differential map, the following equality holds.
\begin{equation}\label{eqn:diffl-morph}
\op D(d')\theta =\theta \big(\op D(d)\acirc 1+1\acirc\del _{\op F}\big)
\end{equation}
Here, the composition rule applied is that of morphisms of right $\op A$-modules.
 
We now define a morphism of symmetric sequences of graded $R$-modules
$$\theta ':\op L(C\oplus \bC)\diamond \op S\to  \op T(C'\oplus \bC)$$
by $\theta' (\iota (c)\otimes z_{m-1})=\iota ^{\otimes m} \theta (c\otimes z_{m-1})$ for all $c\in C$, since $\widetilde\theta$ should extend $\theta$, and $\theta '(\bc\otimes z_{m-1})=\op D_{\iota, \iota}(\sigma)\theta (c\otimes z_{m-1})$.  In other words, 
\begin{equation}\label{eqn:theta-iota}
\theta'(\iota\diamond 1)=\op T(\iota)\theta\text{ and }\theta'(\sigma \diamond 1)=\op D_{\iota, \iota}(\sigma)\theta. 
\end{equation}

Applying Lemma \ref{lem:inducedT},  we obtain another morphism of symmetric sequences
$$\theta'':\op T(C\oplus \bC)\diamond \op S\to \op T(C'\oplus \bC'),$$
defined for all $k$ and for all $w_1,..., w_k\in C\oplus \bC$ by
$$\theta''\big( (w_1\otimes \cdots\otimes w_k)\otimes (z_{n_1-1}\otimes \cdots \otimes z_{n_k-1})\big)=\pm \theta '(w_1\otimes z_{n_1-1})\otimes \cdots \otimes \theta'(w_k\otimes z_{n_k-1}),$$
where the sign is determined by the Koszul rule.

Now use the right $\op A$-module structure of $\op T\big(\mathfrak P(C')\big)$ to extend $\theta''$ to a morphism
$$\widetilde\theta: \op T (C\oplus \bC)\acirc \op F\cong \op T(C\oplus \bC)\diamond \op S\diamond \op A\to \op T(C'\oplus \bC')$$
of right $\op A$-modules of graded $R$-modules. As an easy consequence of (\ref{eqn:theta-iota}), we have that
$$\ttheta \big(\op T(\iota)\acirc 1)=\op T(\iota)\theta\text{ and }\ttheta \big(\op D_{\iota, \iota}(\sigma)\big)=\op D_{\iota, \iota}(\sigma)\theta.$$  

To complete the proof, we need to verify that $\widetilde\theta$ is differential, i.e., that 
$$\widetilde\theta (\op D(d)\acirc 1 + 1\acirc \del _{\op F})=\op D(\td')\widetilde\theta.$$
It is enough to prove that the two sides of the equation are equal when precomposed (as maps of right $\op A$-modules) with either $\op L(\iota)\diamond 1:\op L(C)\diamond\op F \to \op T\big (\mathfrak P(C)\big)\acirc \op F$ or $\op L(\sigma)\diamond 1:\op L(C)\diamond\op F \to \op T\big (\mathfrak P(C)\big)\acirc \op F$.

Observe that 
\begin{align*}
\op D(\td')\ttheta (\op L(\iota)\diamond 1)&=\op D(\td')\op T(\iota)\theta\\
&=\op D_{\iota, \iota}(\td\iota)\theta\\
&=\op D_{\iota, \iota} (\iota d-\sigma)\theta\\
&=\big(\op T(\iota)\op D(d')-\op D_{\iota, \iota} (\sigma)\big)\theta,
\end{align*}
while
\begin{align*}
\ttheta\big(\op D(\td)\acirc 1\big)(\op L(\iota)\diamond 1)&=\ttheta\big(\op T(\iota)\op D(d)\acirc 1-\op D_{\iota, \iota}(\sigma)\acirc 1\big)\\
&=\op T(\iota)\theta \op D(d)-\ttheta\big(\op D_{\iota, \iota}(\sigma)\acirc 1\big).
\end{align*} 
Thus, 
\begin{align*}
\big(\op D(\td')\ttheta-\ttheta(\op D(\td)\acirc 1)\big)(\op L(\iota)\diamond 1)&=\op T(\iota)\big (\op D(d')\theta-\theta \op D(d)\big)\\
&=\op T(\iota)\theta (1\acirc \del_{\op F})\\
&=\ttheta (1\acirc\del _{\op F})(\op L(\iota)\diamond 1),
\end{align*}
by equation (\ref{eqn:diffl-morph}).

Similarly, 
$$\op D(\td')\ttheta (\op L(\sigma)\diamond 1)=\op D(\td')\op D_{\iota, \iota}(\sigma)\theta=-\op D_{\iota, \iota}(\sigma)\op D(d')\theta,$$
and
$$\ttheta\big(\op D(\td)\acirc 1\big)(\op L(\sigma)\diamond 1)=-\ttheta(\op D_{\iota, \iota}(\sigma)\op D(d)\acirc 1)=-(\op D_{\iota, \iota}(\sigma)\theta(\op D(d)\acirc 1).$$
Thus, 
\begin{align*}
\big(\op D(\td')\ttheta-\ttheta(\op D(\td)\acirc 1)\big)(\op L(\sigma)\diamond 1)&=-(\op D_{\iota, \iota}(\sigma)\theta(1\acirc \del_{\op F})\\
&=-\ttheta (\op D_{\iota, \iota}(\sigma)\acirc 1)(1\acirc\del _{\op F})\\
&=\ttheta (1\acirc \del _{\op F})(\op L(\sigma)\diamond 1),
\end{align*}
again by (\ref{eqn:diffl-morph}).
 \end{proof}

Recall that $\mathfrak I_{\op F}:\acoalg\to \Coalg AF$ denotes the ``inclusion" functor (\ref{eqn:emb-fat}).

\begin{cor}\label{cor:aw-cone}  If $(C, \Psi)$ is an object in $\shc$, then the based-path object $\mathfrak P(C)$ admits a natural Alexander-Whitney coalgebra structure map $\widetilde \Psi$, extending $\Psi$.  Furthermore, the morphism of right $\op A$-modules 
$$\mathfrak I_{\op F}(\pi )=\op T(\pi )\acirc \varepsilon :\op T\big(\mathfrak P(C)\big)\acirc \op F\to \op T(C)$$
induced by the natural projection map of chain complexes $\pi:\mathfrak P(C)\to C$ is a morphism in $\cat F$, i.e., 
$$\Psi \mathfrak I_{\op F}(\pi)=\big(\mathfrak I_{\op F}(\pi)\curlywedge \mathfrak I_{\op F}(\pi)\big)\widetilde\Psi,$$
where the composition is calculated in $\Coalg AF$.
\end{cor}

\begin{proof} By Proposition \ref{prop:ext-morph} , the morphism of right $\op A$-modules
$$\Psi :\op T(C)\acirc \op F\to \op T(C\otimes C)$$
gives rise naturally to 
$$\widetilde \Psi:\op T\big(\mathfrak P(C)\big)\acirc \op F\to \op T\big(\mathfrak P(C\otimes C)\big).$$
Since $\mathfrak P(C\otimes C)$ injects into $\mathfrak P(C)\otimes \mathfrak P(C)$, we can look at $\widetilde \Psi $ as a morphism with target $\op T\big(\mathfrak P(C)\otimes \mathfrak P(C)\big)$.

To complete the proof that $\big(\mathfrak P(C), \widetilde \Psi)$ is an Alexander-Whitney coalgebra, we need to check that $q\ind (\widetilde \Psi)$ is coassociative.  By naturality, however, this follows immediately from the coassociativity of $q\ind(\Psi)$.

Verification that $\mathfrak I_{\op F}(\pi)$ is a morphism in $\cat F$ is trivial.
 \end{proof}

\begin{prop}\label{prop:awmorph-cone} Let $\theta :(C,\Psi)\to (C',\Psi')$ be  a morphism in $\cat F$.  Then $$\widetilde \theta:\big(\mathfrak P (C),\widetilde \Psi)\to \big(\mathfrak P(C'), \widetilde \Psi '\big)$$ is also a morphism in $\cat F$, i.e., $(\widetilde \theta\curlywedge \widetilde \theta)\widetilde\Psi=\widetilde\Psi'\widetilde\theta$, where the composition is performed in $\Coalg AF$.\end{prop}

\begin{proof} It suffices to check the desired equality holds when precomposed (as morphisms of right $\op A$-modules) with either $\op L(\iota)\diamond 1$ or $\op L(\sigma)\diamond 1$.  To distinguish between composition as morphisms of right $\op A$-modules and as morphisms in $\Coalg AF$, we denote the first by simple concatentation of symbols and the second by $\circ$.

The definition of $\tPsi$ given in the proofs of Proposition \ref{prop:ext-morph} and Corollary \ref{cor:aw-cone} implies that
\begin{align*}
(\ttheta\twedge\ttheta)\circ \tPsi(\op L(\iota)\diamond 1)&=(\ttheta\twedge\ttheta)\circ \op T(\iota)\Psi\\
&=\big(\ttheta \op T(\iota)\twedge\ttheta \op T(\iota)\big) \circ \Psi\\
&=\op T(\iota)(\theta\twedge \theta)\circ \Psi\\
&=\op T(\iota)\Psi '\circ \theta\qquad\text{since $\theta $ is a morphism in $\cat F$}\\
&=\tPsi '\circ \op T(\iota)\theta\\
&=\tPsi '\ttheta (\op L(\iota)\diamond 1).
\end{align*}
On the other hand,
\begin{align*}
(\ttheta\twedge\ttheta)\circ \tPsi(\op L(\sigma)\diamond 1)&=(\ttheta\twedge\ttheta)\circ \op D_{\iota,\iota}(\sigma)\Psi\\
&=\big(\ttheta \op D_{\iota, \iota}(\sigma)\twedge\ttheta \op T(\iota)+\ttheta \op T(\iota)\twedge\ttheta  \op D_{\iota, \iota}(\sigma)\big) \circ \Psi\\
&= \op D_{\iota, \iota}(\sigma)(\theta\twedge \theta)\circ \Psi\\
&= \op D_{\iota, \iota}(\sigma)\Psi '\circ \theta\qquad\text{since $\theta $ is a morphism in $\cat F$}\\
&=\tPsi '\circ  \op D_{\iota, \iota}(\sigma)\theta\\
&=\tPsi '\ttheta (\op L(\sigma)\diamond 1).
\end{align*}
 \end{proof}

Corollary \ref{cor:aw-cone} and Proposition \ref{prop:awmorph-cone} imply that the following definition makes sense.

\begin{defn} \label{defn:based-path} The \emph {based-path functor} $\widetilde{\mathfrak P}:\mathbf F\to \mathbf F$ is defined on objects by $\widetilde{\mathfrak P}(C,\Psi)=\big(\mathfrak P(C),\widetilde \Psi\big)$ and on morphisms by $\widetilde{\mathfrak P}(\theta)=\widetilde \theta$.
\end{defn}

The second part of Corollary \ref{cor:aw-cone} implies that $\widetilde{\mathfrak P}$ is augmented: the projection $\pi$ serves as a natural transformation $\pi:\widetilde{\mathfrak P}\to Id_{\cat F}$.

For the constructions in the following sections, we need a relative version of the path functor, i.e., a notion of  homotopy fiber in $\cat F$.  We consider first the notion of homotopy fiber in $\cat {Ch}_{R}$.  Any morphism of chain complexes $f:X\to Y$ can be factored as
$$\xymatrix{
X\ar[rr]^f\ar [dr]_{\text{incl.}}^\simeq &&Y\\
&X\oplus \mathfrak P(Y)\ar [ur]_{f+\pi}},$$
so that it is reasonable to define the homotopy fiber $\mathfrak {HF}(f)$ to be the pullback
$$\xymatrix{
\mathfrak{HF}(f)\ar [d]\ar [r]&X\oplus \mathfrak P(Y)\ar [d]^{f+\pi}\\
0\ar [r]&Y,}$$
or, more prosaically, $\mathfrak{HF}(f)$ is the kernel of $f+\pi$.

Analogously, to define homotopy fibers in $\cat F$, we first need to show that $\cat F$ admits coproducts.  It is easy to see, however, that if $(C,\Psi)$ and $(C',\Psi')$ are Alexander-Whitney coalgebras, then their coproduct $(C,\Psi)\coprod(C',\Psi')$ in $\cat F$ is $(C\oplus C', \Psi'')$, where
$$\Psi'': \op T\big (C\oplus C'\big)\acirc \op F\to  \op T \left(\big (C\oplus C'\big)^{\otimes 2}\right)$$
is the morphism of $\op A$-modules specified (as in the proofs of Lemma \ref{lem:inducedT} and Proposition \ref{prop:ext-morph}) by 
$$\Psi''(c\otimes z_{k-1})=\Psi(c\otimes z_{k-1})\qquad\text{and}\qquad\Psi''(c'\otimes z_{k-1})=\widetilde\Psi ( c'\otimes z_{k-1})$$
for all $c\in C$, $c'\in C'$ and $k>0$ .

Let $\theta:(C',\Psi')\to (C,\Psi)$ be a morphism in $\shc$.  There is an obvious factorization in $\cat F$
\begin{equation}\label{eqn:can-fact}
\xymatrix{
(C',\Psi ')\ar[rr]^\theta\ar [dr]_{\text{incl.}}^\simeq &&(C,\Psi)\\
&(C',\Psi')\coprod \mathfrak P(C,\Psi)\ar [ur]_{\theta+\pi}}.
\end{equation} 
Note that $\theta+\pi$, seen simply as a map of coalgebras, admits a (nondifferential) section, the coalgebra map 
$$\mathfrak j: C\to C'\oplus(\cbc),$$ 
which is just the natural inclusion.  Moreover, $\Psi ''(\mathfrak j(c)\otimes z_{k-1})=\mathfrak j ^{\otimes k}\Psi (c\otimes z_{k-1})$ for all $k$, so that the induced algebra map  
$$\Om \mathfrak j:\Om C\to \Om \big(C'\oplus (\cbc)\big)$$ 
commutes with the induced comultiplications, i.e., 
$\Om \mathfrak j$ is a Hopf algebra map, which is a (nondifferential) section of $\Om (\theta+\pi)$.

Recall the functor $\tom: \cat F\to \cat H$ from (\ref{eqn:tom}).

\begin{defn} The \emph {algebraic path-loop functor} 
$$\xymatrix{\pl:\shc\ar [r]&\mathbf H}$$
is the composite $\pl =\tom\circ \widetilde{\mathfrak P}$.  The induced comultiplication on $\pl(C,\Psi)$ is denoted $\tpsi$.
\end{defn}

We prefer the notation $\pl$ for this functor, instead of $\tom\circ \widetilde {\mathfrak P}$, as it reminds us that $\pl(C,\Psi)$ plays the role of the \textbf{p}aths on a \textbf{l}oop space.

Observe that $\pl (C)$ is always acyclic, since $\widetilde{\mathfrak P}(C)$ is acyclic and $1$-connected. 

Consider the natural right $\Om C$-comodule structure on $\pl(C)$ given by the coaction
$$\nu=(1\otimes \Om \pi)\tpsi :\Om (\cbc)\to \Om (\cbc )\otimes \Om C.$$ 
It is important for the proof of Theorem \ref{thm:homol-loopfib} to know that $\pl (C)$ is a cofree right $\Om C$-comodule, which is an immediate consequence of the following more general result.

\begin{prop}\label{prop:mm} Let $p:(H',d')\to (H,d)$ be a surjection of connected chain Hopf algebras, which are free as graded $R$-modules.  Let $\nu=(1\otimes p )\Delta ':H'\to H'\otimes H$ denote the right $H$-coaction on $H'$ induced by $p$.  If $p$ admits a (nondifferential) Hopf algebra section $s$,  then
\begin{enumerate}
\item $H'\square_{H}R$ is a sub chain algebra of $H'$, and
\item  $(H',\nu)$ is cofree as a nondifferential $H$-comodule, with cobasis $H'\square _H R$.
\end{enumerate}
\end{prop}

\begin{proof} Let $(B, d_{B})=H'\square _{H}R$. Theorem 4.6 
in \cite {MM}, 
implies directly that 
$H'\square_{H}R$ is a sub chain  algebra of $H'$.  

Consider the linear map 
$$h=\mu (i\otimes s): B\otimes H\xrightarrow{\cong } H',$$
where $\mu $ denotes the multiplication on $H'$ and $i$ is the canonical inclusion. According to Theorem 4.7 in \cite {MM}, $h$ is an isomorphism of both right $H$-comodules and left 
$B$-modules, since the 
underlying Hopf algebra of $H$ is connected, while the underlying 
algebra of $H'$ is a connected $H$-comodule algebra, and 
all graded $R$-modules in question are free. Since $h$ is an isomorphism, we can use it 
to define a differential $\bar d$ on $B\otimes H$ by 
$$\bar d= h^{-1}d'h.$$
Then $h$ becomes an isomorphism of differential right $H$-comodules and left 
$B$-modules, i.e., $H'$ is cofree.  

Observe that for all $x\in B$, 
\begin{align*}
\bar d (x\otimes 1)&= h^{-1}d'(i(x))\\
&=h^{-1}i(d_{B}x)\\
&=d_{B}x\otimes 1,\end{align*}
i.e., the restriction of $\bar d$ to $B\otimes 1$ is simply $d_{B}\otimes 1$.  In other words, the inclusion of $B$ into $H'$ is a differential map.
 \end{proof}
 
 \begin{cor}\label{cor:pl-cofree}  The path-loop construction $\mathfrak{PL}(C,\Psi)$ on any Alexander-Whitney coalgebra $(C,\Psi)$ is cofree over $\tom(C,\Psi)$, with cobasis $ \mathfrak{PL}(C,\Psi)\square_{\tom(C,\Psi)}R,$ which is a sub chain algebra of $\mathfrak{PL}(C,\Psi)$.
 
 More generally, for any $\theta\in \cat F\big((C',\Psi'), (C,\Psi)\big)$, the coalgebra $\tom\big((C',\Psi')\coprod\mathfrak P(C,\Psi)\big)$ is cofree over $\tom(C,\Psi)$, with cobasis  $\tom\big((C',\Psi')\coprod\mathfrak P(C,\Psi)\big)\square _{\tom(C,\Psi)}R$, which is a sub chain algebra of $\tom\big((C',\Psi')\coprod\mathfrak P(C,\Psi)\big)$.
 \end{cor}
 
 \begin{proof} Since $\tom(\theta+\pi):\tom\big((C',\Psi')\coprod\mathfrak P(C,\Psi)\big)\to \tom (C,\Psi)$ admits a (nondifferential) Hopf algebra section $\Om \mathfrak j$, we can apply Proposition \ref{prop:mm}.
 \end{proof}

To conclude this section we analyze more precisely the nature of 
$\tpsi$, the comultiplication on $\pl(C)$, and $\nu$, the induced $\tom C$-coaction.  Let $\iota $ denote the natural (nondifferential) section $\Om C\hookrightarrow\Om (\cbc)$ of $\Om \pi$.  Let 
$$\kappa: \Om C\to  \Om (\cbc)$$
denote the $(\iota, \iota)$-derivation of degree $-1$ specified by $\kappa (\si 
c)=-\si \bc$, i.e.,
$\kappa \mu =\mu (\kappa\otimes\iota +\iota\otimes \kappa)$.

\begin{lem}\label{lem:kappa} The derivation $\kappa $ satisfies the following properties.
\begin{enumerate}
\item $\kappa$ is a differential map of degree $-1$, i.e., $\kappa d_{\Om}=-\td_{\Om}\kappa$.
\item $\kappa $ is a $(\iota, \iota)$-coderivation, i.e., $$\tpsi\kappa =(\kappa\otimes \iota +\iota\otimes \kappa)\psi:\Om 
C\to \Om (\cbc)\otimes \Om (\cbc).$$
\item $\kappa$ is a map of right $\Om C$-comodules, i.e., $\nu\kappa=(\kappa \otimes 1)\psi :\Om C\to \Om 
(\cbc)\otimes \Om C$.
\end{enumerate}
\end{lem}

\begin{proof} (1) Let $c\in C$, and let $c_{i}\otimes c^i$ denote its 
reduced comultiplication. Then, using the definitions of $\td$ and $\tD$ from the beginning of section \ref{sec:path-loop} as well as the definition of the cobar construction differential from the introduction, we obtain
\begin{align*}
\kappa d_{\Om}(\si c)&=\kappa (-\si dc+(-1)^{c_{i}}\cdot \si c_{i}\si c^i)\\
&=\si \overline {dc}+(-1)^{c_{i}}(-\si \bc _{i}\cdot\si c^i-(-1)^{c_{i}}\si c_{i}\cdot\si 
\bc^i)\\
&=\td _{\Om}\si \bc\\
&=-\td _{\Om}\kappa (\si c).
\end{align*}

\medskip
\noindent (2) We have defined $\tpsi $ so that the desired equality 
obviously holds on the generators $\si C_{+}$.  Thus, to establish that 
the equality holds on all of $\Om C$, we must show that 
$$\tpsi\kappa\mu =(\kappa\otimes \iota +\iota\otimes\kappa)\psi 
\mu:\Om C^{\otimes 2}\to \Om (\cbc )^{\otimes 2}$$
where $\mu $ denotes the multiplication map.  

We verify this equality 
by induction on total length of elements in $\Om C^{\otimes 2}$.  By 
definition of $\tpsi$, the equality holds for total length equal to 
$1$.  Suppose that it holds for all elements of $\Om C^{\otimes 2}$ 
of total length less than $n$   

Let $\tau $ denote the usual twisting isomorphism $\tau: A\otimes B\cong 
B\otimes A$. Observe that on $T^{m}\si C_{+}\otimes T^{n-m}\si C_{+}$,
\begin{align*}
\tpsi\kappa\mu&=\tpsi \mu (\kappa\otimes \iota+\iota\otimes \kappa)\\
&=(\mu\otimes\mu)(1\otimes \tau\otimes 1)(\tpsi\otimes\tpsi)(\kappa\otimes \iota+\iota\otimes \kappa)\\
&=(\mu\otimes\mu)(1\otimes \tau\otimes 1)\left( \big((\kappa\otimes \iota 
+\iota\otimes \kappa)\psi\big)\otimes\psi+ \psi\otimes\big((\kappa\otimes \iota 
+\iota\otimes \kappa)\psi\big)\right)\\
&=(\mu\otimes \mu)(\kappa \otimes \iota ^{\otimes 3}+\iota ^{\otimes 
2}\otimes \kappa\otimes\iota+\iota\otimes\kappa\otimes\iota^{\otimes 2}+\iota 
^{\otimes 3}\otimes \kappa)(1\otimes \tau \otimes 1)(\psi\otimes\psi)\\
&=\left(\big(\mu (\kappa\otimes \iota+\iota\otimes 
\kappa)\big)\otimes\mu +\mu\otimes\big(\mu (\kappa\otimes \iota+\iota\otimes 
\kappa)\big)\right)(1\otimes \tau \otimes 1)(\psi\otimes\psi)\\
&=(\kappa \mu\otimes\mu+\mu\otimes\kappa\mu)(1\otimes \tau \otimes 1)(\psi\otimes\psi)\\
&=(\kappa \otimes\iota+\iota\otimes\kappa)\psi \mu.\end{align*}
The induction hypothesis assures that the third equality in this sequence holds.

The equality of part (2) of the lemma therefore holds for all elements of total length 
$n$.

\medskip
\noindent (3) This is an immediate consequence of (2). \end{proof}


\section {Homology of homotopy fibers in $\mathbf H$}\label{sec:homol-htpyfib}

In this section we describe the homology of the homotopy fiber  $\dl (C,\Psi)$ of the path-loop map $\Om\pi:\pl(C,\Psi)\to \tom(C,\Psi)$ on an object $(C,\Psi)$ of $\shc$.  We show in particular that when a chain Hopf algebra $H$ is endowed with an \emph {Alexander-Whitney model} $\theta :\tom (C,\Psi)\to H$ (see below), then 
$$\H \big(\dl (C,\Psi )\big)\cong \cotor ^H(R, R)$$
as graded algebras. We emphasize that this is a true isomorphism and not merely an isomorphism of associated bigraded complexes: there are no extension problems to solve. More generally, we apply the path-loop construction to  building a model for computation of the algebra structure  of $\cotor ^H(H', R)$,  the homology of the homotopy fiber of  a map of chain Hopf algebras $H'\to H$, which endows  $H'$ with the structure of an $H$-comodule algebra.  In section \ref{sec:topology} we show that
 our terminology is fully justified by its application to chain complexes of simplicial sets.

The chain Hopf algebras that we can study by the methods of this paper possess a model of the following sort.  We use here the notion of pseudo $\op A$-Hopf algebras of Definition \ref{defn:pseudo}.

 \begin{defn}\label{defn:awmodel} Let $H$ be a chain Hopf algebra, seen as a pseudo $\op A$-Hopf algebra, via the ``inclusion" functor $\mathfrak I_{\op F}: \mathbf H\to (\op A,\op F)\text{-}\cat {PsHopf}$. An \emph {Alexander-Whitney model} of $H$ consists of an object $(C,\Psi)$ of $\shc$ together with a morphism 
 $$\Theta\in(\op A,\op F)\text{-}\cat {PsHopf}(\tom (C, \Psi), H)$$
restricting to a quasi-isomorphism of chain algebras
 $$\xymatrix{\theta=\Theta (-\otimes z_0):\tom (C,\Psi)\ar [r]^(0.75)\simeq&H.}$$
\end{defn}

Unrolling the definition,  we see that the existence of $\Theta$ is equivalent to the existence of a family of $R$-linear maps
$$\mathfrak F(\Theta)=\{\theta_{n}=\Theta(-\otimes z_{n}): \tom(C,\Psi)\to H^{\otimes n+1}\}$$
satisfying certain conditions with respect to the differentials (cf., (\ref{eqn:family}) and Proposition \ref{prop:free-ext}).  Furthermore, it follows from the formula for the level comultiplication in $\op F$ (cf., \cite{HPST}, page 854)  that for all $a,b\in \tom(C,\Psi)$,
$$\theta _{n}(ab)= \sum_{\substack{ 1\leq k\leq n+1\\ \vec \imath\in I_{k,n+1}}}\pm \big((\Delta ^{(i_1)}\otimes \cdots \otimes \Delta ^{(i_k)})\theta _{k-1}(a)\big)\bullet \big ((\theta_{i_1 -1}\otimes \cdots \otimes \theta_{i_k-1})\Delta ^{(k)}(b)\big),$$
where $\bullet $ denotes the multiplication in $H^{\otimes n+1}$, $\Delta$ denotes the comultiplication in $\tom(C,\Psi)$ and in $H$ and  the signs follow from the Koszul rule.

Throughout this section we assume that any chain Hopf algebra mentioned is endowed with an Alexander-Whitney model, which we usually denote simply by $\theta:\tom (C,\Psi)\to H$.

Let $H$ be a chain Hopf algebra, and consider the acyclic $H$-comodule algebra $\Om H\otimes _{t_\Om} H$, as constructed in Proposition \ref{prop:cobar-alg}. Let 
$$\xymatrix{\mathfrak p:\Om H\otimes _{t_{\Om}}H\ar[r]&H}$$
denote the natural projection.

We explain first how to lift $\theta\circ \tom \pi$ naturally to a quasi-isomorphism 
$\ttheta:\lp (C,\Psi)\to \Om H\otimes _{t_{\Om}}H$ 
such that the following square commutes.
$$\xymatrix{\lp (C,\Psi )\ar [rr]^-{\ttheta}\ar [d] ^{{\tom} \pi}&&\Om H\otimes _{t_{\Om}}H\ar 
[d]^{\mathfrak p}\\
{\tom} (C,\Psi)\ar [rr]^{\theta }&&H}$$
We begin by defining and studying a certain section of $\mathfrak p$ and 
a derivation homotopy associated with it. Let $\xymatrix@1{\chi 
:H\ar [r]&H\otimes H}$ denote the comultiplication on $H$.

The proof of the following lemma  is an immediate consequence of the definitions.
\begin{lem} Define $\mathfrak s:H\to \Om H\otimes 
_{t_{\Om}}H$ to be the linear map of degree 0 given by $\mathfrak s (w)=1\otimes 
w$. Then 
\begin{enumerate}
\item $\mathfrak p\mathfrak s=1_{H}$;
\item $(D_{\Om}\mathfrak s-\mathfrak s d)(w)=-\si w_{i}\otimes w^i$ 
for all $w\in H$, where $\overline \chi (w)=w_{i}\otimes w^i$;
\item $\mathfrak s$ is a map of graded algebras; and 
\item $\mathfrak s$ is a map of right $H$-comodules.
\end{enumerate}
\end{lem}

Using this knowledge of $\mathfrak s$, we can build an important chain map from $H$ to $\Om H\otimes _{t_{\Om}}H$, as explained in the next lemma.

\begin{lem} Let $\mathfrak h=D_{\Om}\mathfrak s-\mathfrak s d:H\to 
\Om H\otimes _{t_{\Om}}H$. Then
\begin{enumerate}
\item $\mathfrak p\mathfrak h =0$;
\item $\mathfrak h$ is a chain map of degree $-1$, i.e., $D_{\Om}\mathfrak h =-\mathfrak h d$;
\item $\mathfrak h (a\cdot b)=\mathfrak h (a)\cdot \mathfrak s(b)+(-1)^a \mathfrak 
s(a)\cdot \mathfrak h (b)$ for all $a,b\in H$, i.e., $\mathfrak h$ is a 
derivation homotopy from $\mathfrak s$ to itself; and
\item $\mathfrak h$ is a map of right $H$-comodules.
\end{enumerate}
\end{lem}

\begin{proof}  We leave the trivial verifications of (1) and (2) to the reader.

\noindent (3) Observe that \begin{align*}
\mathfrak h (a\cdot b)&=D_{\Om}\mathfrak s(a\cdot b)-\mathfrak s d(a\cdot b)\\
&=D_{\Om}\big (\mathfrak s(a)\cdot 
\mathfrak s(b)\big)-\mathfrak s\big (da\cdot b+(-1)^a a\cdot db\big )\\
&=D_{\Om}\mathfrak s(a)\cdot \mathfrak s(b) +(-1)^a \mathfrak s(a)\cdot D_{\Om}\mathfrak s(b)-\mathfrak sd(a)\cdot \mathfrak s(b)-(-1)^a 
\mathfrak s(a)\cdot \mathfrak sd(b)\\
&= \mathfrak h  (a) \cdot \mathfrak s(b)+(-1)^a \mathfrak s(a)\cdot \mathfrak h  (b).
\end{align*}

\noindent (4) Observe that 
\begin{align*}
(1\otimes \chi )\mathfrak h&=(1\otimes \chi)(D_{\Om }\mathfrak s-\mathfrak sd)\\
&=(D_{\Om }\otimes 1 + 1\otimes d)(1\otimes \chi)\mathfrak s-(\mathfrak s\otimes 
1)\chi d\\
&=(D_{\Om }\otimes 1 + 1\otimes d)(\mathfrak s\otimes 1)\chi-(\mathfrak s\otimes 
1)(d\otimes 1 + 1\otimes d)\chi\\
&=(\mathfrak h \otimes 1)\chi. 
\end{align*}
\end{proof}

We now apply $\mathfrak s$ and $\mathfrak h$ to the construction of the lift 
of $\ttheta$.

\begin{thm} \label{thm:lift-props}Let $\ttheta :T \si (\cbc)_+\to  \Om H\otimes _{t_{\Om}}H$ be the 
graded algebra map specified by $\ttheta (\si c)= \mathfrak s\theta (\si c)$ and $\ttheta 
(\si \bc)= \mathfrak h \theta (\si c)$.
\begin{enumerate}
\item $\ttheta $ is a 
differential map, i.e., $\ttheta \td_{\Om}=D_{\Om}\ttheta $, and is therefore 
a quasi-isomorphism.
\item $\ttheta:\pl (C)\to  \Om H\otimes _{t_{\Om}}H$ is 
a CASH map.
\end{enumerate}
\end{thm}

\begin{proof} (1) Let $c\in C$ and write $\overline\Delta (c)=c_{i}\otimes c^i$.  Then
\begin{align*}
D_{\Om }\ttheta (\si c)&=D_{\Om }\mathfrak s\theta (\si c)\\
&=\mathfrak s d\theta (\si c)+\mathfrak h \theta (\si c)\\
&=\mathfrak s\theta d_{\Om}(\si c)+\ttheta (\si \bc)\\
&=\mathfrak s\theta \big(-\si (dc)+(-1)^{c_{i}}\si c_{i}\si c^i\big )+\ttheta (\si \bc)\\
&=-\ttheta \big(\si (dc)\big)+(-1)^{c_{i}}\mathfrak s\theta (\si c_{i})\mathfrak s\theta (\si c^i)+\ttheta (\si \bc)\\
&=-\ttheta \big(\si (dc)\big)+(-1)^{c_{i}}\ttheta (\si c_{i})\ttheta (\si c^i)+\ttheta (\si \bc)\\
&=\ttheta \big(-\si (dc)+\si \bc+(-1)^{c_{i}}\si c_{i}\si c^i\big)\\
&=\ttheta \big(-\si (\td c)+(-1)^{c_{i}}\si c_{i}\si c^i\big)\\
&=\ttheta \td_{\Om}(\si c).
\end{align*}
Furthermore,
\begin{align*}
D_{\Om }\ttheta (\si \bc)&=D_{\Om }\mathfrak h\theta (\si c)\\
&=-\mathfrak h d\theta (\si c)\\
&=-\mathfrak h\theta d_{\Om}(\si c)\\
&=-\mathfrak h\theta \big(-\si (dc)+(-1)^{c_{i}}\si c_{i}\si c^i\big )\\
&=\ttheta(\si \overline {dc})-(-1)^{c_{i}}\mathfrak h\big (\theta 
(\si c_{i})\theta (\si c^i)\big)\\
&=\ttheta(-\si \td\bc)\\
&\quad -(-1)^{c_{i}}\big (\mathfrak h\theta (\si 
c_{i})\cdot \mathfrak s\theta (\si c^i)+(-1)^{c_{i}+1}\mathfrak 
s\theta(\si c_{i})\cdot\mathfrak h\theta (\si c^i)\big )\\
&=\ttheta(-\si \td\bc)\\
&\quad -(-1)^{c_{i}}\big (\ttheta (\si 
\bc_{i})\cdot \ttheta (\si c^i)-(-1)^{c_{i}}\ttheta(\si 
c_{i})\cdot\ttheta (\si \bc^i)\big )\\
&=\ttheta\big(-\si \td\bc+(-1)^{c_{i}+1}\si \bc _{i}\si c^i 
+(-1)^c\si c_{i}\si \bc^i\big)\\
&=\ttheta \td _{\Om}(\si \bc).
\end{align*}

Observe that since $\ttheta$ is a differential map, it is 
necessarily a quasi-isomorphism, as both $\pl (C)$ and $\Om H\otimes _{t_{\Om}}H$ 
are acyclic.

\medskip
\noindent (2) Let $\Theta\in(\op A,\op F)\text{-}\cat {PsHopf}\big(\tom (C, \Psi), H\big)$ denote the pseudo $\op A$-Hopf algebra map that $\theta $ underlies. Let $\mathfrak F(\Theta )=\{\theta_k :\Om C\to H^{\otimes 
k}\mid k\geq 1\}$, where $\theta _k=\Theta (-\otimes z_{k-1}) $.  For $k\geq 0$, 
define 
$$\tTheta _k :\si (\cbc )\to (\Om H\otimes _{t_{\Om}}H)\otimes 
H^{\otimes k}$$  by
$$\tTheta _k(\si c)=( \mathfrak s\otimes 1^{\otimes k})\theta _{k+1}(\si c)$$
and
$$\tTheta_k (\si \bc)= ( \mathfrak h \otimes 1^{\otimes k})\theta 
_{k+1}(\si c),$$
where $1$ denotes the identity on $H$.

We claim that 
$$\{\tTheta_k :\si (\cbc)\to (\Om H\otimes _{t_{\Om}}H)\otimes 
H^{\otimes k}\mid k\geq 1\}$$ satisfies the hypotheses of Proposition \ref{prop:rel-free-ext} and therefore induces a morphism 
$$\widetilde\Theta\in(\op A,\op F)\text{-}\cat {PsHopf}\big (\pl (C, \Psi), H\otimes _{t_\Om}\Om H\big),$$ i.e., $\ttheta$ is a CASH map. We prove this claim by induction on $k$ and on degree in the Appendix.  

Note that $\widetilde\Theta$ lifts $\Theta$, in the sense that  $ \mathfrak p^{\otimes n}\widetilde \Theta (n)=\Theta (n)(\tom \pi \diamond 1_{\op F})$ for all $n$.   This is a necessary condition for $\ttheta$ to be a CASH map with respect to $\theta$.
 \end{proof}

\begin{summary} Given an object $H$ of $\mathbf H$ and an Alexander-Whitney model of $H$
$$\xymatrix{\theta:\tom (C,\Psi)\ar [r]^(0.6)\simeq & H,}$$
there exists 
a commutative diagram
\begin{equation}\label{diag:summary}
\xymatrix{ \lp (C,\Psi)\ar [d] ^{\tom \pi}\ar [rr]^-{\ttheta 
}&&\Om H\otimes _{t_{\Om}}H\ar [d]^{\mathfrak p}\\ 
{\tom }(C,\Psi)\ar [rr]^-{\theta }&&H} 
\end{equation}
such that
\begin{enumerate}
\item $\tom \pi$ is a strict algebra and 
coalgebra map;
\item $\mathfrak p$ is a strict algebra and right $H$-comodule map;
\item the natural right $\Om C$-comodule structure on $\lp(C,\Psi)$ is 
cofree on $\lp (C,\Psi)\square _{\Om C}R$;
and
\item $\ttheta $ is a quasi-isomorphism that 
is a strict algebra 
map and a CASH map.
\end{enumerate}
Furthermore, this construction is natural in $\theta$.
 \end{summary}

The next theorem, which is the heart of this article, describes how we can use the path-loop construction to compute the multiplicative structure of homotopy fiber homologies in $\mathbf H$.  Recall from Corollary  \ref{cor:cotor-alg} that if $H$ is a chain Hopf algebra and $M$ is an $H$-comodule algebra, then $\cotor ^H(M,R)=\H (M\otimes _{t_\Om} \Om H)$ has a natural graded algebra structure.  Recall furthermore from Corollary \ref{cor:pl-cofree} that $\tom\big((C',\Psi')\coprod\mathfrak P(C,\Psi)\big)\square _{\tom(C,\Psi)}R$ is a sub chain algebra of $\tom\big((C',\Psi')\coprod\mathfrak P(C,\Psi)\big)$ for all $\theta\in \cat F\big((C',\Psi'), (C,\Psi)\big)$.

\begin{thm}\label{thm:homol-loopfib} Let $\vp: H'\to H$ be a map of chain Hopf algebras.  Suppose that there is a a map $\omega: (C',\Psi')\to (C,\Psi)$ in $\shc$ and a commutative diagram
\begin{equation}\label{diag:AWmod}
\xymatrix{
{\tom (C',\Psi')}\ar[d]_{{\tom} \omega}\ar [r]^(0.6){\theta'}_(0.6)\simeq&H'\ar[d]_{\vp}\\
{\tom (C,\Psi)}\ar [r]^(0.6){\theta}_(0.6)\simeq &H}
\end{equation}
in which $\theta'$ and $\theta$ are Alexander-Whitney models.  Let 
$$\hf(\omega)=\tom \big(C'\coprod \mathfrak P(C)\big)\square_{\Om C} R.$$  
Then there is a  a zig-zag of quasi-isomorphisms of chain algebras
$$\hf(\omega)\xleftarrow\simeq \bullet\xrightarrow\simeq \cdots\xleftarrow\simeq\bullet\xrightarrow\simeq H'\otimes _{t_{\Om}}\Om H.$$ 
In particular, $\H \big (\hf(\omega)\big)$ is isomorphic to $\cotor ^H(H', R)$ as graded algebras.
\end{thm}

It is not surprising that there is at least a linear isomorphism between $\H \big (\hf(\omega)\big)$ and $\cotor ^H(H', R)$.  Since $\tom\big((C',\Psi')\coprod\mathfrak P(C,\Psi)\big)$ is a $\tom(C,\Psi)$-cofree resolution of $\tom (C',\Psi')$ and $\cotor$ is the derived functor of the cotensor product, $\cotor^{\tom(C,\Psi)}(\tom(C',\Psi'), R)$ should be the same as $\H \big(\hf(\omega)\big)$.  Naturality then gives rise to the linear isomorphism $\H \big (\hf(\omega)\big)$ and $\cotor ^H(H', R)$.
The challenge lies in showing that the isomorphism is multiplicative.

\begin{proof} The factorization of $\omega$ described in \ref{eqn:can-fact}
$$\xymatrix{C'\ar [rd]\ar [rr]^\omega&& C\\&C'\oplus (\cbc)\ar [ru]_{\omega+\pi} },$$
induces a factorization of $\tom \omega$
$$\xymatrix{\Om C'\ar [rd]_{}\ar [rr]^{\tom \omega }&& \Om C\\&\Om \big(C'\oplus (\cbc)\big)\ar [ru]_{\tom (\omega+\pi)} }.$$
Note that $\tom (\omega+\pi)$ admits a (nondifferential) Hopf algebra section $\Om \mathfrak j$, where $\mathfrak j$ is the natural section of $\omega +\pi$.  We can therefore apply Proposition \ref{prop:mm} to $\tom (\omega+\pi)$.

From Proposition \ref{prop:mm}, we know  that 
there is an injection of chain algebras 
$$i:\hf (\omega)\hookrightarrow \Om 
\big(C'\oplus (\cbc)\big)$$ 
and an isomorphism of left $\hf (\omega)$-modules and right $\Om 
C$-comodules
$$h:\hf (\omega)\otimes \Om C\to \Om\big(C'\oplus (\cbc)\big)$$
defined by $h=\mu (i\otimes\Om \mathfrak j)$.

Next, using the right $\Om C$-comodule structure of $ \Om 
\big(C'\oplus (\cbc)\big)$, 
define the twisted tensor product $ \Om 
\big(C'\oplus (\cbc)\big)\otimes_{t_{\Om}} \Om ^2 
C$ , and let $$\eta : \Om 
\big(C'\oplus (\cbc)\big)\hookrightarrow\Om\big(C'\oplus (\cbc)\big)\otimes_{t_{\Om}} \Om ^2 C$$ 
denote the natural (nondifferential) injection of algebras such that $\eta (w)=w\otimes 
1$.  We claim that the composition 
$$\eta i:\hf(\omega)\to \Om\big (C'\oplus (\cbc)\big)\otimes_{t_{\Om}} \Om ^2 C$$
is a quasi-isomorphism of chain algebras.  

Note that 
$\eta i$ factors as a composite of chain maps
$$\hf (\omega)\xrightarrow{\text{incl.}}(\hf (\omega)\otimes \Om C)\otimes _{t_{\Om}}\Om ^2 C\xrightarrow{h\otimes 1} \Om \big(C'\oplus (\cbc)\big)\otimes _{t_{\Om}}\Om ^2 C.$$
The linear map $h\otimes 1$ is a differential map because $h$ is a map 
of differential $\Om C$-comodules, while the inclusion $\hf(\omega)\hookrightarrow (\hf (\omega)\otimes \Om C)\otimes _{t_{\Om}}\Om ^2 C$ is a differential map since $\hf (\omega)\otimes \Om C$ is cofree. Since $\Om C\otimes _{t_{\Om}}\Om ^2 C$ is 
acyclic, 
the first, inclusion map is a 
quasi-isomorphism.  The second map is also a quasi-isomorphism, 
as $h$ is an isomorphism.  Thus, $\eta i$ is a quasi-isomorphism, as 
claimed.

Similarly, since there is a chain 
subalgebra inclusion 
$$H'\otimes _{t_\Om}\Om H\hookrightarrow H'\otimes _{t_\Om}(\Om H\otimes _{t_{\Om}} H)\otimes _{t_{\Om}} 
\Om H,$$
which is a quasi-isomorphism since $H\otimes _{t_{\Om}} \Om H$ is acyclic.

The conditions on diagrams (\ref{diag:summary}) and (\ref{diag:AWmod}) imply that we have morphisms in $\cat {CASH}$ (cf., Definition \ref{defn:cash})
$$(\theta; \ttheta): \big(\tom(C,\Psi); \pl (C,\Psi)\big)\to \big(H;\Om H\otimes_{t_{\Om}}H\big)$$
and
$$(\theta; \theta'):\big(\tom(C,\Psi); \tom (C',\Psi')\big)\to (H;H').$$
Applying Proposition \ref{prop:cash-coproduct},
we obtain a CASH map
$$\theta '\star \ttheta:\Om \big (C'\oplus (\cbc )\big)\to H'\otimes _{t_\Om}(\Om H\otimes _{t_{\Om}} H),$$
associated to a chain algebra map
$$\zeta:\Om \big (C'\oplus (\cbc )\big)\otimes _{t_{\Om}}\Om ^2 C\longrightarrow H'\otimes _{t_\Om}(\Om H\otimes _{t_{\Om}} H)\otimes _{t_{\Om}} 
\Om H.$$
Since $\Om (\cbc)$ and $\Om H\otimes _{t_\Om} H$ are acyclic, and $\Om \big (C'\oplus (\cbc )\big)\cong \Om C'\coprod \Om (\cbc)$, the vertical arrows in the commuting diagram 
$$\xymatrix{\Om \big (C'\oplus (\cbc )\big)\ar [rr]^{\theta '\star \ttheta}&&H'\otimes _{t_\Om}(\Om H\otimes _{t_{\Om}} H)\ar[d] ^{1\otimes \epsilon}_\simeq\\
\Om C'\ar [u]^{\text{incl.}}_\simeq \ar [rr]^{\theta '}_\simeq &&H',}$$
where $\epsilon$ denotes the augmentation map, are both quasi-isomorphisms. Consequently, $\theta'\star \ttheta$ is a quasi-isomorphism, which implies that $\zeta$ is also a quasi-isomorphism, since the underlying graded modules of all objects involved are assumed to be free over $R$.

We have therefore a zig-zag of quasi-isomorphisms of chain algebras
$$\xymatrix {\hf (\omega)\ar [d]^{\simeq}_{\eta i}&&H'\otimes _{t_\Om}\Om H\ar [d] 
^{\simeq}_{\text{incl.}}\\
\Om\big (C'\oplus (\cbc)\big)\otimes_{t_{\Om}} \Om ^2 C\ar [rr]^-{\zeta}_-{\simeq}&&H'\otimes _{t_\Om}(\Om H\otimes _{t_{\Om}} H)\otimes _{t_{\Om}} 
\Om H.}$$
Consequently, 
$$\H \hf (\omega)\cong \H (H'\otimes _{t_\Om}\Om H)= \cotor ^H(H', R)$$ 
as graded algebras.
 \end{proof}

\begin{cor} Let $\theta :\tom (C,\Psi)\to H$ be an Alexander-Whitney model of a chain Hopf algebra $H$.
If 
$$\dl (C,\Psi)=\lp(C,\Psi)\square _{\tom (C,\Psi)}R,$$
with its natural chain algebra structure inherited from $\tom (C,\Psi)$, then $\H \big(\dl (C,\Psi)\big)$ is isomorphic 
to $\cotor ^H(R, R)$ as graded algebras.\end{cor}

\begin{proof} Apply Theorem \ref{thm:homol-loopfib} to the unit map $R\to H$. \end{proof} 

\begin{defn} Given an object $(C,\Psi)$ in $\shc$, the chain algebra $\dl (C,\Psi )$ is the \emph {double-loop construction} on $C$.  Given a morphism $\omega :(C',\Psi')\to (C,\Psi)$ of $\shc$, the chain  algebra $\hf(\omega)$ is the \emph {loop-homotopy fiber} construction on $\omega$. 
\end{defn}


\begin{rmk} Note that $\dl (C,\Psi)$ is a subalgebra of $\pl (C,\Psi)=\tom\circ \widetilde {\mathfrak P}(C,\Psi)$, which is free as a graded algebra on $C\oplus \overline C$.  In particular, if $C$ admits an $R$-basis of $n$ elements, then $\dl (C,\Psi)$ is a subalgebra of a free algebra on $2n$ generators.  On the other hand, $\Om\big(\tom (C,\Psi)\big)$, which is connected by a zig-zag of  quasi-isomorphisms of chain algebras  to $\dl (C,\Psi)$, is free on an infinite number of generators.

A similar comparison can be made between $\hf (\omega)$ and $\Om\big (C'\oplus (\cbc)\big)\otimes_{t_{\Om}} \Om ^2 C$.
\end{rmk}

\begin{rmk}   Since the path-loop construction is functorial, the double-loop construction is as well, i.e., a morphism $h:(C,\Psi)\to (C',\Psi')$ in $\shc$ induces a chain algebra map $\dl (h):\dl (C,\Psi)\to \dl (C',\Psi')$.  It is evident that if $h$ is a quasi-isomorphism in $\shc$, then $\dl (h)$ is a quasi-isomorphism of chain algebras.

Similarly, the loop-homotopy fiber construction is clearly defines a functor from the category of morphisms in $\shc$ to the category of chain algebras. Furthermore if 
$$\xymatrix{(C',\Psi')\ar [r]^\omega\ar [d]_\simeq ^{h '}&(C,\Psi)\ar [d]_\simeq ^h \\
(B',\Phi')\ar [r]^\zeta&(B,\Phi)}$$
is a commuting diagram in $\shc$, where $h$ and $h '$ are quasi-isomorphisms, then the induced map $\hf (\omega)\to \hf (\zeta)$ is a quasi-isomorphism of chain algebras.
\end{rmk}



\section {The loops on a homotopy fiber}\label{sec:topology}

We are now ready to apply the purely algebraic results above to topology.

\subsection{The loop-homotopy fiber model}

Here we apply the constructions and theorems of the previous two sections to constructing a chain algebra, the homology of which is isomorphic as a graded algebra to $\H ( GF)$, where $F$ is the homotopy fiber of a morphism $g:K\to L$ of $2$-reduced simplicial sets.  

We begin by specifying our input data for the constructions of sections \ref{sec:path-loop} and \ref{sec:homol-htpyfib}: the \emph{canonical enriched Adams-Hilton model} of \cite {HPST}.

Recall from Theorem \ref{thm:loopmodel} that there is a functor $\tC:\mathbf{sSet}_1\to \mathbf F$. In \cite {S} Szczarba gave an explicit formula for a 
natural transformation between functors from $\mathbf{sSet}_1$ to the category of associative chain algebras
$$\theta: \Om C(-)\to C(G(-))$$
such that $\theta_{K}: \Om C(K)\to C(GK)$ is a quasi-isomorphism of 
chain algebras for every $1$-reduced simplicial set $K$.  

Since $\Om$ extends to a functor $\tom:\shc \to  \mathbf H$ (see (1.5)), there is a natural transformation $\psi : \Om C(-)\to \Om C(-)\otimes \Om C(-)$ given for each $1$-reduced simplicial set $K$ by the composition
$$\Om C(K)\xrightarrow{\ind (\Psi _K)}\Om \big (C(K)\otimes C(K)\big )\xrightarrow{q}\Om C(K)\otimes \Om 
C(K),$$
where $\tC(K)=\big(C(K), \Psi _K)$. The comultiplication $\psi _{K}:\Om C(K)\to \Om C(K)\otimes \Om C(K)$ is 
called the \emph { Alexander-Whitney cobar diagonal}. 

 In \cite {HPST} Hess, Parent, Scott and Tonks  established 
that Szczarba's equivalence $\theta_{K}$ underlies a pseudo $\op A$-Hopf map with respect 
to $\psi _{K}$ and the usual comultiplication on $C(GK)$. In other words, for all $1$-reduced $K$, 
 $$\theta_{K}:\tom \tC(K)\xrightarrow{\simeq}C(GK)$$
is an Alexander-Whitney model, called the \emph{canonical enriched Adams-Hilton model} of $K$.
Finally, they showed that $\psi _K$ agrees with the comultiplication on $\Om C(K)$ defined in a purely combinatorial manner by Baues in \cite {Ba}.

We now apply the canonical enriched Adams-Hilton model to modelling the loop homology of homotopy fibers.

\begin{thm}\label{thm:loopfib-model} Let $f:K\to L$ be a morphism of $2$-reduced simplicial sets, and let $F$ be the homotopy fiber of $f$. Then there is a zig-zag of quasi-isomorphisms of chain algebras
$$\hf(f)\xleftarrow\simeq \bullet\xrightarrow\simeq \cdots\xleftarrow\simeq\bullet\xrightarrow\simeq C(GF).$$  
Thus, $$\H \big(\hf (\widetilde C(f))\big)\cong \H (GF)$$
as graded algebras.  In particular, 
$$\H \big(\dl (\widetilde C(L))\big)\cong \H (G^2L)$$
as graded algebras.\end{thm}

\begin{proof} Applying Theorem \ref{thm:homol-loopfib} to the commuting diagram
$$\xymatrix{
\tom \widetilde C(K)\ar [r] ^{\theta _K}_\simeq \ar [d] _{\tom \widetilde C(f)}&C(GK)\ar [d]_{C(Gf)} \\
\tom \widetilde  C(L)\ar [r]^{\theta_L} _\simeq&C(GL)}$$
we obtain a zig-zag of quasi-isomorphisms of chain algebras
$$\hf(f)\xleftarrow\simeq \bullet\xrightarrow\simeq \cdots\xleftarrow\simeq\bullet\xrightarrow\simeq C(GK)\otimes _{t_{\Om}}\Om C(GL).$$ 
By the dual of Theorem 5.1 in \cite{fht}, there is a quasi-isomorphism of chain algebras 
$$C(GF)\xrightarrow\simeq C(GK)\otimes _{t_{\Om}}\Om C(GL),$$
and so we can conclude.
 \end{proof}

\subsection {Double suspensions and formal spaces}\label{subsec:formal}

In this section we provide a more explicit description of $\dl (C,\Psi)$ for Alexander-Whitney coalgebras $(C,\Psi)$ such that $\tom(C,\Psi)$ is primitively generated.  Since $\tom\widetilde C(K)$ is primitively generated for all simplicial double suspensions $K$ \cite {HPS2}, we have good control of the model $\dl (K)$ for a large class of spaces $K$.

 Our description of $\dl (C,\Psi)$ also applies to Alexander-Whitney coalgebras that are ``formal" in some appropriate sense.  The chain coalgebras of numerous interesting spaces satisfy our formality criteria, enabling us to  give a more explicit and computationally amenable formula for  the model of their double loop spaces as well.

Our notion of formality is certainly closely related to that of Ndombol and Thomas \cite {ndombol-thomas}, but we do not know whether the two notions are actually equivalent.

\begin{defn} A morphism $\theta\in \mathbf F\big( (C, \Psi), (C', \Psi')\big)$ is an \emph { $\mathbf F$-quasi-isomorph\-ism} if $\theta (-\otimes z_0):C\to C'$ is a quasi-isomorphism.  Two Alexander-Whitney coalgebras $(C,\Psi)$ and $(C',\Psi')$ are \emph { weakly equivalent} if there is a zig-zag of $\mathbf F$-quasi-isomorphisms in $\shc$
$$\xymatrix{(C,\Psi )&\bullet\ar [l]_(0.25)\simeq\ar [r]^(0.45)\simeq&\cdots&\bullet\ar [r]^(0.25)\simeq\ar[l]_(0.4)\simeq&(C',\Psi')}.$$\end{defn}

Let $(C,\Psi)$ and $(C',\Psi')$  be weakly equivalent Alexander-Whitney coalgebras, with underlying coalgebras $(C,\Delta)$ and $(C',\Delta ')$.  It follows immediately from the definition that $H_*(C,\Delta)\cong H_*(C',\Delta')$ as cocommutative coalgebras and that $H_*\big(\tom (C,\Psi)\big)\cong H_*\big(\tom (C',\Psi')\big)$ as Hopf algebras.

Recall the ``inclusion" functor $\mathfrak I_{\op F}: \acoalg\to \Coalg AF$ from (\ref{eqn:emb-fat}).

\begin{defn} An Alexander-Whitney coalgebra $(C,\Psi)$ is \emph {formal} if it is weakly equivalent to $\mathfrak I _{\op F}\big(\H (C)\big)$.  A simplicial set $K$ is \emph {Alexander-Whitney formal} if $\tC(K)$ formal in $\cat F$.\end{defn}

Spheres are obviously Alexander-Whitney formal spaces.  More generally, if $K$ is an $r$-reduced simplicial set  such that all simplices of dimension greater than $2r$ are degenerate, then $K$ is clearly Alexander-Whitney formal for degree reasons. 

We now recall algebraic notions and results from \cite {CMN} and \cite {T} that enable us to simplify $\dl (C,\Psi )$ in the formal case.   We work henceforth over any integral domain $R$ in which $2$ is a unit or over a field $R$ of characteristic $2$.

\begin{defn}\cite[Section 3] {CMN}\label{defn:homol-hopf-alg}  A \emph{homology Hopf algebra} is a connected, cocommutative, graded Hopf algebra that is free as a graded $R$-module.   Given homology Hopf algebras $H$, $H'$ and $H''$, a sequence $H'\xrightarrow i H\xrightarrow p H''$ of morphisms of Hopf algebras is a \emph{short exact sequence} if
\begin{enumerate}
\item the composite $pi$ is equal to the composite $H' \xrightarrow \varepsilon R \xrightarrow \eta H''$;
\item $i$ is injective, while $p$ is surjective; and
\item the canonical map $\bar\imath : H' \to H\square_{H''} R$ is an isomorphism.
\end{enumerate}
\end{defn}

\begin{prop}\cite[Proposition 3.7]{CMN}\label{prop:cmn} If $L'\to L \to L''$ is a short exact sequence of connected, graded Lie algebras over $R$, then $UL'\to UL\to UL''$ is a short exact sequence of homology Hopf algebras, where $U$ denotes the universal enveloping algebra functor, from graded Lie algebras to graded Hopf algebras.  .
\end{prop}

\begin{prop}\cite [Proposition VI.2 (7)]{T}\label{prop:tanrŽ} Let $p:\mathbb L (V\oplus W) \to \mathbb L(V)$ be the projection map of graded Lie algebras, determined by $p(v)=v$ and $p(w)=0$ for all $v\in V$ and $w\in W$.  There is then a short exact sequence of Lie algebras
$$\mathbb L\big(\mathfrak {A}^{\mathbb N}(V)(W)\big)\xrightarrow i L(V\oplus W) \xrightarrow p \mathbb L(V),$$
where
$$\mathfrak{A}^{\mathbb N}(V)(W):=\big\{[v_1,[v_2,[...[v_k,w]...]]]\mid w\in W, v_i\in V\;\forall i, k\in \mathbb N\big \}$$
and  $[-,-]$ denotes a commutator.  
\end{prop}

We can now apply the results recalled above to determining the underlying graded algebra of $\dl (C,\Psi)$ when $\tom (C,\Psi)$ is primitively generated.

\begin{thm}\label{thm:formal-model} Let $R$ be either an integral domain in which $2$ is a unit or a field of characteristic 2. If $(C,\Psi )$ is an Alexander-Whitney coalgebra over $R$ such that $\tom (C,\Psi)$ is primitively generated with respect to the naturally induced comultiplication $\psi$, then the underlying graded algebra of $\dl (C,\Psi)$ is $T\big(\mathfrak {A}^{\mathbb N}(\si C)(\si \bC)\big)$.\end{thm}

\begin{proof} Since the underlying graded Hopf algebra of $\tom (C, \Psi)$ is primitively generated and therefore cocommutative, so is the underlying graded Hopf algebra of $\mathfrak{PL} (C, \Psi)$.  It follows that there is a short exact sequence of homology Hopf algebras
\begin{equation}\label{eqn:sequence}
T\si (C_{+}\oplus \overline C_{+})\square _{T\si C_{+}}R \to  T\si (C_{+}\oplus \overline C_{+}) \xrightarrow {\Om \pi} T\si C_{+},
\end{equation}
where $T\si (C_{+}\oplus \overline C_{+})\square _{T\si C_{+}}R$ is the graded Hopf algebra underlying $\dl (C, \Psi)$. 

On the other hand, since each of the Hopf algebras above is primitively generated, $\Om \pi=Up$, where $p:\mathbb L\si (C_{+}\oplus \overline C_{+}) \xrightarrow {p} \mathbb L\si C_{+}$ is the usual projection map of graded Lie algebras.  Combining Propositions \ref{prop:cmn} and \ref{prop:tanrŽ}, we obtain another short exact sequence of homology Hopf algebras
$$T\big(\mathfrak {A}^{\mathbb N}(\si C)(\si \bC)\big)\to T\si (C_{+}\oplus \overline C_{+}) \xrightarrow {\Om \pi} T\si C_{+}.$$
Comparing with sequence (\ref{eqn:sequence}), we conclude that 
$$T\si (C_{+}\oplus \overline C_{+})\square _{T\si C_{+}}R\cong T\big(\mathfrak {A}^{\mathbb N}(\si C)(\si \bC)\big).$$
\end{proof}

\begin{cor}  Let $R$ be either an integral domain in which $2$ is a unit or a field of characteristic 2. If $K=E^2L$ is a simplicial double suspension, then the graded algebra underlying $\dl(K)$ is $T\big(\mathfrak {A}^{\mathbb N}(C_{+}EL)(C_{+}L)\big)$.
\end{cor}

\begin{proof}  In \cite{HPS2} the authors proved that $\tom\widetilde C(E^2 L)$ was primitively generated for all simplicial sets $L$.  To conclude, observe that $\si C_{+}EX\cong C_{+}X$ for all simplicial sets $X$.
\end{proof}

\begin{cor} Let $R$ be either an integral domain in which $2$ is a unit or a field of characteristic 2. If $(C,\Psi)$ is  formal, then $\dl (C,\Psi )$ is weakly equivalent to a chain algebra with underlying graded algebra $T\big(\mathfrak {A}^{\mathbb N}(\si H)(\si \overline H)\big)$, where $H=\H (C)$.\end{cor}

\begin{proof} Writing $\mathfrak I(H)=(H,\Psi_H)$, we have by definition that $\Psi _{H}= \op T(\Delta _{H})\acirc \varepsilon$.  The induced comultiplication on $\tom (H,\Psi _H)$ is such that $\tom (H,\Psi _H)$ is primitively generated.  Now apply Theorem \ref{thm:formal-model}.
 \end{proof}

If $(C,\Psi)$ is formal and  all elements of $\H (C)$ are primitive, then the differential $\td _\Om$ on $\dl (H,\Psi _H)$ is given explicitly by
\begin{multline*}
\td_\Om \big([\si x_1,[\si x_2,[...[\si x_m, \si \bar y]...]]]\big)\\=\sum _{1\leq i\leq m} \pm [\si x_1,[\si x_2,[...[\si \bar x_i,[ ...[\si x_m, \si \bar y]...]]],\end{multline*}
 for all $x_1,..., x_m, y\in H$, where the sign is determined by the Koszul convention.  It is not too difficult to see in this case that  if $B$ is a basis of $\H (C,\mathbb F_2)$, then
 $$\H \big(\dl (C,\Psi );\mathbb F_2\big)\cong  \mathbb F_2[\operatorname {ad}^{2^k-1}(x)(\bar x)\mid x\in B, k\geq 1].$$ 
A similar result holds mod $p$.


\section{Appendix: Technical proofs}

\subsection{Proof of Proposition \ref{prop:free-ext}}
Since $\op Q$ is free as a $\op P$-bimodule, there is a symmetric sequence $\op X$ of graded $R$-modules such that $\op Q\cong \op P\diamond \op X\diamond \op P$ in ${}_{\op P}\cat {Mod}_{\op P}$.

The morphism $\theta$ consists of a morphism of right $\op P$-modules
$$\theta:\op T(H)\pcirc \op Q\to \op T(H')$$ 
such that 
$$\theta\mathfrak I_{\op Q}(\mu)=\mathfrak I_{\op Q}(\mu ')(\theta \curlywedge \theta):\op T(H\wedge H)\pcirc \op Q\to  \op T(H'),$$
where the $\op P$-coalgebra maps $\mu $ and $\mu'$ are the multiplication maps on $H$ and $H'$, respectively.

Let $\rho: \op T(H\coprod Tv)\diamond \op P\to \op T(H\coprod Tv)$ be the right $\op P$-module structure corresponding to the $\op P$-coalgebra structure of $H\coprod Tv$.    Since $H\coprod Tv$ is an extension of $H$, the action $\rho$ restricts to $\rho: \op L(R\cdot v)\diamond \op P\to  \op T(H\oplus R\cdot v)$.  In particular, the direct sum of graded modules $H\oplus R\cdot v$ underlies a sub $\op P$-coalgebra of $H\coprod Tv$.

Taken together, $\lambda: \op L(R\cdot v)\diamond \op X\to  \op T(H')$ and the restriction $\theta:\op L(H)\diamond \op X\to \op T(H')$ give rise in the obvious way to a morphism of symmetric sequences  from $\op L(H\oplus R\cdot v)\diamond \op X$ to $\op T(H')$, which we call $\theta\boxplus \lambda$, and then, by Lemma \ref{lem:inducedT}, to 
$$(\theta\boxplus \lambda)\breve{}:\op T(H\oplus R\cdot v)\diamond \op X\to \op T(H').$$
Furthermore, from $(\theta\boxplus \lambda)\breve{}$ and the restriction of $\rho$, we obtain a morphism $\lambda ':\op L(R\cdot v)\diamond \op Q\to \op T(H')$ of right $\op P$-modules as the composite
$$\op L(R\cdot v)\diamond \op P\diamond \op X\diamond \op P\xrightarrow{\rho\diamond 1\diamond 1}\op T(H\oplus R\cdot v)\diamond\op X\diamond \op P\xrightarrow{(\theta\boxplus \lambda)\breve{}\diamond 1}\op T(H')\diamond \op P\xrightarrow{\rho '}\op T(H'),$$
where $\rho'$ is the right $\op P$-module structure map corresponding to the $\op P$-coalgebra structure of $H'$.

We next recursively define a filtration of $H\coprod Tv$ by sub $\op P$-coalgebras, then construct $\widehat\theta$ by induction on filtration degree.  Set $F^0=H$ and, for $m>0$, $F^m$ is the image of the iterated multiplication map restricted to  $F^{m-1}\otimes (R\oplus R\cdot v)\otimes H$.  By definition of a free algebraic $\op P$-Hopf extension, each $F^m$ is a sub $\op P$-coalgebra of $H\coprod Tv$.   Furthermore the multiplication $\mu$ on $H\coprod Tv$ restricts to a morphism of $\op P$-coalgebras 
$$\mu _{k,m}:F^k\otimes F^{m-k}\to F^m,$$ 
for all $0\leq k\leq m$ and for all $m\geq 0$.

Let $\widehat\theta^{(0)}=\theta$.  Suppose that for some $m\geq 0$ and for all $k\leq m$, there exists 
$$\widehat \theta^{(k)}:\op T(F^k)\pcirc \op Q\to \op T(H')$$
such that for all $j\leq k$, $\widehat \theta^{(k)}$ agrees with $\widehat \theta^{(j)}$ on $\op T(F^j)\pcirc \op Q$ and 
\begin{equation}\label{eqn:skeleton}
\widehat \theta^{(k)} \mathfrak I(\mu _{j,k})=\mathfrak I(\mu ')\big(\widehat \theta^{(j)}\curlywedge\widehat \theta^{(k-j)}\big):\op T(F^j\wedge F^{k-j})\pcirc \op Q\to  \op T(H'). 
\end{equation}

To construct $\widehat \theta^{(m+1)}$, we begin by defining a morphism $\theta ^{(m+1)}$ in $\cat{Mod}_R^\Sigma$ from $\op L(F^{m+1})\diamond \op X$ to $\op T(H')$, to which we then apply Lemma \ref{lem:inducedT}.  The morphism $\theta ^{(m+1)}$ is defined to be the following composite.
$$\xymatrix{ {\op L}(F^{m+1})\diamond \op X\ar [ddd]_{\theta ^{(m+1)}}\ar [r]&\big(\op L(F^m)\wedge \op L(R\cdot v)\wedge \op L(H)\big)\diamond \op X\ar [d] ^{1\diamond \Delta _{\op Q}^{(2)}}\\&\big(\op L(F^m)\wedge \op L(R\cdot v)\wedge \op L(H)\big)\diamond \op Q^{\wedge 3}\ar[d]^{\mathfrak i}\\
&\big(\op L(F^m)\diamond \op Q\big)\wedge\big(\op L(R\cdot v)\diamond \op Q\big)\wedge \big(\op L(H)\diamond \op Q\big) \ar [d]^{\widehat \theta ^{(m)}\wedge \lambda '\wedge \widehat\theta}\\
\op T(H')&\op T\big((H')^{\otimes 3}\big)\cong \op T(H')^{\wedge 3}\ar [l]_{\op T\big((\mu ')^{(2)}\big)}}$$
Here $\widehat\theta ^{(m)}$ is slightly abusive shorthand for the composite
$$\op L(F^m)\diamond \op Q\hookrightarrow \op T(F^m)\diamond \op Q\to \op T(F^m)\pcirc \op Q\xrightarrow{\widehat \theta^{(m)}} \op T(H').$$

Applying Lemma \ref{lem:inducedT}, we obtain $\breve{\theta} ^{(m+1)}:\op T(F^{m+1})\diamond \op X\to  \op T(H')$. The composite
$$\op T(F^{m+1})\diamond \op Q\cong \op T(F^{m+1})\diamond \op P\diamond\op X\diamond\op P\xrightarrow{\rho\diamond 1\diamond 1}\op T(F^{m+1})\diamond\op X\diamond\op P\xrightarrow{\breve{\theta} ^{(m+1)}\diamond \op P}\op T(H')\diamond \op P\xrightarrow{\rho'}\op T(H')$$
then induces the desired map $\widehat\theta ^{(m+1)}:\op T(F^{m+1})\pcirc \op Q\to \op T(H')$, since $\op Q$ is a free bimodule.  By construction, equality (\ref{eqn:skeleton}) now holds for all $k\leq m+1$ and  $j\leq k$.  

To complete the proof, set $\widehat \theta=\operatorname{colim}_m \widehat\theta ^{(m)}$.

\subsection{Complete proof of Theorem \ref{thm:lift-props}(2)}
Recall that we have fixed $\mathfrak F(\Theta )=\{\theta _{k}:\Om C\to H^{\otimes 
k}\mid k\geq 1\}.$  
 Furthermore, for $k\geq 0$, we have
defined 
$$\tTheta _{k}:\si (\cbc )\to (\Om H\otimes _{t_{\Om}}H)\otimes 
H^{\otimes k}$$  by
$$\tTheta _{k}(\si c)=( \mathfrak s\otimes 1^{\otimes k})\theta _{k+1}(\si c)$$
and
$$\tTheta _{k}(\si \bc)= ( \mathfrak h \otimes 1^{\otimes k})\theta 
_{k+1}(\si c),$$
where $1$ denotes the identity on $H$.

We claim that 
$$\{\tTheta_k :\si (\cbc)\to (\Om H\otimes _{t_{\Om}}H)\otimes 
H^{\otimes k}\mid k\geq 1\}$$ satisfies the hypotheses of Proposition 
\ref{prop:rel-free-ext} and therefore induces a morphism $$\widetilde\Theta\in(\op A,\op F)\text{-}\cat {PsHopf}(\pl (C, \Psi), H\otimes _{t_\Om}\Om H),$$
i.e., $\ttheta$ is a CASH map.
We have already dealt with 
the case $k=0$, since $\tTheta _{0}=\ttheta $. Suppose that the 
claim is true for all $k<m$ and for $\tTheta _{m}$ restricted to 
$\si (\cbc)_{<n}$.

Before proving the claim for $\tTheta _{m}$ applied to $\si 
(\cbc)_{n}$, we establish some useful notation.  For all $\mathbf j\in 
J_{k,m}$, let 
$$\hth _{\mathbf j}=\tTheta _{j_{0}}\otimes \theta _{j_{1}}\otimes 
\cdots\otimes \theta _{j_{k}}\text{ and }\theta _{\mathbf j}=\theta _{j_{0}}\otimes \theta _{j_{1}}\otimes 
\cdots\otimes \theta _{j_{k}}$$
as maps from $\Om (\cbc)\otimes \Om C^{\otimes 
k}$ to $\big(\Om H\otimes_{t_{\Om}} H\big)\otimes H 
^{\otimes m}$ and from $\Om C ^{\otimes k+1}$ to $H^{\otimes m}$.
Furthermore, let
$$\hn ^{\mathbf j}=\nu ^{j_{0}}\otimes \chi ^{(j_{1}-1)}\otimes 
\cdots\otimes\chi ^{(j_{k}-1)} \text{ and } \chi ^{\mathbf j}=\chi
^{(j_{0}-1)}\otimes 
\cdots\otimes\chi ^{(j_{k}-1)}$$
as maps from $\big(\Om H\otimes_{t_{\Om}} H\big)\otimes H 
^{\otimes k}$ to $\big(\Om H\otimes_{t_{\Om}} H\big)\otimes H 
^{\otimes m}$ and from $H^{\otimes k+1}$ to $H^{\otimes m}$.  Recall that $\chi $ is the comultiplication on $H$ and that $\nu$ is the right $\Om C$-comodule action on $\pl (C)$.

In terms of this notation $\tTheta _{m}$ 
can be defined recursively by
$$\tTheta _{m}\mu =\mu \left (\sum _{\substack{ 0\leq k\leq m\\\mathbf j\in 
J_{k,m}}} \hn ^{\mathbf j}\tTheta _{k}\otimes \hth _{\mathbf j}\nu 
^{(k)}\right ),$$
where all multiplication maps are denoted $\mu$ and all coactions are 
denoted $\nu$. 

Observe that since $\mathfrak s$ and $\mathfrak h$ are  maps of comodules, $\nu ^{(j)}\mathfrak s=(\mathfrak 
s\otimes 1^{\otimes j})\chi ^{(j-1)}$ and $\nu ^{(j)}\mathfrak h=(\mathfrak 
h\otimes 1^{\otimes j})\chi ^{(j-1)}$ for all $j$.

Let $\si c\in \si (\cbc)_{n}$. Let $c_{i}\otimes c^i$ denote the image of $c$ under the 
reduced comultiplication in $C$.  Then
\begin{align*}
\tTheta _{m}d_{\Om}(\si c)=&\tTheta _{m}\big(-\si 
(dc)+\si \bc+(-1)^{c_{i}}\si c_{i}\si c^i\big)\\
=&-(\mathfrak s\otimes 1^{\otimes m}) \theta _{m+1}(\si (dc))+(\mathfrak h\otimes 1^{\otimes m}) \theta _{m+1}(\si c)\\
&+(-1)^{c_{i}}\mu \biggl (\sum _{\substack{ 0\leq k\leq m\\ \mathbf j\in 
J_{k,m}}} \hn ^{\mathbf j}\tTheta _{k}(\si c_{i})\otimes \hth _{\mathbf j}\nu 
^{(k)}(\si c^i)\biggr )\\
=&-(\mathfrak s\otimes 1^{\otimes m}) \theta _{m+1}(\si (dc))+(\mathfrak h\otimes 1^{\otimes m}) \theta _{m+1}(\si c)\\
&+(-1)^{c_{i}}\mu \biggl (\sum _{\substack{ 0\leq k\leq m\\\mathbf j\in 
J_{k,m}}} \hn ^{\mathbf j}(\mathfrak s\otimes 1^{\otimes k})\theta _{k}(\si 
c_{i})\otimes \hth _{\mathbf j}\psi 
^{(k)}(\si c^i)\biggr )\\
=&-(\mathfrak s\otimes 1^{\otimes m}) \theta _{m+1}(\si (dc))+(\mathfrak h\otimes 1^{\otimes m}) \theta _{m+1}(\si c)\\
&+(-1)^{c_{i}}\mu \biggl (\sum _{\substack{ 0\leq k\leq m\\\mathbf j\in 
J_{k,m}}} (\mathfrak s\otimes 1^{\otimes m})\chi ^{\mathbf j}\theta _{k}(\si 
c_{i})\otimes (\mathfrak s\otimes 1^{\otimes m})\theta _{\mathbf j}\psi 
^{(k)}(\si c^i)\biggr )\\
=&-(\mathfrak s\otimes 1^{\otimes m}) \theta _{m+1}(\si (dc))+(\mathfrak h\otimes 1^{\otimes m}) \theta _{m+1}(\si c)\\
&+(-1)^{c_{i}}(\mathfrak s\otimes 1^{\otimes m})\mu \biggl (\sum _{\substack{ 0\leq k\leq m\\\mathbf j\in 
J_{k,m}}} \chi ^{\mathbf j}\theta _{k}(\si 
c_{i})\otimes \theta _{\mathbf j}\psi 
^{(k)}(\si c^i)\biggr )\\
=&-(\mathfrak s\otimes 1^{\otimes m}) \theta _{m+1}(\si (dc))+(\mathfrak h\otimes 1^{\otimes m}) \theta _{m+1}(\si c)\\
&+(-1)^{c_{i}}(\mathfrak s\otimes 1^{\otimes m})\theta _{m+1}(\si c_{i}\cdot 
\si s^i)\\
=&(\mathfrak s\otimes 1^{\otimes m}) \theta _{m+1}d_{\Om}(\si c)+(\mathfrak h\otimes 1^{\otimes m}) \theta _{m+1}(\si c).
\end{align*}
Next,
\begin{align*}
(D_{\Om}\otimes 1^{\otimes m})\tTheta _{m}(\si c)=&(D_{\Om}\mathfrak s\otimes 1^{\otimes 
m})\theta _{m+1}(\si c)\\
=&(\mathfrak s d\otimes 1^{\otimes m})\theta _{m+1}(\si c)+(\mathfrak h\otimes 
1^{\otimes m})\theta _{m+1}(\si c)\\
=&(\mathfrak s \otimes 1^{\otimes m})(d\otimes 1^{\otimes m})\theta 
_{m+1}(\si c)\\
&+(\mathfrak h\otimes 
1^{\otimes m})\theta _{m+1}(\si c)
\end{align*}
and
\begin{align*}
\big (1\otimes \sum _{r+s=m-1} 1^{\otimes r}\otimes d\otimes 
&1^{\otimes s}\big)\tTheta _{m}(\si c)\\
&=\big (\mathfrak s\otimes \sum _{r+s=m-1} 1^{\otimes r}\otimes d\otimes 
1^{\otimes s}\big)\theta _{m+1}(\si c)\\
&=(\mathfrak s\otimes 1^{\otimes m})\left(\sum _{r+s=m-1} (1^{\otimes r+1}\otimes d\otimes 
1^{\otimes s}\big)\theta _{m+1}(\si c)\right).
\end{align*}
On the other hand
\begin{align*}
(\nu \otimes 1^{\otimes m-1})\tTheta _{m-1}(\si c)&=(\nu\mathfrak s \otimes 
1^{\otimes m-1})\theta _{m}(\si c)\\
&=(\mathfrak s\otimes 1^{\otimes m})(\chi\otimes 1^{\otimes m-1})\theta 
_{m}(\si c)
\end{align*}
and 
\begin{align*}
\big (1\otimes &\sum _{r+s=m-2} 1^{\otimes r}\otimes \chi\otimes 
1^{\otimes s}\big)\tTheta _{m-1}(\si c)\\
=&\big (\mathfrak s\otimes \sum _{r+s=m-2} 1^{\otimes r}\otimes \chi\otimes 
1^{\otimes s}\big)\theta _{m}(\si c)\\
=&(\mathfrak s\otimes 1^{\otimes m})\left( \sum _{r+s=m-1} (1^{\otimes r}\otimes 
\chi\otimes 
1^{\otimes s}\big)\theta _{m}(\si c)\right).
\end{align*}
Finally
\begin{align*}
\sum _{r+s =m}(\tTheta _{r}\otimes \theta _{s})\nu (\si c)=&\sum _{r+s 
=m}(\tTheta _{r}\otimes \theta _{s})(1\otimes \Om \pi)\tpsi (\si c)\\
=&\sum _{r+s =m}(\tTheta _{r}\otimes \theta _{s})(1\otimes \Om \pi)\psi 
(\si c)\\
=&\sum _{r+s =m}(\tTheta _{r}\otimes \theta _{s})\psi (\si c)\\
=&\sum _{r+s =m}\big((\mathfrak s\otimes 1^{\otimes r})\theta _{r+1}\otimes 
\theta _{s}\big)\psi (\si c)\\
=&(\mathfrak s\otimes 1^{\otimes m})\left(\sum _{r+s =m}\big(\theta _{r+1}\otimes 
\theta _{s}\big)\psi (\si c)\right).
\end{align*}
The sum of the terms above, with appropriate signs, yields $\mathfrak s\otimes 1 ^{\otimes 
m}$ applied to the difference of the two sides of equation (1.7) for 
$\theta _{m+1}$.  The sum is therefore zero, as desired, i.e., 
the hypothesis of Proposition \ref{prop:rel-free-ext} holds for $\si c\in \si (\cbc)_{n}$.

Before showing that the same condition holds for the remaining 
generators, we consider the relation between the degree $-1$ map $\kappa:\Om 
C\to \Om (\cbc)$ of Lemma \ref{lem:kappa} and $\tTheta _{k}$.  Note that
$\tTheta _{k}$ has been defined precisely so that 
$$\tTheta _{k}\circ\kappa=(\mathfrak h\otimes 1^{\otimes k})\theta _{k+1}:\si 
C\to (\Om H\otimes _{t_{\Om}}H)\otimes H^{\otimes k}$$
for all $k$.
We show now by induction on $k$ and on wordlength in $\Om C$ that this equality holds in fact on all of $\Om C$, for 
all $k$.

Suppose that  $\tTheta _{k}\kappa=(\mathfrak h\otimes 1^{\otimes k})\theta 
_{k+1}$ everywhere in $\Om C$ for all $k<m$ and that $\tTheta 
_{m}\kappa=(\mathfrak h\otimes 1^{\otimes m})\theta 
_{m+1}$ on $T^{<\ell}\si C_{+}$.  Then on $\bigoplus _{a+b=\ell}T^a\si 
C_{+}\otimes T^b\si C_{+}$ we have that
\begin{align*}
\tTheta _{m}\kappa \mu&=\tTheta _{m}\mu (\kappa\otimes\iota 
+\iota\otimes \kappa)\\ 
&=\mu \left(\sum _{\substack{ 1\leq k\leq m\\ \mathbf j\in J_{k,m}} }(\nu 
^{\mathbf j}\tTheta _{k}\kappa\otimes\tTheta _{\mathbf j}\nu ^{(k)}\iota 
+\nu 
^{\mathbf j}\tTheta _{k}\iota\otimes\tTheta _{\mathbf j}\nu 
^{(k)}\kappa)\right)\\
&=\mu \big(\sum _{\substack{ 1\leq k\leq m\\ \mathbf j\in J_{k,m}}} (\nu 
^{\mathbf j}(\mathfrak h\otimes 1^{\otimes k})\theta _{k+1}\otimes(\mathfrak s\otimes 
1^{\otimes m})\theta _{\mathbf j}\psi ^{(k)}\\
&\qquad\qquad\qquad+(\mathfrak s\otimes 1^{\otimes m})\chi 
^{\mathbf j}\theta _{k+1}\otimes\tTheta _{\mathbf j}(\kappa\otimes 
1^{\otimes k})\psi 
^{(k)})\big)\\
&=\mu \big(\sum _{\substack{ 1\leq k\leq m\\ \mathbf j\in J_{k,m}} }((\mathfrak 
h\otimes 1^{\otimes m})\chi ^{\mathbf j}\theta _{k+1}\otimes(\mathfrak s\otimes 
1^{\otimes m})\theta _{\mathbf j}\psi ^{(k)}\\
&\qquad\qquad\qquad+(\mathfrak s\otimes 1^{\otimes m})\chi 
^{\mathbf j}\theta _{k+1}\otimes(\mathfrak h\otimes 
1^{\otimes m})\theta _{\mathbf j}\psi 
^{(k)})\big)\\
&=(\mathfrak h\otimes 
1^{\otimes m})\mu \left(\sum _{\substack{ 1\leq k\leq m\\ \mathbf j\in J_{k,m}}} 
(\chi ^{\mathbf j}\theta _{k+1}\otimes\theta _{\mathbf j}\psi ^{(k)}+\chi 
^{\mathbf j}\theta _{k+1}\otimes\theta _{\mathbf j}\psi 
^{(k)})\right)\\
&=(\mathfrak h\otimes 1^{\otimes m})\theta _{m+1}\mu.
\end{align*}
Thus $\tTheta _{m}\kappa=(\mathfrak h\otimes 1^{\otimes m})\theta 
_{m+1}$ on $T^{\leq\ell}\si C_{+}$.  

We can therefore conclude by 
induction that 
$$\tTheta _{m}\kappa=(\mathfrak h\otimes 1^{\otimes m})\theta 
_{m+1}:\Om C\to (\Om H\otimes _{t_{\Om}}H)\otimes H^{\otimes 
m}$$
for all $m$.

Suppose now that $\si \bc \in \si (\cbc)_{n}$. Let $c_{i}\otimes c^i$ denote the image of $c$ under the
reduced comultiplication in $C$.  Then
\begin{align*}
\tTheta _{m}d_{\Om}(\si \bc)=&\tTheta _{m}\big(\si 
(\overline {dc})-(-1)^{c_{i}}\si \bc_{i}\si c^i+\si c_{i}\si \bc^i\big)\\
=&(\mathfrak h\otimes 1^{\otimes m}) \theta _{m+1}(\si (dc))\\
&-(-1)^{c_{i}}\mu \left (\sum _{\substack{ 0\leq k\leq m\\\mathbf j\in 
J_{k,m}}} \hn ^{\mathbf j}\tTheta _{k}(\si \bc_{i})\otimes \hth _{\mathbf j}\nu 
^{(k)}(\si c^i)\right )\\
&+\mu \left (\sum _{\substack{ 0\leq k\leq m\\\mathbf j\in 
J_{k,m}}} \hn ^{\mathbf j}\tTheta _{k}(\si c_{i})\otimes \hth _{\mathbf j}\nu 
^{(k)}(\si \bc^i)\right )\\
=&(\mathfrak h\otimes 1^{\otimes m}) \theta _{m+1}(\si (dc))\\
&-(-1)^{c_{i}}\mu \left (\sum _{\substack{ 0\leq k\leq m\\\mathbf j\in 
J_{k,m}}} \hn ^{\mathbf j}(\mathfrak h\otimes 1^{\otimes k})\theta _{k}(\si 
c_{i})\otimes \hth _{\mathbf j}\psi 
^{(k)}(\si c^i)\right )\\
&+\mu \left (\sum _{\substack{ 0\leq k\leq m\\\mathbf j\in 
J_{k,m}}} \hn ^{\mathbf j}(\mathfrak s\otimes 1^{\otimes k})\theta _{k}(\si 
c_{i})\otimes \hth _{\mathbf j}(\kappa\otimes 1^{\otimes k})\psi 
^{(k)}(\si c^i)\right )\\
=&(\mathfrak h\otimes 1^{\otimes m}) \theta _{m+1}(\si (dc))\\
&-(-1)^{c_{i}}\mu \left (\sum _{\substack{ 0\leq k\leq m\\\mathbf j\in 
J_{k,m}}} (\mathfrak h\otimes 1^{\otimes m})\chi ^{\mathbf j}\theta _{k}(\si 
c_{i})\otimes (\mathfrak s\otimes 1^{\otimes m})\theta _{\mathbf j}\psi 
^{(k)}(\si c^i)\right )\\
&+\mu \left (\sum _{\substack{ 0\leq k\leq m\\\mathbf j\in 
J_{k,m}}} (\mathfrak s\otimes 1^{\otimes m})\chi ^{\mathbf j}\theta _{k}(\si 
c_{i})\otimes (\mathfrak h\otimes 1^{\otimes m})\theta _{\mathbf j}\psi 
^{(k)}(\si c^i)\right )\\
=&(\mathfrak h\otimes 1^{\otimes m}) \theta _{m+1}(\si (dc))\\
&+(\mathfrak h\otimes 1^{\otimes m})\mu \left (\sum _{\substack{ 0\leq k\leq m\\\mathbf j\in 
J_{k,m}}} \chi ^{\mathbf j}\theta _{k}(\si 
c_{i})\otimes \theta _{\mathbf j}\psi 
^{(k)}(\si c^i)\right )\\
=&(\mathfrak h\otimes 1^{\otimes m}) \theta _{m+1}(\si (dc))+(\mathfrak h\otimes 1^{\otimes m})\theta _{m+1}(\si c_{i}\si c^i)\\
=&(\mathfrak h\otimes 1^{\otimes m}) \theta _{m+1}d_{\Om}(\si c).
\end{align*}

Next,
\begin{align*}
(D_{\Om}\otimes 1^{\otimes m})\tTheta _{m}(\si \bc)=&(D_{\Om}\mathfrak h\otimes 1^{\otimes 
m})\theta _{m+1}(\si c)\\
=&-(\mathfrak h d\otimes 1^{\otimes m})\theta _{m+1}(\si c)\\
=&-(\mathfrak h \otimes 1^{\otimes m})(d\otimes 1^{\otimes m})\theta 
_{m+1}(\si c)\end{align*}
and
\begin{align*}
\big (1\otimes \sum _{r+s=m-1} 1^{\otimes r}\otimes d\otimes 
&1^{\otimes s}\big)\tTheta _{m}(\si \bc)\\
&=\big (\mathfrak h\otimes \sum _{r+s=m-1} 1^{\otimes r}\otimes d\otimes 
1^{\otimes s}\big)\theta _{m+1}(\si c)\\
&=(\mathfrak h\otimes 1^{\otimes m})\left(\sum _{r+s=m-1} (1^{\otimes r+1}\otimes d\otimes 
1^{\otimes s}\big)\theta _{m+1}(\si c)\right).
\end{align*}
On the other hand,
\begin{align*}
(\nu \otimes 1^{\otimes m-1})\tTheta _{m-1}(\si \bc)&=(\nu\mathfrak h \otimes 
1^{\otimes m-1})\theta _{m}(\si c)\\
&=(\mathfrak h\otimes 1^{\otimes m})(\chi\otimes 1^{\otimes m-1})\theta 
_{m}(\si c)
\end{align*}
and 
\begin{align*}
\big (1\otimes &\sum _{r+s=m-2} 1^{\otimes r}\otimes \chi\otimes 
1^{\otimes s}\big)\tTheta _{m-1}(\si \bc)\\
=&\big (\mathfrak h\otimes \sum _{r+s=m-2} 1^{\otimes r}\otimes \chi\otimes 
1^{\otimes s}\big)\theta _{m}(\si c)\\
=&(\mathfrak h\otimes 1^{\otimes m})\left( \sum _{r+s=m-1} (1^{\otimes r}\otimes 
\chi\otimes 
1^{\otimes s}\big)\theta _{m}(\si c)\right).
\end{align*}
Finally
\begin{align*}
\sum _{r+s =m}(\tTheta _{r}\otimes \theta _{s})\nu (\si \bc)=&\sum _{r+s 
=m}(\tTheta _{r}\otimes \theta _{s})\nu\kappa (\si c)\\
=&\sum _{r+s =m}(\tTheta _{r}\otimes \theta _{s})(\kappa\otimes 1)\psi 
(\si c)\\
=&\sum _{r+s =m}\big((\mathfrak h\otimes 1^{\otimes r})\theta _{r+1}\otimes 
\theta _{s}\big)\psi (\si c)\\
=&(\mathfrak h\otimes 1^{\otimes m})\left(\sum _{r+s =m}\big(\theta _{r+1}\otimes 
\theta _{s}\big)\psi (\si c)\right).
\end{align*}
The sum of the above terms, with appropriate signs, yields $\mathfrak h\otimes 1 ^{\otimes 
m}$ applied to the difference of the two sides of equation (1.7) for 
$\theta _{m+1}$.  The sum is therefore zero, as desired, i.e., 
the hypothesis of Proposition \ref{prop:rel-free-ext} holds for $\si \bc\in \si (\cbc)_{n}$.  

Since we can prove by induction that (1.7)
holds for all $m$, we can conclude that $\ttheta $ is indeed a CASH map. \hfill $\Box$

\bibliographystyle{gtart}
\bibliography{dls}

\begin{thebibliography}{}
\providecommand\bibmarginpar{\leavevmode\marginpar}
\def\urlstyle#1{{\tt #1}}

\bibitem{Ba}
\textbf{H-J Baues}, \emph{The cobar construction as a {H}opf algebra}, Invent.
  Math. 132 (1998) 467--489

\bibitem{CMN}
\textbf{F\,R Cohen}, \textbf{J\,C Moore}, \textbf{J\,A Neisendorfer},
  \emph{Torsion in homotopy groups}, Ann. of Math. (2) 109 (1979) 121--168

\bibitem{EM}
\textbf{S Eilenberg}, \textbf{J\,C Moore}, \emph{Homology and fibrations. {I}.
  {C}oalgebras, cotensor product and its derived functors}, Comment. Math.
  Helv. 40 (1966) 199--236

\bibitem{fht}
\textbf{Y F{\'e}lix}, \textbf{S Halperin}, \textbf{J-C Thomas},
  \emph{Differential graded algebras in topology}, from: ``Handbook of
  algebraic topology'', North-Holland, Amsterdam (1995)  829--865

\bibitem{GM}
\textbf{V\,K A\,M Gugenheim}, \textbf{H\,J Munkholm}, \emph{On the extended
  functoriality of {T}or and {C}otor}, J. Pure Appl. Algebra 4 (1974) 9--29

\bibitem{HPS}
\textbf{K Hess}, \textbf{P-E Parent}, \textbf{J Scott}, \emph{{Co-rings over
  operads characterize morphisms}}Preprint arXiv:math.AT/0505559

\bibitem{HPS2}
\textbf{K Hess}, \textbf{P-E Parent}, \textbf{J Scott}, \emph{A chain coalgebra
  model for the James map}, Homology, Homotopy Appl. 9 (2007) 209--231

\bibitem{HPST}
\textbf{K Hess}, \textbf{P-E Parent}, \textbf{J Scott}, \textbf{A Tonks},
  \emph{A canonical enriched {A}dams-{H}ilton model for simplicial sets}, Adv.
  Math. 207 (2006) 847--875

\bibitem{Markl}
\textbf{M Markl}, \emph{{Operads and {PROP}s}}Preprint arXiv:math.AT/0601129

\bibitem{MSS}
\textbf{M Markl}, \textbf{S Shnider}, \textbf{J Stasheff}, \emph{Operads in
  algebra, topology and physics}, volume~96 of \emph{Mathematical Surveys and
  Monographs}, American Mathematical Society, Providence, RI (2002)

\bibitem{Mi}
\textbf{R\,J Milgram}, \emph{Iterated loop spaces}, Ann. of Math. (2) 84 (1966)
  386--403

\bibitem{Mil}
\textbf{H\,R Miller}, \emph{A localization theorem in homological algebra},
  Math. Proc. Cambridge Philos. Soc. 84 (1978) 73--84

\bibitem{MM}
\textbf{J\,W Milnor}, \textbf{J\,C Moore}, \emph{On the structure of {H}opf
  algebras}, Ann. of Math. (2) 81 (1965) 211--264

\bibitem{ndombol-thomas}
\textbf{B Ndombol}, \textbf{J-C Thomas}, \emph{On the cohomology algebra of
  free loop spaces}, Topology 41 (2002) 85--106

\bibitem{HS}
\textbf{J Stasheff}, \textbf{S Halperin}, \emph{Differential algebra in its own
  rite}, from: ``Proceedings of the Advanced Study Institute on Algebraic
  Topology (Aarhus Univ., Aarhus 1970), Vol. III'', Mat. Inst., Aarhus Univ.,
  Aarhus (1970)  567--577. Various Publ. Ser., No. 13

\bibitem{S}
\textbf{R\,H Szczarba}, \emph{The homology of twisted cartesian products},
  Trans. Amer. Math. Soc. 100 (1961) 197--216

\bibitem{T}
\textbf{D Tanr{\'e}}, \emph{Homotopie rationnelle: mod\`eles de {C}hen,
  {Q}uillen, {S}ullivan}, volume 1025 of \emph{Lecture Notes in Mathematics},
  Springer-Verlag, Berlin (1983)

\end{thebibliography}

\end{document}